  \documentstyle{article}  

\newcommand{\ra}{\rightarrow} 
\newcommand{\ca}{\leadsto}

\newtheorem{prop}{Proposition}[section]
\newtheorem{th}[prop]{Theorem}
\newtheorem{cor}[prop]{Corollary}

\newtheorem{rem}[prop]{Remark} 
\newtheorem{rems}[prop]{Remarks}
\newtheorem{df}[prop]{Definition} 

\newtheorem{ex}[prop]{Example}
 
\newtheorem{op}[prop]{Open problem}

\begin{document}  

\title{New generalizations of BCI, BCK and Hilbert algebras }
\author {Afrodita  Iorgulescu  \\
\footnotesize Department of Computer Science and Cybernetics  \\ 
 \footnotesize Academy of Economic Studies\\
 \footnotesize               E-mail:  afrodita.iorgulescu@ase.ro}
\date{December 15, 2013}
\maketitle

\begin{abstract}  
 We introduce more generalizations of BCI, BCK  and  of Hilbert algebras, with proper examples, 
 and  show the  hierarchies existing between  all these
algebras, old and new ones. Namely, we found  thirty one new generalizations of BCI and BCK algebras 
and twenty generalizations of Hilbert algebras.


{\bf AMS classification (2000):} 06F15, 06F35, 06D35
 \end{abstract}

Hilbert algebras are particular cases of BCK algebras \cite{book}, while BCK algberas
are particular cases of BCI algebras.

{\it Hilbert algebras} were introduced in 1950,  in a dual form, by Henkin \cite{Henkin}, 
under the name ``implicative model", as a model of positive implicative 
propositional calculus - 
an important  fragment of classical propositional calculus introduced by Hilbert \cite{Hilbert},
\cite{Hilbert-Bernays}.   Cf. A. Diego \cite{Diego}, it was Antonio Monteiro who
has given the name  ``Hilbert algebras"
 to the dual algebras of Henkin's implicative models.

 {\it BCK algebras} and {\it BCI algebras} were introduced in 1966 by  K.  Is\'{e}ki \cite{Iseki}, as algebraic models of 
BCK-logic and of BCI-logic, respectively. The axioms of the propositional calculus of the BCK-logic are the following:\\
(B) $(\psi \ra \chi) \ra ((\varphi \ra \psi) \ra (\varphi \ra \chi))$\\
(C) $(\varphi \ra (\psi \ra \chi))\ra (\psi \ra (\varphi \ra \chi))$\\
(K) $\varphi \ra (\psi \ra \varphi)$.\\
The axioms of the propositional calculus of the BCI-logic are the above (B), (C) and the following (I):\\
(I) $\varphi \ra \varphi$.\\
 Hundred of papers were writen on BCK and BCI algebras, and  
the books \cite{Meng-Jun} and \cite{book} on BCK algebras and the book \cite{Huang} on BCI algebras.\\
Most of the commutative algebras of logic (such as  residuated lattices,
 Boolean algebras, MV algebras, Wajsberg algebras, BL algebras, G\"{o}del algebras, product algebras, 
 Hilbert algebras, Heyting algebras,
NM algebras, MTL algebras, IMTL algebras, R$_0$ algebras, weak-R$_0$ algebras etc.) 
 can be expressed as particular cases of BCK algebras (more precisely, of reversed left-BCK algebras) (see \cite{book}). \\
The BCK algebras and the commutative groups are particular cases of BCI algebras. 

Several generalizations of BCI and BCK algebras were introduced in time, namely:

 {\it BCH algebras} were introduced in 1983  by  Q.P. Hu and  X. Li \cite{Hu-Li-83}; 
their examples of proper BCH algebras given in 
\cite{Hu-Li-85} are in fact BCI algebras. There are many papers on BCH algebras since then,
 but the exact connection between BCH algebras  and BCI algebras was not found;
 we shall establish this connection in this paper.

  {\it  BCC algebras},
 also called {\it BIK$^+$ algebras}, were introduced in 1984 by Y. Komori \cite{Komori-83}, \cite{Komori-84} (see \cite{Zhang}).

{\it  BZ algebras}, also called  {\it weak-BCC algebras},
 were introduced  in 1995  by  X.H. Zhang and  R. Ye  \cite{Zhang-Ye}.

  {\it BE algebras} were introduced in 2006   by  H.S. Kim and  Y.H. Kim \cite{Kim-Kim}.

 {\it  Pre-BCK algebras} were introduced  in 2010 by  D. Bu\c sneag and  S. Rudeanu \cite{Busneag-Rudeanu}.
In pre-BCK algebras, the binary relation
$\leq$ is only a pre-order (i.e. reflexive and transitive). A BCK algebra is  
just a pre-BCK algebra verifying also the antisymmetry.  And as it can be noticed from \cite{Busneag-Rudeanu}, a
pre-BCK algebra is a BE algebra verifying an additional property.\\

In this paper, we introduce more generalizations of BCI, BCK  and  of Hilbert algebras, with proper examples, 
 and  show the  hierarchies existing between  all these
algebras, old and new ones. Namely, we found  31 new generalizations of BCI and BCK algebras 
and 20 generalizations of Hilbert algebras.
 With these new generalizations, we   better understand the BCI, BCK and the Hilbert algebras.

The paper is organized as follows:

In Section 1, we present a list of properties of BCK algebras which is the starting point of the research; we 
define the seven old algebras BCI, BCK, BCH, BCC, BZ, BE and pre-BCK in an unifying way by the properties of the list and
we draw the hierarchy of the seven old algebras.

In Section 2, which is the core of the paper, we present connections between the properties in the list.

In Section 3, we present new equivalent definitions of BCI and BCK algebras comming from the corresponding logics.

In Section 4, first  we  find the connection between BCH and BCI algebras and then we introduce step by step
the first nine new generalizations of BCI and BCK algebras:  the pre-BCC,
 aBE,  RM, pre-BZ, aRM, RME, pre-BCI, RML and aRML algebras. We present the connections 
 between  the sixteen old and new algebras.
 All the mentioned algebras are particular cases of RM  and RML algebras.

In Section 5, we introduce other twenty two new RM and RML algebras,
 generalizations of BCI and BCK algebras respectively, and we present 
their connections with the previous defined algebras.

In section 6, we introduce twenty generalizations of Hilbert algebras; all are particular cases of RML algebras.

In Section 7, we define the proper algebras mentioned in the paper.

In Sections 8, 9, 10 we present  examples of the proper old and new algebras defined in Sections 4, 5, 6 respectively.

In Section 11, we present final remarks. 

           \section{The list. The seven old  algebras. Hierarchy 0}

Let ${\cal A}=(A, \ra,1)$ be an algebra of type $(2,0)$ through this paper, where a binary relation $\leq$ can be defined by: for all $x,y$,
 $$ (\#) \qquad x \leq y \stackrel{def.}{\Longleftrightarrow}  x \ra y=1.$$
Equivalently, \\let  ${\cal A}=(A, \leq, \ra,1)$ be a structure where $\leq$ is a binary relation on $A$,
 $\ra$ is a binary operation on $A$ and $1 \in A$, all connected  by: $x \leq y  \Longleftrightarrow x \ra y=1$.

\subsection{ The list of properties}
Consider the following list of properties that can be satisfied  by ${\cal A}$ (in fact, the properties in the list 
are the most important properties 
 satisfied by a  BCK algebra),  where each property is  presented 
  in  two equivalent forms, determined by the corresponding two equivalent above definitions of ${\cal A}$:\\\\
(An) (Antisymmetry) $x \ra y=1=y \ra x \; \Longrightarrow \; x=y,$ \\
(An') (Antisymmetry) $ x \leq y, \; y \leq x \; \Longrightarrow \; x=y;$\\\\
(B) $ (y \ra z) \ra [(x \ra y) \ra (x \ra z)]=1,$\\
(B') $ y \ra z \leq (x \ra y) \ra (x \ra z)$, \\\\
(BB) $(y \ra z) \ra [(z \ra x) \ra (y \ra x)]=1$, \\
(BB') $y \ra z \leq  (z \ra x) \ra (y \ra x)$;\\\\
(*) $y \ra z=1 \;  \Longrightarrow \;  (x \ra y) \ra (x \ra z)=1$,\\
 (*') $y \leq z \;  \Longrightarrow \;  x \ra y \leq x \ra z$;\\\\         
(**) $y \ra z=1 \;  \Longrightarrow \; (z \ra x) \ra (y \ra x)=1$,\\
(**') $y \leq z \;  \Longrightarrow \; z \ra x \leq y \ra x$;\\\\
(C) $[x \ra (y \ra z)]\ra [y \ra (x \ra z)]=1$,\\
(C') $x \ra (y \ra z)\leq y \ra (x \ra z)$;\\\\
(D) $ y \ra [(y \ra x) \ra x]=1$,\\
(D') $ y \leq (y \ra x) \ra x$;\\\\
(Ex) (Exchange) $x \ra (y \ra z)=y \ra (x \ra z)$;\\\\
(K) $ x \ra (y \ra x)=1$,\\
(K') $ x \leq y \ra x$;\\\\
(L) (Last element) $x \ra 1=1$,\\ 
(L') (Last element) $x \leq 1$;\\\\
(M) $ 1 \ra x=x$;\\\\
(N) $1 \ra x=1 \;  \Longrightarrow \; x=1$,\\
(N') $1 \leq x \;  \Longrightarrow \; x=1$;\\\\
(Re) (Reflexivity) $x \ra x=1$ (we prefer here notation (Re) instead of original (I)),\\
(Re') (Reflexivity) $x \leq x$;\\\\
(S) $x = y \;  \Longrightarrow \; x \ra y =1$,\\  
(S') $x = y \;  \Longrightarrow \; x \leq y$;\\ \\
(Tr) (Transitivity) $x \ra y=1=y \ra z \; \Longrightarrow \; x \ra z=1$,\\
(Tr') (Transitivity) $x \leq y, \; y \leq z \; \Longrightarrow \; x \leq z$;\\\\
(U) $((y \ra x) \ra x) \ra x=y \ra x$. 

\begin{rem}\em Consider the following property:\\
(MP) (Modus Ponens)   $x=1$ and $x \ra y=1 \;  \Longrightarrow \; y=1$,\\
(MP') (Modus Ponens)   $x=1$ and $ x \leq y \;  \Longrightarrow \; y=1$.\\
Note that (MP) is in fact a reformulation of   (N), and viceversa.
\end{rem}

    \subsection{The seven old  algebras: \\BCI, BCK, BCH, BCC, BZ, BE, pre-BCK }

Recall now the following definitions:
\begin{df}\em An algebra $(A,\ra,1)$ is a:

O1.   {\it BCI algebra}  if it verifies  the axioms  
(BB), (D), (Re), (N), (An) or, equivalently, (BB), (D), (Re), (An), or, equivalently, (BB), (M), (An).

O2. {\it BCK algebra}   if it verifies   the axioms  (BB), (D),  (Re), (L), (An), or, equivalently,  (BB), (M), (L), (An).

O3.  {\it BCH algebra} if it verifies  the axioms (Re), (Ex), (An). 

O4. {\it BCC algebra}  if it verifies  the axioms (Re), (M), (L), (B), (An).

O5.  {\it BZ algebra} if it verifies  the axioms (Re), (M), (B), (An).

O6.  {\it BE algebra} if it verifies   the axioms (Re), (M), (L), (Ex).

O7. {\it pre-BCK algebra} if it verifies   the axioms (Re), (M), (L), (Ex), (*).  
\end{df} 
Note that, obviously, there are equivalent definitions as structures $(A, \leq, \ra,1)$.

It is known that the binary relation $\leq$ is an {\it order} relation 
 in BZ, BCI, BCC and BCK algebras and that  it is only a  {\it pre-order} in pre-BCK algebras.

Denote by {\bf BCI, BCK, BCH, BCC, BZ, BE, pre-BCK} the classes of BCI, BCK, BCH, BCC, BZ, BE, pre-BCK 
algebras respectively.
 
It is known  that:\\
 {\bf BCI} $\subseteq$ {\bf BCH}; \\
 {\bf BCI} + (L) =  {\bf   BCK}; $\quad$
  {\bf BCC } +  (Ex) = {\bf BCK};\\
 {\bf BZ} +  (Ex) = {\bf BCI};  $\qquad$ {\bf BZ} + (L) = {\bf  BCC};\\
   {\bf BE} + (*) = {\bf  pre-BCK};
 {\bf pre-BCK} + (An) = {\bf BCK}.

It is known  that: 
\begin{prop}\label{propBCI}
   In any BCI algebra, the following properties hold: (Ex), (U), (B), (M), (*) and (**).
\end{prop}

          \subsection{Hierarchy 0, of the seven old algebras}

We introduce now the following definition.
\begin{df}\em Let ${\cal A}=(A, \ra,1)$ be an algebra of type $(2,0)$ and
 let  $\leq$  be the associated binary relation defined by ($\#$).

1) We shall say that ${\cal A}$ is {\it reflexive} if $\leq$ is reflexive (i.e. it satisfies  property (Re) or (Re').
                                
2) We shall say that ${\cal A}$ is {\it antisymmetrique} if $\leq$ is antisymmetrique  (i.e.
 it satisfies  property (An) or (An')).

3) We shall say that ${\cal A}$ is {\it transitive} if $\leq$ is transitive (i.e. it satisfies  property (Tr) or (Tr')).

4) We shall say that ${\cal A}$ is {\it pre-ordered} if $\leq$ is a pre-order relation
 (i.e. it is reflexive and transitive).

5) We shall say that ${\cal A}$ is {\it ordered} if $\leq$ is a partial-order relation 
(i.e. it is reflexive, antisymmetrique and transitive).

6) We shall say that ${\cal A}$ is a {\it lattice} if $\leq$ is a lattice-order relation 
(i.e. it is a partial-order such that there exists sup(x,y) and inf(x,y) for each $x,y \in A$); we shall use the notation
$x \vee y$ for sup(x,y) and $x \wedge y$ for inf(x,y), 
with $x \leq y \Leftrightarrow x \vee y=y  \Leftrightarrow x \wedge y=x$.
\end{df}

\begin{rem}\em
In a hierarchy of classes of algebras, we shall represent:\\
- a  class of  {\it reflexive} algebras    by $\bigcirc$\\
- a  class of {\it antisymmetrique} algebras   by $\circ$\\
- a  class of {\it transitive} algebras    by $\bullet$\\
Consequently, we shall represent:\\
- a class of algebras where $\leq$ is {\it  reflexive and antisymetrique}  
and\\
- a class of  {\it pre-ordered}   algebras (i.e. $\leq$ is reflexive and transitive), respectively, by: \\

\begin{figure}[htbp]
\begin{center}
\begin{picture}(100,-180)(0,0) 
  \put(10,0){\makebox(2,2){$\bigcirc$}} \put(10,0){\makebox(2,2){$\circ$}} 
\put(40,0){\makebox(2,2){$\bigcirc$}} \put(40,0){\makebox(2,2){$\bullet$}} 
\end{picture}
 \end{center}
\end{figure}

${}$\\
-  a  class of {\it ordered} algebras  by: 

${} $ $ \qquad \qquad \qquad \qquad \; \quad \qquad \qquad \qquad \qquad \; \quad $   \circle*{11} 
\end{rem}

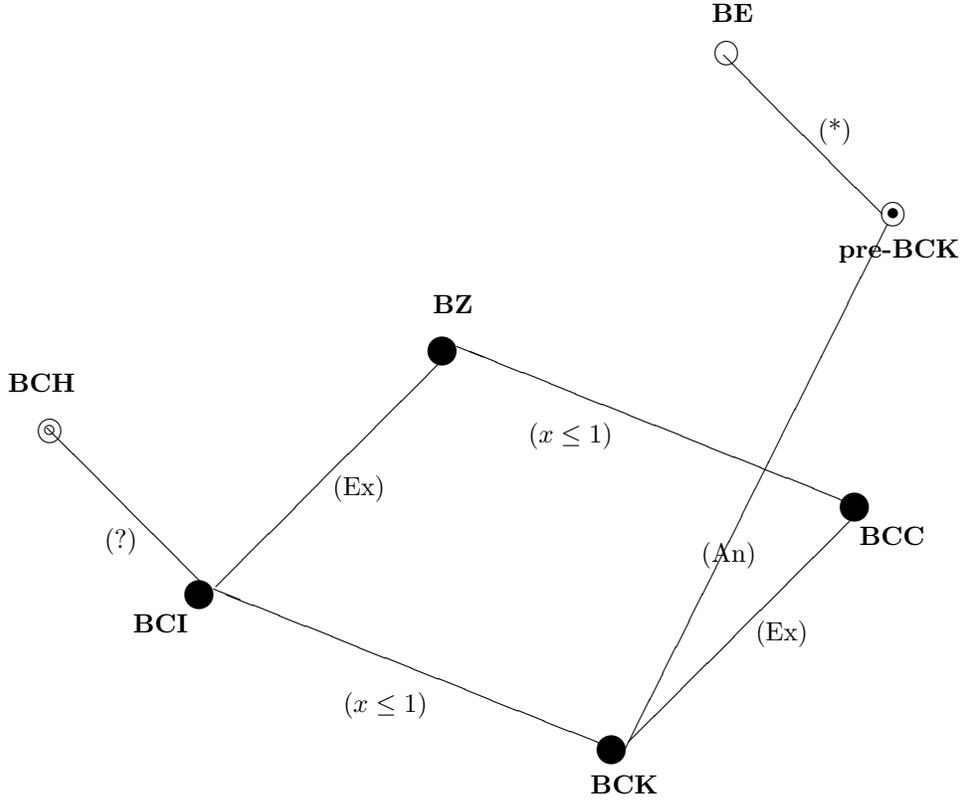
\begin{figure}[htbp]
\begin{center}
\begin{picture}(420,260)(-180,25) 
  \put(137,119){\line(-1,-1){90}}    
     \put(45,25){\line(1,2){100}} 
\put(83,289){\line(1,-1){60}}   

\put(83,289){\makebox(2,2){$\bigcirc$}}  \put(83,300){\makebox(10,10){{\bf BE} }}
\put(146,228){\makebox(2,2){$\bigcirc$}} \put(146,228){\makebox(2,2){$\bullet$}}\put(146,210){\makebox(10,10){{\bf pre-BCK} }}
\put(80,95){\makebox(10,10){(An)}}
\put(120,255){\makebox(10,10){(*)}}
\put(45,25){\makebox(2,2){\circle*{11}}}
 
\put(100,65){\makebox(10,10){(Ex)}}

\put(32,8){\makebox(30,10){{\bf BCK} }}

\put(137,117){\makebox(2,2){\circle*{11}}} 
                                          
\put(130,102){\makebox(30,10){  {\bf BCC}}}

   \put(-19,178){\line(-1,-1){90}}    
     
       \put(-173,148){\line(1,-1){60}}    

  \put(45,25){\line(-5,2){155}}
 
  \put(137,117){\line(-5,2){155}}

  
\put(-111,84){\makebox(2,2){\circle*{11}}}

\put(-150,69){ \makebox(30,10){ {\bf BCI}}}
\put(-173,146){\makebox(2,2){$\bigcirc$}} \put(-173,146){\makebox(2,2){$\circ$}}

 \put(-60,120){\makebox(10,10){(Ex)}}

                           \put(-150,100){\makebox(10,10){(?)}}

\put(-190,160){\makebox(30,10){{\bf BCH}}}

\put(10,140){\makebox(30,10){ ($x \leq 1$) }}
\put(-60,38){\makebox(30,10){($x \leq 1$)}}

\put(-19,176){\makebox(2,2){\circle*{11}}}
 
\put(-36,190){\makebox(30,10){ {\bf BZ}}}                                         

\end{picture}
 \end{center}
\caption{ {\bf Hierarchy 0} of the seven old algebras   }
\label{fig:fig0}
\end{figure}

Hence, we can draw the following hierarchy, called   {\it Hierarchy 0 (zero)}, connecting the seven old algebras,  
  as shown in Figure \ref{fig:fig0}, where ``?" means ``unknown".

         \section{Connections between the properties in the list}

             \subsection{Connections}

We shall establish several connections   between the properties in the above list. These connections 
are the core of the paper. Based on them, we then introduced the new generalizations claimed in the title of the paper.

\begin{prop}\label{propp} Let $(A,\ra,1)$ be an algebra of type $(2,0)$. Then the following are true:\\
 
(0.) (Re) implies (S);

(00.) (M) implies (N);\\

(1.)  (L) + (An)  imply (N);

(2.) (K) + (An)  imply (N);\\

(3.)  (C) + (An)  imply (Ex); $\quad$ (3'.) (Ex) + (Re) imply (C);\\

(4.) (Re) + (Ex) imply  (D);

(5.) (Re) + (Ex) + (An) imply (M);  $\quad$ (5'.) (Re) + (Ex) + (An) imply (N);\\

(6.) (Re) + (K) imply (L);

 (7.)    (N) + (K) imply (L); $\quad$  (7'.)    (M) + (K) imply (L);\\

(8.) (Re) + (L) + (Ex) imply (K);


(9.) (M) + (L) + (B) imply (K);  $\quad$ (9'.) (M) + (L) + (**) imply (K); \\

(10.) (Ex) implies (B) $\Leftrightarrow$ (BB); 

(10'.) (Ex) + (B) imply (BB); $\quad$ (10".) (Ex) + (BB) imply (B);\\

(11.) (Re) + (Ex) + (*) imply (BB);\\

(12.) (N) + (B) imply (*); $\quad$ (12'.) (M) + (B) imply (*);

(13.) (N) + (*) imply (Tr); $\quad$  (13'.) (M) + (*) imply (Tr);

(14.)  (N) + (B) imply (Tr); $\quad$  (14'.) (M) + (B) imply (Tr);\\

(15.) (N) + (BB) imply (**); $\quad$ (15'.) (M) + (BB) imply (**);

(16.) (N) + (**) imply (Tr); $\quad$ (16'.) (M) + (**) imply (Tr);

(17.) (N) + (BB) imply (Tr); $\quad$ (17'.) (M) + (BB) imply (Tr);\\

(18.)  (M) + (BB) imply (Re); $\quad$ (18')  (M) + (BB) imply (D);

(19.)  (M) + (B) imply (Re);\\

(20.) (BB) + (D) + (N) imply (C); $\quad$ (20') {\bf (M) + (BB) imply (C)};

(21.) (BB) + (D) + (N) + (An) imply (Ex);

(21'.) (BB) + (D) + (L) + (An) imply (Ex); $\quad$ (21''.) {\bf (M) + (BB) + (An) imply (Ex)};\\

(22.) (B) + (C) + (K) + (An) imply (Re);

(23.) (BB) + (D) + (Re) + (An) imply (N);\\

(24.) (Re) + (Ex) + (Tr) imply (**).
\end{prop} 

{\bf Proof.}
(0.):  Suppose $x=y$; then $x \ra y=y \ra y \stackrel{(Re)}{=} 1$; thus, (S) holds.

(00.): Suppose $1 \ra x=1$. Then, by (M), we get $x=1$, i.e. (N) holds.

(1.): Suppose $1 \ra x=1$. By (L), we also have $x \ra 1=1$. Hence, by (An), we get $x=1$, i.e. (N) holds.

(2.): Suppose $1 \ra x=1$; by (K), we  have $ x \ra (1 \ra x) =1$, then $x \ra 1=1$. Hence, by (An), $x=1$, i.e. (N) holds.
  
(3.): By (C), we have: $[x \ra (y \ra z)] \ra [y \ra (x \ra z)]=1$
 and also $[y \ra (x \ra z)] \ra [x \ra (y \ra z)]=1$; hence, by (An), we obtain that: $x \ra (y \ra z)=y \ra (x \ra z)$,
i.e. (Ex) holds.


(3'.): By (Ex), we have 
\begin{equation}\label{3-1}
x \ra (y \ra z)=y \ra (x \ra z);
\end{equation}
 by (0.), (Re) implies (S); hence,  by (S), (\ref{3-1}) implies  $[x \ra (y \ra z)] \ra [ y \ra (x \ra z)]=1$, i.e. (C) holds.
  
(4.): $y \ra [(y \ra x) \ra x] \stackrel{(Ex)}{=} (y \ra x) \ra (y \ra x) \stackrel{(Re)}{=} 1$, i.e. (D) holds.

(5.): $x \ra (1 \ra x) \stackrel{(Ex)}{=} 1 \ra (x \ra x) \stackrel{(Re)}{=} 1 \ra 1 \stackrel{(Re)}{=}1$.

We shall prove now that $(1 \ra x) \ra x =1$ also. Indeed, by (5), (D) holds, consequently, $1 \ra [(1 \ra x) \ra x]=1$. Since by (4), (N) holds too,
we obtain that $(1 \ra x) \ra x =1$. \\
Applying now (An), we obtain that $1 \ra x=x$, i.e. (M) holds.


  (5'.): By above (5), (Re) + (Ex) + (An) imply (M); and by above (0), (M) implies (N); hence, (Re) + (Ex) + (An) imply (N).

(6.): In (K) ($x \ra (y \ra x)=1$), take $y=x$: we get $1=x \ra (x \ra x) \stackrel{(Re)}{=}x \ra 1$, i.e. (L) holds. 


(7.): By (K), we have: $1 \ra (x \ra 1)=1$; hence, by (N), we obtain that $x \ra 1=1$, i.e. (L) holds.


(7'.): By  above (00.), (M) implies (N); and by (7), (N) + (K) imply (L); thus, (L) holds.


(8.): In (Ex) ($x \ra (y \ra z)=y \ra (x \ra z)$), take  $z=x$: we get $x \ra (y \ra x)=y \ra (x \ra x) \stackrel{(Re)}{=} y \ra 1 \stackrel{(L)}{=}1$,
i.e. (K) holds.


(9.): Take $y=1$ in (B) ($(y \ra z) \ra [(x \ra y) \ra (x \ra z)]=1$); we obtain:
 $(1 \ra z) \ra [(x \ra 1) \ra (x \ra z)]=1$; then, by (M), we obtain: $z \ra [(x \ra 1) \ra (x \ra z)]=1$; then by
(L) and (M) again, we obtain  $z \ra (x \ra z)=1$, i.e. (K) holds.


(9'.): By (L'), $y \leq 1$ is true and hence, by (**), we get: $1 \ra x \leq y \ra x$,
 which by (M) means that $x \leq y \ra x$, i.e. (K') holds.

(10.): $ (y \ra z) \ra [(x \ra y) \ra (x \ra z)]  \stackrel{(Ex)}{=}
 (x \ra y) \ra [(y \ra z) \ra (x \ra z)]$. Hence, (B) $\Leftrightarrow$ (BB).


(10'.) Obviously, by (10.).


(10".) Obviously, by (10.).


(11.): Since (Re) + (Ex) imply (D), i.e $y \ra [(y \ra z) \ra z]=1$, we apply (*) and we obtain: 
$$(x \ra y) \ra (x \ra [(y \ra z) \ra z])=1.$$
Then, by (Ex), we obtain:
$$(x \ra y) \ra [(y \ra z) \ra (x \ra z)]=1,$$
i.e. (BB) holds. 

(12.): Suppose $y \ra z=1$; then, by (B) ($ (y \ra z) \ra [(x \ra y) \ra (x \ra z)]=1$), 
it follows that $ 1 \ra [(x \ra y) \ra (x \ra z)]=1$;
hence, by (N), we obtain that $(x \ra y) \ra (x \ra z)=1$, i.e. (*) holds.


(12'.): By (00.), (M) implies (N); then apply above (12).

(13.): Suppose $x \ra y =1= y \ra z$, or $y \ra z=1=x \ra y$;   then, by (*), we obtain: $1 \ra (x \ra z)=1$;
 hence by (N), we obtain: $x \ra z=1$, i.e. (Tr) holds.


(13'.): By (00.), (M) implies (N); then apply (13).


(14.): Suppose $y \ra z=1$ and $x \ra y=1$; then, by (B), we obtain: $1 \ra [1 \ra (x \ra z)]=1$; then, by (N), we obtain:
$1 \ra (x \ra z)=1$; by (N) again, we obtain: $x \ra z=1$. Thus, (Tr) holds.


(14'.): By (00.), (M) implies (N); then, apply (14).

(15.): Suppose $y \ra z=1$; then, by (BB) ($(y \ra z) \ra [(z \ra x) \ra (y \ra x)]=1$, 
it follows that $1 \ra [(z \ra x) \ra (y \ra x)]=1$; hence,
by (N), we obtain $(z \ra x) \ra (y \ra x)=1$, i.e. (**) holds.


(15'.): By (00.), (M) implies (N); then apply (15).


(16.): Suppose $y \ra z =1=z \ra x$;   then, by (**), we obtain: $1 \ra (y \ra x)=1$; hence,
 by (N), we obtain: $y \ra x=1$, i.e. (Tr) holds.

(16'.): By (00.), (M) implies (N); then apply (16).

 
(17.):  Suppose $y \ra z =1$ and $z \ra x=1$; then, by (BB),  we obtain: $1 \ra [1 \ra (y \ra x)]=1$; then, 
applying (N) twice, we obtain $y \ra x=1$. Thus, (Tr) holds.


(17'.): By (00.), (M) implies (N); then apply (17).


(18.): In (BB) ($(y \ra z) \ra [(z \ra x) \ra (y \ra x)]=1$), take $y=z=1$; we obtain: $(1 \ra 1) \ra [(1 \ra x) \ra (1 \ra x)]=1$, hence, by (M),
$1 \ra [x \ra x]=1$, hence by (M) again, $x \ra x=1$, i.e. (Re) holds.\\


(18'.): In (BB) ($(y \ra z) \ra [(z \ra x) \ra (y \ra x)]=1$), take $y=1$; we obtain: $(1 \ra z) \ra [(z \ra x) \ra (1 \ra x)]=1$, i.e.
$z \ra [(z \ra x) \ra x]=1$), by (M); thus, (D) holds.\\


(19.)  In (B) ($ (y \ra z) \ra [(x \ra y) \ra (x \ra z)]=1$), take $x=y=1$; we obtain: $ (1 \ra z) \ra [(1 \ra 1) \ra (1 \ra z)]=1$, hence, by (M),
 $ z \ra [1 \ra  z]=1$,  by (M) again, $z \ra z=1$; thus, (Re) holds.\\   


(20.): (see \cite{Iseki-Tanaka}, Theorem 1) 
We shall use $\leq$ for a better understanding.\\
 By (BB'), we have: $y \ra z \leq (z \ra x) \ra (y \ra x)$. By (15.), (BB) + (N) imply (**).  Then, by (**'), we obtain:
\begin{equation}\label{18}
[(z \ra x) \ra (y \ra x)] \ra u \leq (y \ra z) \ra u.
\end{equation}

We substitute in (\ref{18}): $x$ by  $u \ra x$, $z$ by $z \ra x$; $u$ by $(u \ra z) \ra (y \ra (u \ra x))$. Then, we obtain:

$$V \stackrel{notation}{=}[((z \ra x) \ra (u \ra x)) \ra (y \ra (u \ra x))] \ra 
[(u \ra z) \ra (y \ra (u \ra x))] \leq $$
$$\leq (y \ra (z \ra x)) \ra [(u \ra z) \ra (y \ra (u \ra x))] \stackrel{notation}{=} W.$$
Then, the left side $V=1$, by (\ref{18}) (with $Y=u$, $U=y \ra (u \ra x)$). Thus, $1 \leq W$; then, by (N'), $W=1$, i.e.

\begin{equation}\label{18-1} 
y \ra (z \ra x) \leq (u \ra z) \ra (y \ra (u \ra x)).
\end{equation}

Take now $u=z$, $z=y \ra x$ in (\ref{18-1}). Then, we obtain:
\begin{equation}\label{18-2}
y \ra ((y \ra x) \ra x) \leq [z \ra (y \ra x)] \ra [y \ra (z \ra x)].
\end{equation}

By (D), from (\ref{18-2}) we obtain: 
\begin{equation}\label{18-3}
1 \leq [z \ra (y \ra x)] \ra [y \ra (z \ra x)].
\end{equation}

From (\ref{18-2}), we obtain, by (N') again, that:
$[z \ra (y \ra x)] \ra [y \ra (z \ra x)]=1,$  i.e.  (C) holds.


(20'.): (M) implies (N), by above (00.), and (M) + (BB) imply (D), by above (18'.). Hence, (M) + (BB) imply 
(BB) + (D) + (N), which imply (C), by above (20.) 


(21.): By above (20.), (BB) + (D) + (N) imply (C); and (C) implies that:

\begin{equation}\label{18-4} 
[z \ra (y \ra x)] \ra [y \ra (z \ra x)]=1.
\end{equation}
 and also  that:
\begin{equation}\label{18-5}
[y \ra (z \ra x)] \ra [z \ra (y \ra x)] =1.
\end{equation}
From (\ref{18-4}) and (\ref{18-5}), by applying (An), we obtain that $y \ra (z \ra x) = z \ra (y \ra x)$, i.e. (Ex) holds.


(21'.): By  above (1.), (L) + (An) imply (N); then apply (21.).


(21''.): By above (20'.), (M) + (BB) imply (C); by above (3.), (C) + (An) imply (Ex); thus, (M) + (BB) + (An) imply (Ex).


(22.): Take $z=1$ and $y=x$ in (B) ($ (y \ra z) \ra [(x \ra y) \ra (x \ra z)]=1$); we obtain:
\begin{equation}\label{19}
 (x \ra 1) \ra [(x \ra x) \ra (x \ra 1)]=1.
\end{equation}
By above (2.),  (K) + (An)  imply (N); by above (7.),     (N) + (K) imply (L); by above (3.),   (C) + (An)  imply (Ex).
Then, from (\ref{19}), we obtain, by (L):
\begin{equation}\label{19-1}
 1 \ra [(x \ra x) \ra 1]=1.
\end{equation}
From (\ref{19-1}), by (N), we obtain: 
\begin{equation}\label{19-2}
 (x \ra x) \ra 1=1.
\end{equation}
On the other hand, $1 \ra (x \ra x) \stackrel{(Ex)}{=} x \ra (1 \ra x) \stackrel{(K)}{=} 1$, hence
\begin{equation}\label{19-3}
 1 \ra (x \ra x)=1.
\end{equation}
From (\ref{19-2}) and  (\ref{19-3}), by applying (An), we obtain $x \ra x=1$, i.e. (Re) holds.


(23.): Suppose that:
 \begin{equation}\label{23-1}
1 \ra x=1;
\end{equation} 
we must prove that $x=1$. Indeed, in  (BB) ($(y \ra z) \ra [(z \ra x) \ra (y \ra x)]=1$), take $y=x$ and $z=1$; we obtain:
$(x \ra 1) \ra [(1 \ra x) \ra (x \ra x)]=1$, hence,  by (\ref{23-1}) and (Re), $(x \ra 1) \ra [1 \ra 1]=1$, hence, by (Re) again, 
\begin{equation}\label{23-2}
(x \ra 1) \ra 1=1.
\end{equation}
Since by (D) we have: $x \ra [(x \ra 1) \ra 1] =1$, it follows, by  (\ref{23-2}), that: 
\begin{equation}\label{23-3}
x \ra 1=1.
\end{equation}
From  (\ref{23-1}), (\ref{23-3})  and (An), we obtain $x=1$.  Thus, (N) holds.

(24.): Suppose (Tr) holds, i.e. $X \ra Y=1, \; Y \ra Z=1$ imply $X \ra Z=1$.\\
We must prove that (**) holds, i.e. $y \ra z=1$ implies $(z \ra x) \ra (y \ra x)=1$.\\
Suppose that  $y \ra z=1$; we must prove that $H \stackrel{notation}{=}(z \ra x) \ra (y \ra x)=1$.\\
 Indeed, $H \stackrel{(Ex)}{=}y \ra ((z \ra x) \ra x)$. 
Take  $X=y$, $Z=(z \ra x) \ra x$ and $Y=z$; then we have $H=X \ra Z$ and:\\
$X \ra Y=y \ra z=1$, by hypothesis;\\
$Y \ra Z=z \ra [(z \ra x) \ra x] \stackrel{(Ex)}{=} (z \ra x) \ra (z \ra x) \stackrel{(Re)}{=}1$.\\
 Hence, by (Tr), it follows that $X \ra Z=1$, i.e. $H=1$.
\hfill $\Box$\\

Now we prove two important results. 
\begin{th}\label{th1} (Generalization of (\cite{Busneag-Rudeanu}, Lemma 1.2 and Proposition 1.3))

 If an  algebra $(A, \ra, 1)$ verifies properties (Re), (M), (Ex), then:
$$ (B) \; \Leftrightarrow \; (BB) \; \Leftrightarrow \; (*).$$
\end{th}
{\bf Proof.} \\
By Proposition \ref{propp} (9.), (Ex) implies that (B) $\Leftrightarrow $ (BB).\\
By Proposition \ref{propp} (11.), (M) + (B) implies  (*).\\
By Proposition \ref{propp} (10.), (Re) + (Ex) + (*) implies (BB). Hence, we have:
$$ (*) \Rightarrow (BB) \; \Leftrightarrow \; (B) \; \Rightarrow (*),$$
thus $ (B) \; \Leftrightarrow \; (BB) \; \Leftrightarrow \; (*).$
\hfill $\Box$

\begin{th}\label{th2}${}$

 If an algebra  $(A,\ra,1)$ verifies properties (Re), (M), (Ex), then:
$$  (**) \;  \Leftrightarrow \; (Tr).$$
\end{th}
{\bf Proof.}\\ 
By Proposition \ref{propp} (16'), (M) + (**) imply (Tr).\\
By Proposition \ref{propp} (24),  (Re) + (Ex) + (Tr) imply (**). 
\hfill$\Box$

\begin{th}\label{th3} ${}$

If an algebra  $(A,\ra,1)$ verifies properties (Re), (M), (B), (An), then:
$$  (Ex) \;  \Leftrightarrow \; (BB).$$
\end{th}
{\bf Proof.} By Proposition \ref{propp} (10.), (Ex)
implies (B) $\Leftrightarrow$ (BB); hence, (Ex) implies (BB). \\
By Proposition \ref{propp} (21''.), (M) + (BB) + (An)  imply (Ex), hence (BB) implies (Ex).
\hfill$\Box$

It follows immediately that:
\begin{cor}\label{corA} In BZ, BCC algebras  we have $(Ex) \;  \Leftrightarrow \; (BB)$.
\end{cor}

The next theorem was proved by professor Michael Kinyon, Department of Mathematics, University of Denver,  following
our open problems announced in the preprint on arXiv; he proved (i)
first by using the automated theorem proving tool Prover9.
\begin{th}\label{th4}  (Michael Kinyon)
In any algebra $(A, \ra,1)$ we have:\\
(i) (M) + (BB) imply (B),\\
(ii) (M) + (B) imply (**).
\end{th}
{\bf Proof.} (i): By Proposition \ref{propp} (18'.), we have  (M) + (BB) imply (D).\\ Next, if (BB) is
$$ (x \ra y) \ra [(y \ra z) \ra (x \ra z)]=1,$$
in (BB) set $x=u$ and $y=(u \ra v) \ra v$, to get:\\
$(u \ra [(u \ra v) \ra v]) \ra [(((u \ra v) \ra v) \ra z) \ra (u \ra z)] \stackrel{(D)}{=}$\\
$ 1 \ra [(((u \ra v) \ra v) \ra z) \ra (u \ra z)] \stackrel{(M)}{=}$ \\
$(((u \ra v) \ra v) \ra z) \ra (u \ra z)=1$. \\After renaming variables, we get:
$$ (a) \qquad  (((x \ra y) \ra y) \ra z) \ra (x \ra z)=1.$$
Next, in (BB) set $x=u \ra v$ and $y=(v \ra w) \ra (u \ra w)$,  to get:\\
$((u \ra v) \ra [(v \ra w) \ra (u \ra w)]) \ra [(((v \ra w) \ra (u \ra w)) \ra z) \ra ((u \ra v) \ra z)] \stackrel{(BB)}{=}$\\
$1 \ra [(((v \ra w) \ra (u \ra w)) \ra z) \ra ((u \ra v) \ra z)] \stackrel{(M)}{=}$\\
$(((v \ra w) \ra (u \ra w)) \ra z) \ra ((u \ra v) \ra z)=1$. \\
After renaming variables, we get:
$$ (b) \qquad (((x \ra y) \ra (u \ra y)) \ra z) \ra ((u \ra x) \ra z)=1.$$
Taking $z=u \ra y$ in (b), we get:
$$ (c) \qquad (((x \ra y) \ra (u \ra y)) \ra (u \ra y)) \ra ((u \ra x) \ra (u \ra y))=1.$$
Now, in (a) set $x=v \ra w$, $y=t \ra w$, $z=(t \ra v) \ra (t \ra w)$ to get:\\    
$[(((v \ra w) \ra (t \ra w)) \ra (t \ra w)) \ra ((t \ra v) \ra (t \ra w))] \ra ((v \ra w) \ra ((t \ra v) \ra (t \ra w)))
\stackrel{(c)}{=}$\\
$1 \ra ((v \ra w) \ra ((t \ra v) \ra (t \ra w))) \stackrel{(M)}{=}$\\
$(v \ra w) \ra ((t \ra v) \ra (t \ra w))=1$, i.e. (B) holds.

(ii): Suppose (B) is 
$$(y \ra z) \ra [(x \ra y) \ra (x \ra z)]=1.$$
If $x \ra y=1$,  then we get from (B): \\
$(y \ra z) \ra [1 \ra (x \ra z)] \stackrel{(M)}{=} (y \ra z) \ra (x \ra z)=1,$ i.e. (**) holds.
\hfill$\Box$

    \subsection{Generalities}
Let as above ${\cal A}= (A, \ra,1)$  verifying $(\#)$  or ${\cal A}= (A, \leq, \ra,1)$.

\begin{df} \em ${}$

(i)  If there exists an element $0 \in A$ such that $0 \leq x$ for all $x \in A$, then  we shall call it 
 {\it zero} of ${\cal A}$.
 
(ii) If there  exists a {\it zero} of ${\cal A}$ and it is unique, denote it by $0$,
and if  property (L) holds, then we shall say that ${\cal A}$ is bounded; it will be then denoted by $(A, \ra, 0,1)$.

(iii) If ${\cal A}=(A, \ra, 0,1)$ is bounded,  then we define
a new operation ${}^-$ on $A$, called {\it negation}, by: for all $x \in A$,
$$x^- \stackrel{def.}{=}x \ra 0.$$

(iv) If the negation ${}^-$ verifies the property: \\
(DN) (Double negation) for all $x \in A$,  $(x^-)^-=x$,\\
we shall say that ${\cal A}$ is {\it involutive} or {\it with (DN)}.
\end{df}

\begin{prop} If an algebra  $(A,\ra,1)$ verifies (Ex), then:
$$ (G1) \quad  x \ra y^-=y \ra x^-.$$ 
\end{prop}
{\bf Proof.} $x \ra y^-=x \ra (y \ra 0) \stackrel{(Ex)}{=} y \ra (x \ra 0) =y \ra x^-$.
\hfill $\Box$

\begin{prop} If a bounded  algebra  $(A,\ra,0,1)$ verifies (Ex) and (DN), then:
$$ (G2) \quad  x \ra y=y^- \ra x^-,$$ 
$$ (G3) \quad y^- \ra x=x^- \ra y.$$
\end{prop}
{\bf Proof.} $x \ra y \stackrel{(DN)}{=}x \ra (y^-)^- \stackrel{(G1)} y^- \ra x^-$.\\
$y^- \ra x \stackrel{(DN)}{=}y^- \ra (x^-)^- \stackrel{(G2)}{=}x^- \ra y.$
\hfill $\Box$

\begin{prop} If an algebra  $(A,\ra,1)$ verifies (D), then:
$$  (G4)  \quad x \leq (x^-)^-.$$ 
\end{prop}
{\bf Proof.} $(x^-)^-=(x \ra 0) \ra 0) \stackrel{(D)}{\geq} x $.
\hfill $\Box$

\begin{prop} If an algebra  $(A,\ra,1)$ verifies (BB), then:
$$ (G5) \quad x \ra y \leq y^- \ra x^-.$$ 
\end{prop}
{\bf Proof.} $y^- \ra x^- = (y \ra 0) \ra (x \ra 0)  \stackrel{(BB')}{\geq} x \ra y$.
\hfill $\Box$

\begin{prop} If an algebra  $(A,\ra,1)$ verifies (**), then:
$$  (G6) \quad x \leq y \Longrightarrow  y^- \leq x^-.$$ 
\end{prop}
{\bf Proof.} $x \leq y \stackrel{(**')}{\Longrightarrow} (y \ra 0) \leq (x \ra 0) $, i.e.
$x \leq y \Longrightarrow y^- \leq x^-$.
\hfill $\Box$

\begin{prop} If a bounded  algebra  $(A,\ra,0,1)$ verifies (**) and (DN), then:
$$ (G7) \quad x \leq y \Longleftrightarrow  y^- \leq x^-.$$ 
\end{prop}
{\bf Proof.} We have by (G6)  $x \leq y \Longrightarrow y^- \leq x^-$ and, similarly, 
 $y^- \leq x^- \Longrightarrow (x^-)^- \leq (y^-)^-$, i.e. $y^- \leq x^- \Longrightarrow x  \leq y$, by (DN).
Hence, $x \leq y \Longleftrightarrow y^- \leq x^-$.
\hfill $\Box$

\begin{prop} If an algebra  $(A,\ra,1)$ verifies (U), then:
$$ (G8) \quad ((x^-)^-)^-=x^-.$$ 
\end{prop}
{\bf Proof.} $((x^-)^-)^-=((x \ra 0) \ra 0) \ra 0 \stackrel{(U)}{=} x \ra 0=x^-$.
\hfill $\Box$

\section{New equivalent definitions of BCI and  of BCK  algebras, comming from logic} 

\subsection{A new equivalent definition of  BCI  algebras comming from logic}   

We shall present now a new equivalent definition of BCI algebras,  starting from the axioms (B), (C), (I)  of the BCI-logic.

\begin{prop}\label{def11} Let  ${\cal A}=(A, \ra,1)$ be an algebra 
 of type (2,0) verifying the following axioms: for all $x,y,z \in A$, \\
(B) $(y \ra z) \ra ((x \ra y) \ra (x \ra z))=1$,\\
(C) $(x \ra (y \ra z))\ra (y \ra (x \ra z))=1$,\\
(Re) $x \ra x =1$, \\
(An) $x \ra y=1=y \ra x$ imply $x=y$,\\
where the name (Re), comming from "Reflexivity", is used  instead of the initial name (I)
 for the corresponding  logical axiom of BCI logic.

Then ${\cal A}$ is a BCI-algebra.
\end{prop}
{\bf Proof.} We shall prove that axioms (BB), (D), (Re), (N), (An) from the first definition of BCI algebras  hold. Indeed, \\
- by Proposition \ref{propp} (3.), (C) + (An) imply (Ex); then, by Proposition \ref{propp} (10'.), (B) + (Ex) imply (BB); thus, (BB)  holds;\\
- by Proposition \ref {propp} (4.), (Re) + (Ex) imply (D); thus, (D) holds;\\
- by Proposition \ref {propp} (5'.), (Re) + (Ex) + (An) imply (N); thus, (N) holds.
\hfill $\Box$

\begin{prop}\label{def11-b} A BCI algebra verifies (B), (C), (Re), (An).
\end{prop}
{\bf Proof.} Suppose that  axioms (BB), (D), (Re), (N), (An) from the first definition of BCI algebras  hold. Then:\\
 - by Proposition \ref{propp} (20.), (BB) + (D) + (N)  imply (C); thus, (C) holds;\\
- by Proposition \ref{propp} (3.), (C) + (An) imply (Ex);  then, by Proposition \ref{propp} (10".), (Ex) + (BB) imply (B); thus, (B) holds.
\hfill $\Box$\\

By Propositions \ref{def11} and \ref{def11-b}, an algebra $(A, \ra,1)$ verifying the axioms (B), (C), (Re),  (An) is an 
equivalent definition,  coming from BCI logic, of a BCI algebra.
Hence, we have the following

\begin{th}
An algebra $(A,\ra,1)$ of type $(2,0)$ is a BCI algebra if and only if properties (B), (C), (Re), (An) are satisfied.
\end{th}
\subsection{A new equivalent definition of BCK  algebras comming from logic} 

We shall present now a new equivalent definition of BCK algebras,  starting from the axioms (B), (C), (K)  of the BCK-logic.

\begin{prop}\label{def10} Let  ${\cal A}=(A, \ra,1)$ be an algebra 
 of type (2,0) verifying the following axioms: for all $x,y,z \in A$, \\
(B) $(y \ra z) \ra ((x \ra y) \ra (x \ra z))=1$,\\
(C) $(x \ra (y \ra z))\ra (y \ra (x \ra z))=1$,\\
(K) $x \ra (y \ra x)=1$,\\
(An) $x \ra y=1=y \ra x$ imply $x=y$.

Then ${\cal A}$ is a BCK-algebra.
\end{prop}
{\bf Proof.} We shall prove that axioms (BB), (D), (Re), (L), (An) from the first definition of BCK algebras hold. Indeed,\\
- by Proposition \ref{propp} (2.), (K) + (An) imply (N) and (N) + (K) imply (L); thus, (L) holds.\\
- by Proposition \ref {propp} (3.), (C) + (An) imply (Ex) and by Proposition \ref {propp} (10'.), (Ex) + (B) imply (BB);
thus, (BB) holds.\\
- by Proposition \ref {propp} (22.), (B) + (C) + (K) + (An) imply (Re); thus, (Re) holds.\\
- Proposition \ref {propp} (4.), (Re) + (Ex) imply (D); thus, (D) holds.\\
- (V') is  (An).
\hfill $\Box$

We obviously have:
\begin{prop}\label{def10-b} A BCK algebra verifies  properties (B), (C), (K), (An).
\end{prop}

By Propositions \ref{def10} and \ref{def10-b}, an algebra $(A, \ra,1)$ verifying the axioms (B), (C), (K),  (An) is an 
equivalent definition,  coming from  BCK logic, of a BCK algebra. Hence, we have the following

\begin{th}
An algebra $(A,\ra,1)$ of type $(2,0)$ is a BCK algebra if and only if properties (B), (C), (K), (An) are satisfied.
\end{th}

                \section{The first nine new  generalizations of BCI and BCK algebras}

          \subsection {Hierarchy 1, including the old algebras BZ, BCC, BCH, BCI, BCK }

It is known that:
\begin{prop}\label{propBCH}
In any BCH algebra,  the following properties hold \cite{C}:
 (D), (M), (N) (see Proposition \ref{propp} (5), (6), (4) respectively).
\end{prop}

We obtain immediately the following corollary (which gives an equivalent, longer definition of BCH algebras):
\begin{cor} For any algebra $(A, \ra, 1)$, we have the equivalence:
 $$(Re) \; + \; (Ex) \;  + \; (An) \; \Leftrightarrow \;  (Re) \;  + \; (M) \;  + \; (Ex) \;  + \; (An). $$
\end{cor}

 We  prove now   the following important result:
\begin{th} 
$$  {\bf BCH} \;   + \; (B) \;  = \;  {\bf BCI}.$$
\end{th}
 
{\bf Proof.}\\
We must prove that (using the first definition of BCI algebras):
$$((Re) \; + \; (Ex) \; + \; (An)) \; + \; (B) \; \Leftrightarrow (BB) \; + \; (D) \; + \; (Re) \; + \; (N) \; + \; (An).$$
$\Longrightarrow$: Suppose (Re) + (Ex) + (An) + (B) holds. By Proposition \ref{propp} (10'.), (Ex) + (B) $\Longrightarrow$ (BB).
By Proposition \ref{propp} (4.), (Re) + (Ex) $\Longrightarrow$ (D). 
By Proposition \ref{propp} (5'.), (Re) + (Ex)  + (An) $\Longrightarrow$ (N). 
Thus, (BB) + (D) + (Re) +  (N) + (An) holds.\\
$\Longleftarrow$: Suppose now that (BB) + (D) + (Re) +  (N) + (An) holds. By Proposition \ref{propp} (21.), 
(BB) + (D) +  (N) + (An) imply (Ex). By Proposition \ref{propp} (10".), (Ex) + (BB) $\Longrightarrow$ (B).
\hfill $\Box$

\begin{rem}\em  We also have, by   Proposition \ref{propp} (10.):
$$  {\bf BCH} \;   + \; (BB) \;  = \;  {\bf BCI}.$$
\end{rem}

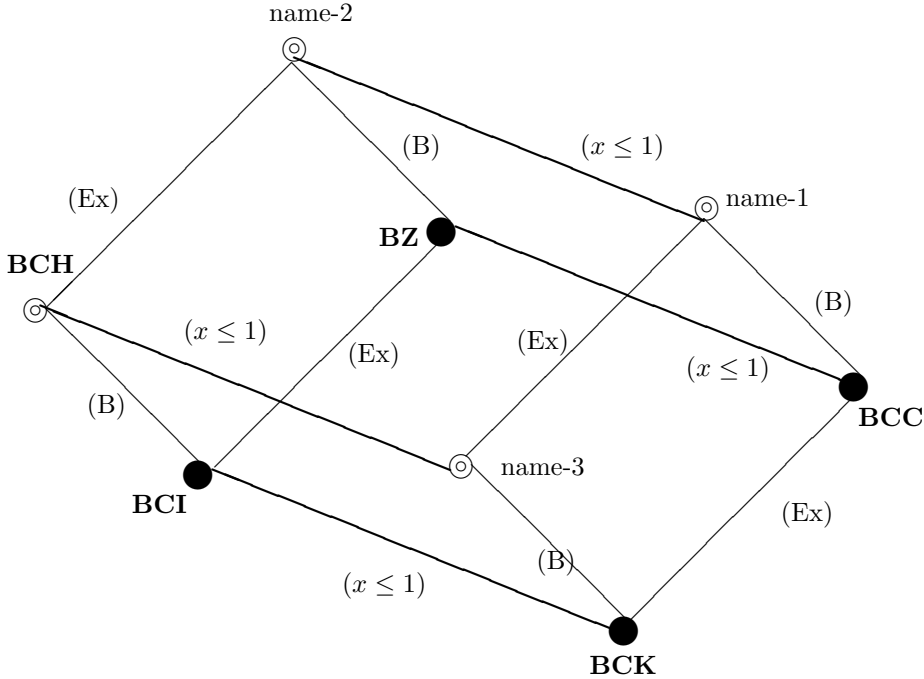
\begin{figure}[htbp]
\begin{center}
\begin{picture}(420,220)(-200,25) 
\put(76,182){\line(-1,-1){90}}    \put(137,119){\line(-1,-1){90}}    
     \put(76,182){\line(1,-1){60}}        
  \put(-12,89){\line(1,-1){60}}    

  
\put(50,25){\makebox(2,2){\circle*{11}}}
 \put(15,47){\makebox(10,10){(B)}}
\put(110,65){\makebox(10,10){(Ex)}}

\put(32,8){\makebox(30,10){{\bf BCK} }}

\put(-17,87){\makebox(2,2){$\bigcirc$}}  \put(-17,87){\makebox(2,2){$\circ$}} 
\put(10,130){\makebox(10,10){(Ex)}}
\put(-54,125){\makebox(10,10){(Ex)}}

\put(0,83){\makebox(30,10){name-3}}

\put(76,185){\makebox(2,2){$\bigcirc$}}  \put(76,185){\makebox(2,2){$\circ$}}

\put(85,185){\makebox(30,10){ name-1 }}
\put(137,117){\makebox(2,2){\circle*{11}}} \put(120,145){\makebox(10,10){(B)}}
                                          
\put(130,102){\makebox(30,10){  {\bf BCC}}}

\put(-80,241){\line(-1,-1){93}}    \put(-19,178){\line(-1,-1){90}}    
     \put(-80,241){\line(1,-1){60}}         \put(-173,148){\line(1,-1){60}}    

\thicklines \put(76,181){\line(-5,2){155}}
\thicklines \put(45,25){\line(-5,2){155}}
\thicklines \put(-20,87){\line(-5,2){155}}
\thicklines \put(137,117){\line(-5,2){155}}

  
\put(-111,84){\makebox(2,2){\circle*{11}}} 

 \put(-155,106){\makebox(10,10){(B)}}

\put(-150,69){ \makebox(30,10){ {\bf BCI}}}
\put(-178,146){\makebox(2,2){$\bigcirc$}}  \put(-178,146){\makebox(2,2){$\circ$}}

 \put(-160,184){\makebox(10,10){(Ex)}}

\put(-190,160){\makebox(30,10){{\bf BCH}}}
 
\put(-80,245){\makebox(2,2){$\bigcirc$}}   \put(-80,245){\makebox(2,2){$\circ$}}

\put(-86,255){\makebox(30,10){name-2 }}

\put(30,204){\makebox(30,10){ ($x \leq 1$) }}
\put(-120,134){\makebox(30,10){($x \leq 1$)}}
\put(70,120){\makebox(30,10){($x \leq 1$)}}
\put(-60,38){\makebox(30,10){($x \leq 1$)}}

\put(-19,176){\makebox(2,2){\circle*{11}}}
 \put(-36,204){\makebox(10,10){(B)}}
\put(-56,170){\makebox(30,10){ {\bf BZ}}}                                         

\end{picture}
 \end{center}
\caption{ Hierarchy 1   }
\label{fig:fig1}
\end{figure}

It follows, based on Hierarchy 0 (from above Figure \ref{fig:fig0}), that we have the 
 hierarchy, called {\it Hierarchy 1}, from Figure \ref{fig:fig1},
 where three new algebras, denoted for the moment by ``name-1", ``name-2", ``name-3", have to exist. Note that
since property (B) implies transitivity (Tr), it follows that the  BCH algebras and these
 three new algebras are not transitive.

             \subsection{Hierarchy 2, including the old algebras BE, pre-BCK, BCC, BCK}
Recall that
 pre-BCK algebra is a BE algebra verifying  property (*). 
 By (\cite{Busneag-Rudeanu}, Proposition 1.3), an equivalent definition of pre-BCK algebras
is given by the axioms: (M), (L), (Ex), (BB).  By (\cite{Busneag-Rudeanu}, Corollary 1.2),
 the class of pre-BCK algebras is equational.
  As the name tells us, pre-BCK algebras are reflexive and transitive algebras, i.e.
the binary relation $\leq$ defined by $x \leq y \Leftrightarrow x \ra y=1$ is 
a pre-order relation (\cite{Busneag-Rudeanu}, Corollary 1.1).

By  (\cite{Busneag-Rudeanu}, Proposition 1.2'), BCK algebras are the same as pre-BCK algebras that satisfy property (An).
By (\cite{Busneag-Rudeanu}, Proposition 1.4), an algebra $(A, \ra, 1)$ is a BCK algebra 
if and only if  it is both a pre-BCK algebra and a BCC algebra. 

\begin{figure}[htbp]
\begin{center}
\begin{picture}(420,420)(-180,25) 
\put(76,182){\line(-1,-1){90}} 
   \put(137,119){\line(-1,-1){90}}    
     \put(76,182){\line(1,-1){60}}  
       \put(-17,89){\line(1,-1){60}}    

  \put(45,25){\line(1,2){100}} 
\put(137,116){\line(1,2){100}} 
\put(-17,89){\line(1,2){100}} 
\put(76,182){\line(1,2){100}} 

\put(176,382){\line(-1,-1){90}} 
  \put(237,316){\line(-1,-1){88}}  
\put(176,382){\line(1,-1){60}}  
\put(83,289){\line(1,-1){60}}   

\put(83,289){\makebox(2,2){$\bigcirc$}}  \put(83,300){\makebox(10,10){{\bf BE} }}
\put(176,382){\makebox(2,2){$\bigcirc$}} \put(176,390){\makebox(10,10){name-4 }}
\put(237,319){\makebox(2,2){$\bigcirc$}} \put(237,319){\makebox(2,2){$\bullet$}}\put(200,315){\makebox(10,10){{\bf pre-BCC} }}
\put(146,228){\makebox(2,2){$\bigcirc$}} \put(146,228){\makebox(2,2){$\bullet$}}\put(146,210){\makebox(10,10){{\bf pre-BCK} }}

\put(30,250){\makebox(30,10){ (An)}}
\put(205,250){\makebox(30,10){ (An)}}
\put(140,300){\makebox(30,10){ (An)}}
\put(80,95){\makebox(10,10){(An)}}

\put(210,350){\makebox(10,10){(B)}}
\put(130,340){\makebox(10,10){(Ex)}}
\put(120,240){\makebox(10,10){(B)}}
\put(178,270){\makebox(10,10){(Ex)}}

\put(45,25){\makebox(2,2){\circle*{11}}}
 \put(-5,47){\makebox(10,10){(B)}}
\put(100,65){\makebox(10,10){(Ex)}}

\put(32,8){\makebox(30,10){{\bf BCK} }}
\put(-72,8){\makebox(30,10){{\it Hierarchy 1} }}

\put(-17,87){\makebox(2,2){$\bigcirc$}}  \put(-17,87){\makebox(2,2){$\circ$}} 

\put(30,120){\makebox(10,10){(Ex)}}

\put(-10,83){\makebox(30,10){{\bf aBE}}}

\put(76,181){\makebox(2,2){$\bigcirc$}} \put(76,181){\makebox(2,2){$\circ$}}
  \put(76,181){\line(-5,2){155}}
\put(85,175){\makebox(30,10){ name-1 }}
\put(137,117){\makebox(2,2){\circle*{11}}} \put(120,145){\makebox(10,10){(B)}}
                                          
\put(130,102){\makebox(30,10){  {\bf BCC}}}
                 \put(180,192){\makebox(30,10){ {\it  Hierarchy 2}}}

\put(-80,241){\line(-1,-1){90}}    \put(-19,178){\line(-1,-1){90}}    
     \put(-80,241){\line(1,-1){60}}         \put(-173,148){\line(1,-1){60}}    

 \thicklines \put(76,181){\line(-5,2){155}}
 \thicklines \put(45,25){\line(-5,2){155}}
 \thicklines \put(-17,87){\line(-5,2){155}}
 \thicklines \put(137,117){\line(-5,2){155}}

  
\put(-111,84){\makebox(2,2){\circle*{11}}} 

 \put(-155,106){\makebox(10,10){(B)}}

\put(-150,69){ \makebox(30,10){ {\bf BCI}}}
\put(-173,146){\makebox(2,2){$\bigcirc$}} \put(-173,146){\makebox(2,2){$\circ$}}

 \put(-150,184){\makebox(10,10){(Ex)}}

\put(-190,160){\makebox(30,10){{\bf BCH}}}
 
\put(-80,240){\makebox(2,2){$\bigcirc$}} \put(-80,240){\makebox(2,2){$\circ$}}

\put(-86,250){\makebox(30,10){name-2 }}

\put(-10,224){\makebox(30,10){ ($x \leq 1$) )}}
\put(-120,134){\makebox(30,10){($x \leq 1$)}}
\put(-60,38){\makebox(30,10){($x \leq 1$)}}

\put(-19,176){\makebox(2,2){\circle*{11}}}
 \put(-36,204){\makebox(10,10){(B)}}
\put(-56,170){\makebox(30,10){ {\bf BZ}}}                                         

\end{picture}
 \end{center}
\caption{ Hierarchies 1 and 2, including Hierarchy 0   }
\label{fig:fig12}
\end{figure}
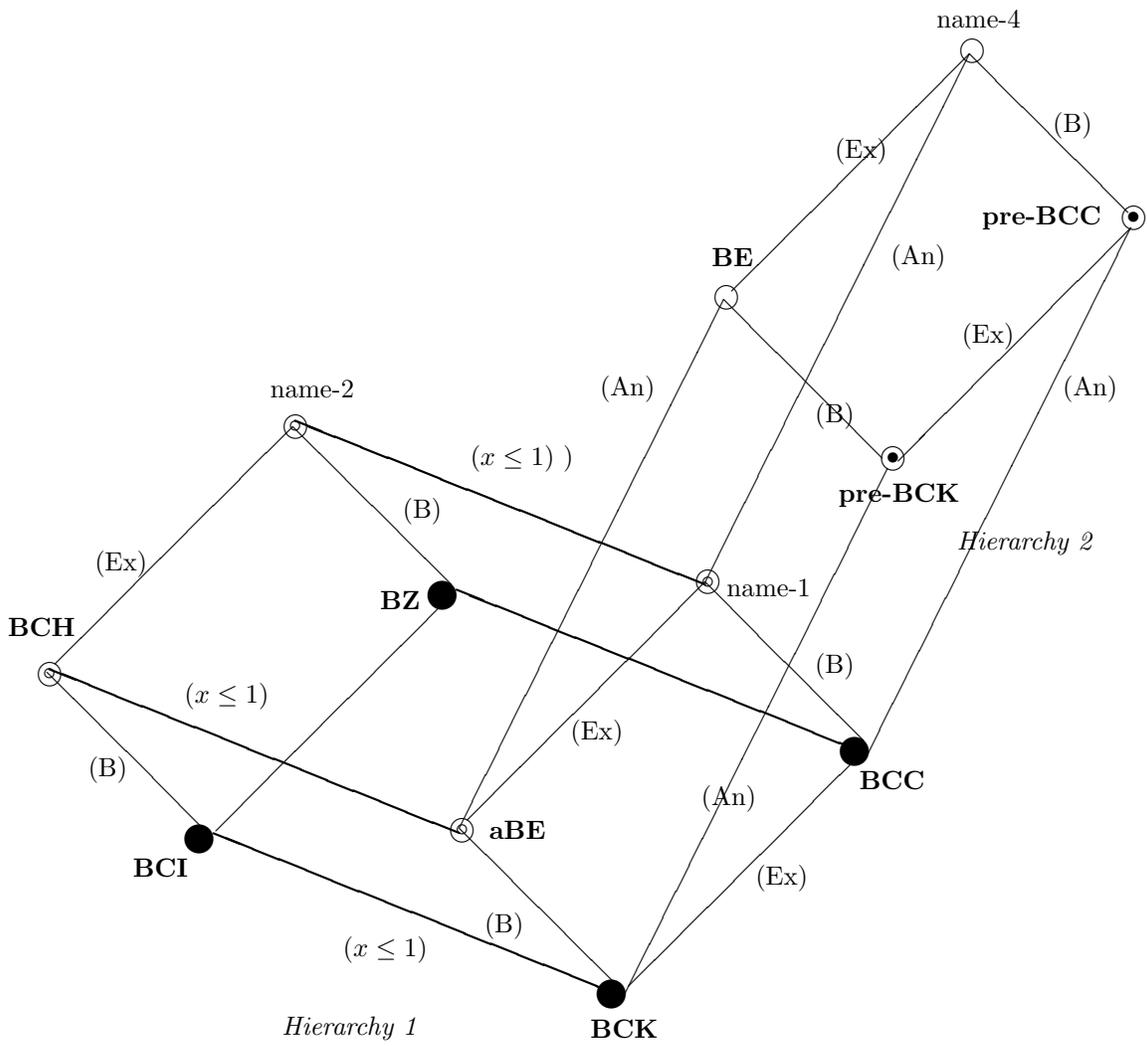

Consequently,
 it is quite obvious that we have a new hierarchy, called {\it Hierarchy 2},
 including the old algebras BE, pre-BCK, BCC and BCK, the algebras ``name-1" and ``name-3", and  two new algebras:
 one is denoted  for the moment by  ``name-4",
 the second is denoted naturally by
``pre-BCC".
 Note that the algebra denoted by ``name-3" in Hierarchy 1 (Figure \ref{fig:fig1}) will receive now naturally 
the name ``aBE". Denote by {\bf pre-BCC, aBE} the classes of pre-BCC and aBE algebras respectively.
 Hierarchy 2  is drawn connected with Hierarchy 1 in Figure \ref{fig:fig12}.  

Hence, we have  the following definitions:
\begin{df}\label{def1}\em ${}$

1. A {\it pre-BCC algebra} is an algebra $(A, \ra,1)$ verifying the axioms: (Re), (M), (L), (B).

2. An {\it aBE algebra} is a BE algebra that is antisymmetrique, i.e.  
is an algebra $(A, \ra,1)$ verifying the axioms: (Re), (M), (L), (Ex), (An). 
\end{df}

As the name says, the binary relation $\leq$ in a pre-BCC algebra is a pre-order too.
 
             \subsection{Hierarchy 3, including the old algebras   BCH, BZ, BCI}

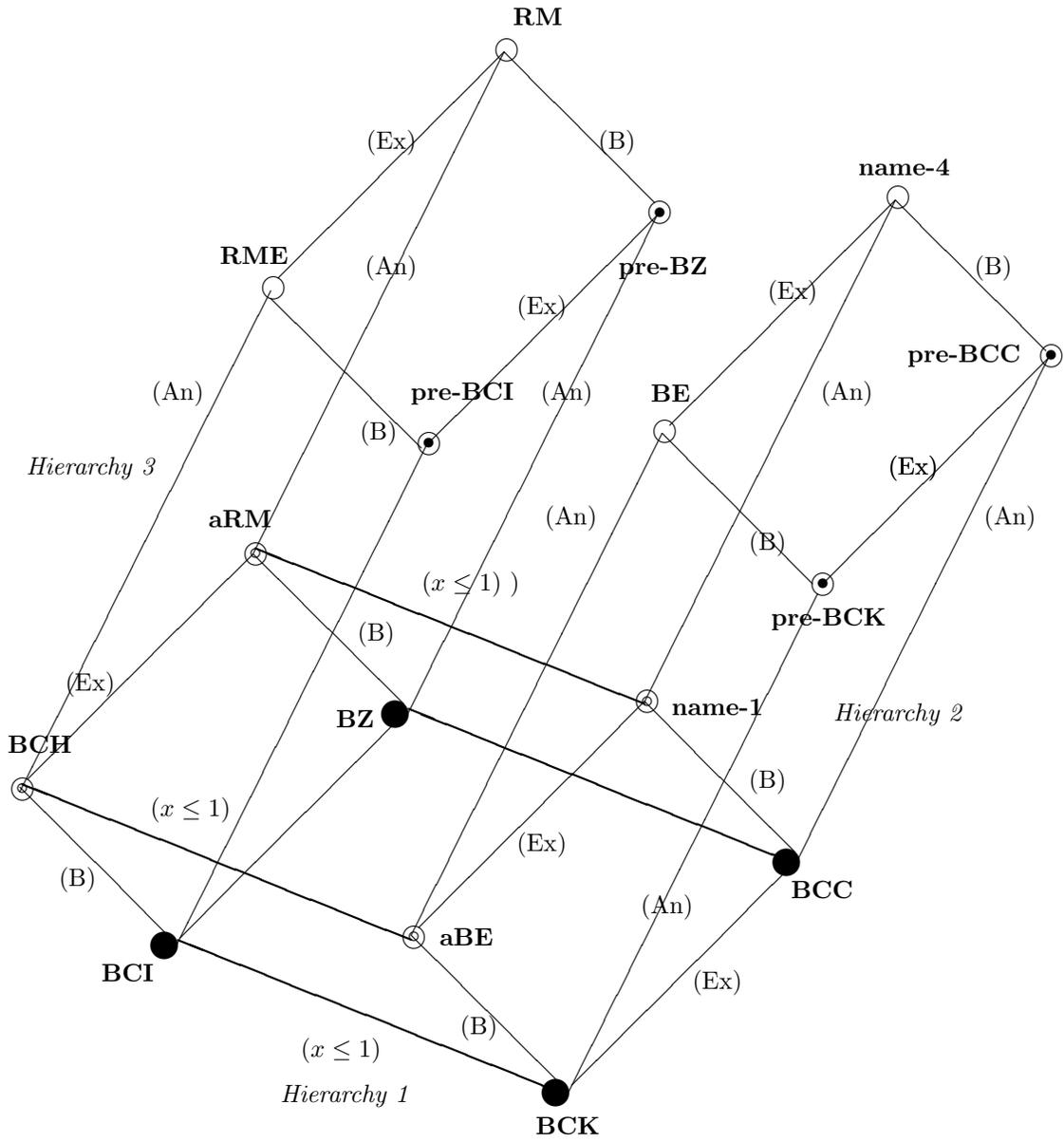
\begin{figure}[htbp]
\begin{center}
\begin{picture}(420,420)(-180,25) 
\put(76,182){\line(-1,-1){90}} 
   \put(137,119){\line(-1,-1){90}}    
     \put(76,182){\line(1,-1){60}}  
       \put(-17,89){\line(1,-1){60}}    

  \put(45,25){\line(1,2){100}} 
\put(137,119){\line(1,2){100}} 
\put(-17,89){\line(1,2){100}} 
\put(76,182){\line(1,2){100}} 

\put(176,382){\line(-1,-1){90}} 
  \put(237,319){\line(-1,-1){90}}  
\put(176,382){\line(1,-1){60}}  
\put(83,289){\line(1,-1){60}}   

\put(83,289){\makebox(2,2){$\bigcirc$}}  \put(83,300){\makebox(10,10){{\bf BE} }}
\put(176,382){\makebox(2,2){$\bigcirc$}} \put(176,390){\makebox(10,10){{\bf name-4} }}
\put(237,319){\makebox(2,2){$\bigcirc$}} \put(237,319){\makebox(2,2){$\bullet$}} \put(200,315){\makebox(10,10){{\bf pre-BCC} }}
\put(146,228){\makebox(2,2){$\bigcirc$}}  \put(146,228){\makebox(2,2){$\bullet$}}  \put(146,210){\makebox(10,10){{\bf pre-BCK} }}

\put(30,250){\makebox(30,10){ (An)}}
\put(205,250){\makebox(30,10){ (An)}}
\put(140,300){\makebox(30,10){ (An)}}
\put(80,95){\makebox(10,10){(An)}}

\put(210,350){\makebox(10,10){(B)}}
\put(130,340){\makebox(10,10){(Ex)}}
\put(120,240){\makebox(10,10){(B)}}
\put(178,270){\makebox(10,10){(Ex)}}

 \put(-173,146){\line(1,2){100}} 
\put(-80,241){\line(1,2){100}} 
\put(-111,84){\line(1,2){100}} 
\put(-19,176){\line(1,2){100}} 

\put(-73,346){\makebox(2,2){$\bigcirc$}}  \put(-83,355){\makebox(10,10){{\bf RME} }}
\put(20,441){\makebox(2,2){$\bigcirc$}} \put(30,450){\makebox(10,10){{\bf RM} }}
\put(-11,284){\makebox(2,2){$\bigcirc$}} \put(-11,284){\makebox(2,2){$\bullet$}} \put(0,300){\makebox(10,10){{\bf pre-BCI} }}
\put(81,376){\makebox(2,2){$\bigcirc$}}  \put(81,376){\makebox(2,2){$\bullet$}} \put(80,350){\makebox(10,10){{\bf pre-BZ} }}

\put(20,441){\line(-1,-1){90}} 
  \put(81,376){\line(-1,-1){90}}  
\put(20,441){\line(1,-1){60}}  
\put(-73,343){\line(1,-1){60}}   

\put(60,400){\makebox(10,10){(B)}} \put(30,334){\makebox(10,10){(Ex)}}
\put(-30,400){\makebox(10,10){(Ex)}}  \put(-30,350){\makebox(10,10){(An)}}
\put(-115,300){\makebox(10,10){(An)}}
\put(40,300){\makebox(10,10){(An)}} \put(-35,284){\makebox(10,10){(B)}}

\put(178,270){\makebox(10,10){(Ex)}}

\put(45,25){\makebox(2,2){\circle*{11}}}
 \put(5,47){\makebox(10,10){(B)}}
\put(100,65){\makebox(10,10){(Ex)}}

                             \put(-52,20){\makebox(20,10){{\it Hierarchy 1} }}
   \put(32,8){\makebox(30,10){{\bf BCK} }}

\put(-17,87){\makebox(2,2){$\bigcirc$}}  \put(-17,87){\makebox(2,2){$\circ$}} 

\put(30,120){\makebox(10,10){(Ex)}}

\put(-10,83){\makebox(30,10){{\bf aBE}}}

\put(76,181){\makebox(2,2){$\bigcirc$}} \put(76,181){\makebox(2,2){$\circ$}}
 
\put(90,175){\makebox(30,10){ {\bf name-1} }}
\put(137,117){\makebox(2,2){\circle*{11}}} 
\put(120,145){\makebox(10,10){(B)}}
                                          
\put(130,102){\makebox(30,10){  {\bf BCC}}}
                  \put(160,172){\makebox(30,10){ {\it  Hierarchy 2}}}

\put(-80,241){\line(-1,-1){90}} 
   \put(-19,178){\line(-1,-1){90}}    
     \put(-80,241){\line(1,-1){60}} 
        \put(-173,148){\line(1,-1){60}}    

\thicklines \put(76,181){\line(-5,2){155}}
 \thicklines \put(45,25){\line(-5,2){155}}
 \thicklines \put(-17,87){\line(-5,2){155}}
 \thicklines \put(137,117){\line(-5,2){155}}

  
\put(-111,84){\makebox(2,2){\circle*{11}}} 

 \put(-155,106){\makebox(10,10){(B)}}

\put(-150,69){ \makebox(30,10){ {\bf BCI}}}
\put(-173,146){\makebox(2,2){$\bigcirc$}} \put(-173,146){\makebox(2,2){$\circ$}}

 \put(-150,184){\makebox(10,10){(Ex)}}

\put(-180,160){\makebox(30,10){{\bf BCH}}}
                       \put(-160,270){\makebox(30,10){{\it Hierarchy 3}}}

\put(-80,240){\makebox(2,2){$\bigcirc$}} \put(-80,240){\makebox(2,2){$\circ$}}

\put(-100,250){\makebox(30,10){ {\bf aRM }}}

\put(-10,224){\makebox(30,10){ ($x \leq 1$) )}}
\put(-120,134){\makebox(30,10){($x \leq 1$)}}
\put(-60,38){\makebox(30,10){($x \leq 1$)}}

\put(-19,176){\makebox(2,2){\circle*{11}}}
 \put(-36,204){\makebox(10,10){(B)}}
\put(-56,170){\makebox(30,10){ {\bf BZ}}}                                         

\end{picture}
 \end{center}
\caption{ Hierarchies 1,2,3, including Hierarchy 0  }
\label{fig:fig123}
\end{figure}

From above  new algebras and above Hierachies 0, 1, 2, we can  
 obtain obviously new algebras and hierarchies.
Thus, we obtain the new hierarchy  called {\it Hierarchy 3}, including: the old algebras   BCH, BZ and BCI,
the algebra ``name-2", and four new algebras, called ``RM", ``RME", ``pre-BZ" and ``pre-BCI".
 Note that the algebra denoted by ``name-2" in Hierarchy 1 (Figure \ref{fig:fig1}) will receive now naturally 
the name ``aRM". Denote by {\bf RM, RME, pre-BZ, pre-BCI, aRM} the classes of  RM, RME, pre-BZ, pre-BCI, aRM 
algebras respectively.
 Hierarchy 3 is drawn connected with Hierarchies 1 and 2 in Figure \ref{fig:fig123}.

Hence, we have  the following obvious  definitions:

\begin{df}\label{def2}\em ${}$

3. A {\it RM algebra} is an algebra $(A, \ra,1)$ verifying the axioms: (Re), (M).

4. A {\it pre-BZ algebra} is an algebra $(A, \ra,1)$ verifying the axioms: (Re), (M), (B).

5. An {\it aRM algebra} is an algebra $(A, \ra,1)$ verifying the axioms: (Re), (M), (An).

6. A  {\it RME algebra} is an algebra $(A, \ra,1)$ verifying the axioms: (Re), (M), (Ex).

7. A {\it pre-BCI algebra} is an algebra $(A, \ra,1)$ verifying the axioms: (Re), (M), (Ex), (B).
\end{df}

We check and we have, indeed, that:

{\bf Pre-BZ} + (An) = {\bf BZ},

{\bf Pre-BCI} + (An) = {\bf BCI},

{\bf RME} + (An) = {\bf BCH},

{\bf aRM} + (B) = {\bf BZ}, 

{\bf aRM} + (Ex) = {\bf BCH}.
                              \subsection{Hierarchy 4 and the final, global  hierarchy}

We can now complete  the Hierarchies 0, 1, 2, 3  by adding  a new hierarchy, called {\it Hierarchy 4}, 
which connects obviously 
(by property (L)) the algebras RM, pre-BZ, RME, pre-BCI from Hierarchy 3 and the algebras name-4, pre-BCC,
BE, pre-BCK from Hierarchy 2. 
Note that the algebra denoted by ``name-4" in Hierarchy 2 (Figure \ref{fig:fig12}) will receive now naturally 
the name ``RML", hence the algebra denoted by ``name-1" receives naturally the name ``aRML".
 Denote by {\bf RML, aRML} the class of  RML and aRML algebras, respectively.
 Hierarchy 4 can now be easily drawn. 

Before drawing the final, global hierarchy, connecting the Hierarchies 0, 1, 2, 3, 4, recall
that by Corollary \ref{corA}, in BZ, BCC algebras we have:  $ (Ex) \; \Leftrightarrow \; (BB) $\\
and  that by Theorem \ref{th1}, we  have:
  (Re) + (M) + (Ex) imply $ (B) \; \Leftrightarrow \; (BB) \; \Leftrightarrow \; (*).$

Consequently, we have:
\begin{cor}
In RM and RML  algebras satisfying property (Ex), i.e. in RME, BCH and BE, aBE algebras, we have:
$$ (B) \; \Leftrightarrow \; (BB) \; \Leftrightarrow \; (*).$$
\end{cor}

Note that this corollary gives us equivalent definitions for the four algebras 
pre-BCI, BCI and pre-BCK, BCK.

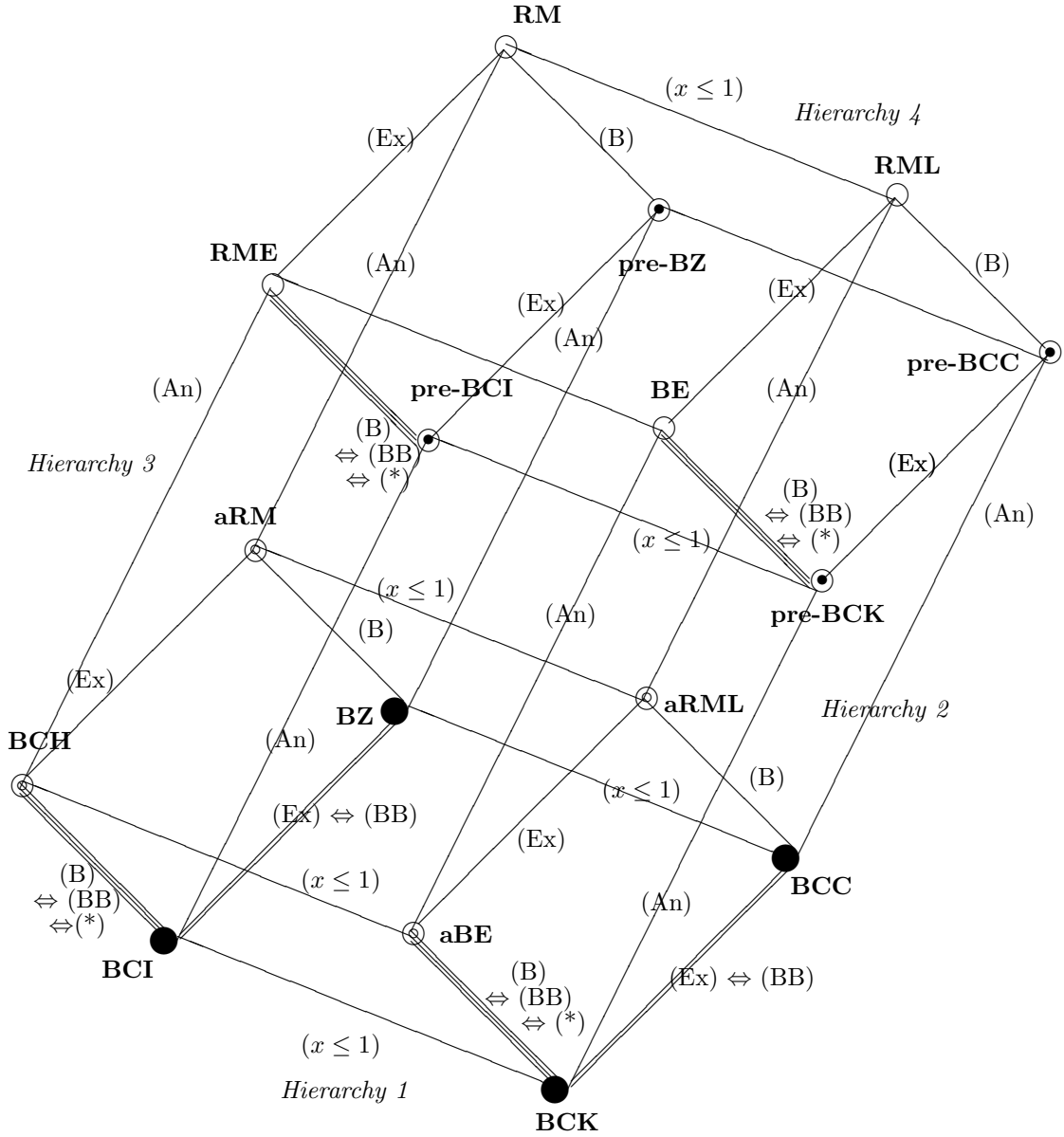
\begin{figure}[htbp]
\begin{center}
\begin{picture}(420,420)(-180,25) 
 \put(-52,20){\makebox(20,10){{\it Hierarchy 1} }}
\put(155,172){\makebox(30,10){ {\it  Hierarchy 2}}}
\put(-160,270){\makebox(30,10){{\it Hierarchy 3}}}
\put(156,410){\makebox(10,10){{\it Hierarchy 4}}}

\put(76,182){\line(-1,-1){90}} 
                                              
     \put(76,182){\line(1,-1){60}}  
           \put(-17,89){\line(1,-1){60}}    
           \put(-17,87){\line(1,-1){60}}    
           \put(-17,85){\line(1,-1){60}}    

  \put(76,181){\line(-5,2){155}}
  \put(45,25){\line(-5,2){155}}
 \put(-17,87){\line(-5,2){155}}
  \put(137,117){\line(-5,2){155}}
  \put(45,25){\line(1,2){100}} 
\put(137,119){\line(1,2){100}} 
\put(-17,89){\line(1,2){100}} 
\put(76,182){\line(1,2){100}} 

\put(176,382){\line(-1,-1){90}} 
  \put(237,319){\line(-1,-1){90}}  
\put(176,382){\line(1,-1){60}}  

\put(83,289){\line(1,-1){59}}   
\put(83,287){\line(1,-1){59}}   
\put(83,285){\line(1,-1){59}}

\put(83,289){\makebox(2,2){$\bigcirc$}}  \put(83,300){\makebox(10,10){{\bf BE} }}
\put(176,382){\makebox(2,2){$\bigcirc$}} \put(176,390){\makebox(10,10){{\bf RML}}}
\put(237,319){\makebox(2,2){$\bigcirc$}} \put(237,319){\makebox(2,2){$\bullet$}} \put(200,310){\makebox(10,10){{\bf pre-BCC} }}
\put(146,228){\makebox(2,2){$\bigcirc$}}  \put(146,228){\makebox(2,2){$\bullet$}}  \put(146,210){\makebox(10,10){{\bf pre-BCK} }}

\put(30,210){\makebox(30,10){ (An)}}
\put(205,250){\makebox(30,10){ (An)}}
\put(80,95){\makebox(10,10){(An)}}

\put(130,300){\makebox(10,10){(An)}}

\put(210,350){\makebox(10,10){(B)}}
\put(130,340){\makebox(10,10){(Ex)}}

\put(82,240){\makebox(10,10){($x \leq 1$)}}

\put(135,260){\makebox(10,10){(B) }}
\put(135,250){\makebox(10,10){ $\Leftrightarrow$ (BB)}}
\put(135,240){\makebox(10,10){ $\Leftrightarrow$  (*)}}

\put(178,270){\makebox(10,10){(Ex)}}
 \put(-173,146){\line(1,2){100}} 
\put(-80,241){\line(1,2){100}} 
\put(-111,84){\line(1,2){100}} 
\put(-19,176){\line(1,2){100}} 

\put(20,441){\line(-1,-1){90}} 
  \put(81,376){\line(-1,-1){90}}  
\put(20,441){\line(1,-1){60}}  
         \put(-73,343){\line(1,-1){58}}   
         \put(-73,341){\line(1,-1){58}}   
         \put(-73,345){\line(1,-1){58}}

\put(-73,346){\makebox(2,2){$\bigcirc$}}  \put(-87,355){\makebox(10,10){{\bf RME} }}
\put(20,441){\makebox(2,2){$\bigcirc$}} \put(30,450){\makebox(10,10){{\bf RM} }}
\put(-11,284){\makebox(2,2){$\bigcirc$}} \put(-11,284){\makebox(2,2){$\bullet$}} \put(0,300){\makebox(10,10){{\bf pre-BCI} }}
\put(81,376){\makebox(2,2){$\bigcirc$}}  \put(81,376){\makebox(2,2){$\bullet$}} \put(80,350){\makebox(10,10){{\bf pre-BZ} }}

\put(60,400){\makebox(10,10){(B)}}
 \put(30,334){\makebox(10,10){(Ex)}}
\put(-30,400){\makebox(10,10){(Ex)}}  
\put(-30,350){\makebox(10,10){(An)}}
\put(-115,300){\makebox(10,10){(An)}}
\put(45,320){\makebox(10,10){(An)}} 

\put(-37,284){\makebox(10,10){(B)}}
\put(-35,274){\makebox(10,10){$\Leftrightarrow$ (BB)}}
\put(-35,264){\makebox(10,10){$\Leftrightarrow$ (*)}}

\put(178,270){\makebox(10,10){(Ex)}}
 \put(176,381){\line(-5,2){155}}
  \put(145,225){\line(-5,2){155}}
 \put(83,289){\line(-5,2){155}}
  \put(237,317){\line(-5,2){155}}

\put(85,420){\makebox(30,10){ ($x \leq 1$) }}

\put(45,25){\makebox(2,2){\circle*{11}}}

 \put(25,67){\makebox(10,10){(B)}}
\put(25,57){\makebox(10,10){$\Leftrightarrow$ (BB)}}
\put(35,47){\makebox(10,10){$\Leftrightarrow$ (*)}}

\put(32,8){\makebox(30,10){{\bf BCK} }}

\put(-17,87){\makebox(2,2){$\bigcirc$}}  \put(-17,87){\makebox(2,2){$\circ$}} 

\put(30,120){\makebox(10,10){(Ex)}}

\put(-10,83){\makebox(30,10){{\bf aBE}}}

\put(76,181){\makebox(2,2){$\bigcirc$}} \put(76,181){\makebox(2,2){$\circ$}}
 
\put(85,175){\makebox(30,10){ {\bf  aRML} }}
\put(137,117){\makebox(2,2){\circle*{11}}} 
\put(120,145){\makebox(10,10){(B)}}
                                          
\put(130,102){\makebox(30,10){  {\bf BCC}}}

\put(-80,241){\line(-1,-1){90}} 
                                                      
     \put(-80,241){\line(1,-1){60}} 
        \put(-173,148){\line(1,-1){60}}    
        \put(-173,146){\line(1,-1){60}}    
        \put(-173,144){\line(1,-1){60}}    

  
\put(-111,84){\makebox(2,2){\circle*{11}}} 

 \put(-155,106){\makebox(10,10){(B)}}
\put(-155,96){\makebox(10,10){$\Leftrightarrow$ (BB)}}
\put(-155,86){\makebox(10,10){$\Leftrightarrow$(*)}}

\put(-150,69){ \makebox(30,10){ {\bf BCI}}}
\put(-173,146){\makebox(2,2){$\bigcirc$}} \put(-173,146){\makebox(2,2){$\circ$}}

 \put(-150,184){\makebox(10,10){(Ex)}}

\put(-180,160){\makebox(30,10){{\bf BCH}}}
 
\put(-80,240){\makebox(2,2){$\bigcirc$}} \put(-80,240){\makebox(2,2){$\circ$}}

\put(-96,250){\makebox(30,10){{\bf aRM }}}

\put(-80,160){\makebox(30,10){ (An) }}

\put(-30,220){\makebox(30,10){ ($x \leq 1$) }}
\put(-60,104){\makebox(30,10){($x \leq 1$)}}
\put(60,140){\makebox(30,10){($x \leq 1$)}}
\put(-60,38){\makebox(30,10){($x \leq 1$)}}

\put(-19,176){\makebox(2,2){\circle*{11}}}
 \put(-36,204){\makebox(10,10){(B)}}
\put(-56,170){\makebox(30,10){ {\bf BZ}}}      
                                                 \put(-19,178){\line(-1,-1){90}}
                                                 \put(-19,176){\line(-1,-1){90}} 
                                                \put(137,119){\line(-1,-1){90}} 
                                             \put(137,117){\line(-1,-1){90}}                                      
                                              \put(-60,130){\makebox(30,10){ (Ex) $\Leftrightarrow$  (BB)}}
                                              \put(110,65){\makebox(10,10){(Ex) $\Leftrightarrow $ (BB)}}

\end{picture}
 \end{center}
\caption{ The global hierarchy (connecting Hierarchies 0, 1, 2, 3, 4) }
\label{fig:fig1234}
\end{figure}

We are now in position to present the final, global hierarchy (which connects the Hierachies 0, 1, 2, 3, 4)
in Figure \ref{fig:fig1234}.

Hence, we can introduce  obviously the following  new definitions:
\begin{df}\label{def3}\em ${}$

8. A {\it RML algebra} is an algebra $(A, \ra,1)$ verifying the axioms: (Re), (M), (L).

9. An {\it aRML algebra} is an algebra $(A, \ra,1)$ verifying the axioms: (Re), (M), (L), (An), i.e. is a RML algebra verifying (An) (Antisymmetry).
\end{df}

\begin{rems} \em ${}$

(i) All the RM algebras from Figure \ref{fig:fig1234} are generalizations of BCK algebras. More precisely,
 the RM algebras from Hierarchy 3 are generalizations of BCI algebras, while the RM algebras  from Hierarchy 2 (i.e. the RML algebras)
 are called {\it proper}
  generalizations of BCK algebras.

(ii) All RM algebras are at least reflexive.

(iii) While the hierarchy of the four old  ordered algebras (BZ, BCI, BCC, BCK) is presented in Figure \ref{fig:fig0},
we  can remark the similar  hierarchy of the four pre-ordered  
corresponding algebras (pre-BZ, pre-BCI, pre-BCC, pre-BCK) in Figure \ref{fig:fig1234}.
\end{rems}

                \section {Other twenty two new RM and RML algebras,\\
 generalizations of BCI and BCK algebras respectively}

Looking for proper examples of the above mentioned RM and RML algebras,
 we discovered that the classes  of RM and RML algebras are much richer. Thus,
we have introduced the following   eleven  new RM algebras  and  eleven new RML algebras.

            \subsection{New  RM and RML algebras, without  property (Ex)}

                       {\bf $\bullet$ New RM algebras, without (Ex)}

Define the following  nine  new algebras:
\begin{df}\em${}$

10. A   {\it tRM algebra}   is a RM algebra verifying (Tr).

11. A   {\it *RM algebra}   is a RM algebra verifying (*).

12. A   {\it RM** algebra}   is a RM algebra verifying (**).

13. A   {\it *RM** algebra}   is a RM algebra verifying (*), (**).

14. A   {\it pre-BBBZ algebra}   is a RM algebra verifying (B), (BB) (hence  also (*), (**), (Tr)).
\end{df}

\begin{df}\em${}$

15. An   {\it oRM algebra}   is an  aRM algebra verifying (Tr).

16. An   {\it *aRM algebra}   is an  aRM algebra verifying (*).

17. An   {\it aRM** algebra}   is an  aRM algebra verifying (**).

18. An   {\it *aRM** algebra}   is an  aRM algebra verifying (*), (**).

\end{df}

{\bf $\bullet$  New  RML algebras, without (Ex)}

Define the following  corresponding  (connected by (L) with the previous ones) nine new algebras:

\begin{df}\em${}$

19. A   {\it tRML algebra}   is a RML algebra verifying (Tr).

20. A   {\it *RML algebra}   is a RML algebra verifying (*).

21. A   {\it RML** algebra}   is a RML algebra verifying (**).

22. A   {\it *RML** algebra}   is a RML algebra verifying (*), (**).

23. A  {\it pre-BBBCC algebra} is a RML algebra verifying (B), (BB) (hence also (*), (**), (Tr)).
\end{df}

\begin{df}\em${}$

24. An   {\it oRML algebra}   is an aRML algebra verifying (Tr).

25. An   {\it *aRML algebra}   is an aRML algebra verifying (*).

26. An   {\it aRML** algebra}   is an aRML algebra verifying (**).

27. An   {\it *aRML** algebra}   is an aRML algebra verifying (*), (**).
\end{df}

\begin{cor}\em  By Corollary  \ref{corA}, we  have:\\
{\bf pre-BBBZ} + (An) = {\bf pre-BZ} + (BB) + (An) = {\bf BZ} + (BB) =  {\bf BZ} + (Ex) = {\bf BCI},\\
{\bf pre-BBBCC} + (An) = {\bf pre-BCC} + (BB) + (An) = {\bf BCC} + (BB)  = {\bf BCC} + (Ex) = {\bf BCK}.
\end{cor}
 
Hence, we have, for example, the Hierarchies 3.1 and 2.1 
(determined by RM and by RML algebras, considering  (B), (*) and (Tr) properties),  from the next  Figure \ref{fig:fig23sT},
 and we also have 
the Hierarchies 3.2 and 2.2 (determined by RM and by RML algebras, considering  (*) and (**) properties)
 from the next Figure \ref{fig:figrml***}. In both Figures, the connection by (L) between the Hierarchies 3.1 and 2.1
and between the Hierarchies 3.2 and 2.2 is not drawn, in order to simplify the figure.

\begin{figure}[htbp]
\begin{center}
\begin{picture}(420,400)(-180,-25)

       \put(-17,89){\line(1,-1){60}}    

  \put(45,25){\line(1,2){100}}                               
\put(-17,89){\line(1,2){100}}

      \put(83,289){\line(1,-1){97}}   
\put(83,289){\line(2,-1){85}}   
\put(83,289){\line(3,-1){102}}   
  \put(146,228){\line(2,1){46}}  
\put(171,244){\makebox(2,2){$\bigcirc$}} \put(171,244){\makebox(2,2){$\bullet$}}\put(185,230){\makebox(10,10){{\bf *RML} }}
\put(190,254){\makebox(2,2){$\bigcirc$}} \put(190,254){\makebox(2,2){$\bullet$}}\put(210,245){\makebox(10,10){{\bf tRML} }}
\put(139,270){\makebox(30,10){ (Tr)}}
\put(125,245){\makebox(30,10){ (*)}}
\put(-21,89){\line(2,-1){85}}   
\put(-21,89){\line(3,-1){102}}   
  \put(45,25){\line(2,1){46}}    
\put(90,50){\makebox(2,2){\circle*{11}}}
\put(70,41){\makebox(2,2){\circle*{11}}}
  \put(90,50){\line(1,2){100}} 
 \put(70,41){\line(1,2){100}} 
\put(95,40){\makebox(30,10){{\bf oRML} }}
\put(70,21){\makebox(30,10){{\bf *aRML} }}
\put(35,69){\makebox(30,10){ (Tr)}}
\put(26,44){\makebox(30,10){ (*)}}

\put(83,291){\makebox(2,2){$\bigcirc$}}  \put(83,300){\makebox(10,10){{\bf RML} }}
\put(146,228){\makebox(2,2){$\bigcirc$}} \put(146,228){\makebox(2,2){$\bullet$}} 
                              \put(110,220){\makebox(10,10){{\bf pre-BCC} }}
\put(182,187){\makebox(2,2){$\bigcirc$}}
\put(182,187){\makebox(2,2){$\bullet$}}
\put(180,167){\makebox(10,10){{\bf pre-BBBCC}}}
\put(152,197){\makebox(10,10){(BB)}}

\put(10,200){\makebox(30,10){ (An)}}
\put(105,130){\makebox(30,10){ (An)}}
                  \put(95,180){\makebox(30,10){ (An)}}
\put(125,95){\makebox(10,10){(An)}}

\put(95,250){\makebox(10,10){(B)}}

\put(45,25){\makebox(2,2){\circle*{11}}}
 \put(-5,50){\makebox(10,10){(B)}}


\put(32,8){\makebox(30,10){{\bf BCC} }}

\put(-21,87){\makebox(2,2){$\bigcirc$}}  \put(-21,87){\makebox(2,2){$\circ$}} 
\put(-45,70){\makebox(30,10){{\bf aRML}}}
 \put(-172,149){\line(1,-1){60}}    

  \put(-110,85){\line(1,2){100}}                               
\put(-172,149){\line(1,2){100}}

      \put(-72,349){\line(1,-1){97}}   

\put(-72,349){\line(2,-1){85}}   
\put(-72,349){\line(3,-1){102}}   
  \put(-9,288){\line(2,1){46}}  
\put(16,304){\makebox(2,2){$\bigcirc$}} \put(16,304){\makebox(2,2){$\bullet$}}\put(30,290){\makebox(10,10){{\bf *RM} }}
\put(35,314){\makebox(2,2){$\bigcirc$}} \put(35,314){\makebox(2,2){$\bullet$}}\put(55,315){\makebox(10,10){{\bf tRM} }}
\put(-16,330){\makebox(30,10){ (Tr)}}
\put(-30,305){\makebox(30,10){ (*)}}
\put(-176,149){\line(2,-1){85}}   
\put(-176,149){\line(3,-1){102}}   
  \put(-110,85){\line(2,1){46}}    
\put(-65,110){\makebox(2,2){\circle*{11}}}
\put(-85,100){\makebox(2,2){\circle*{11}}}
  \put(-65,110){\line(1,2){100}} 
 \put(-85,101){\line(1,2){100}} 
\put(-60,100){\makebox(30,10){{\bf oRM} }}
\put(-75,91){\makebox(30,10){{\bf *aRM} }}
\put(-120,129){\makebox(30,10){ (Tr)}}
\put(-129,104){\makebox(30,10){ (*)}}

\put(-72,351){\makebox(2,2){$\bigcirc$}}  \put(-72,360){\makebox(10,10){{\bf RM} }}
\put(-9,288){\makebox(2,2){$\bigcirc$}} \put(-9,288){\makebox(2,2){$\bullet$}}\put(-40,280){\makebox(10,10){{\bf pre-BZ} }}

\put(27,248){\makebox(2,2){$\bigcirc$}}
 \put(27,248){\makebox(2,2){$\bullet$}}
\put(22,228){\makebox(10,10){{\bf pre-BBBZ} }}
\put(-4,258){\makebox(10,10){ (BB) }}

\put(-145,260){\makebox(30,10){ (An)}}
\put(-40,190){\makebox(30,10){ (An)}}
                  \put(-60,240){\makebox(30,10){ (An)}}
\put(-30,155){\makebox(10,10){(An)}}

\put(-60,310){\makebox(10,10){(B)}}

\put(-103,85){\makebox(2,2){\circle*{11}}}
 \put(-160,110){\makebox(10,10){(B)}}


\put(-123,68){\makebox(30,10){{\bf BZ} }}

\put(-176,147){\makebox(2,2){$\bigcirc$}}  \put(-176,147){\makebox(2,2){$\circ$}} 
\put(-190,130){\makebox(30,10){{\bf aRM}}}
                                         \put(182,187){\line(-1,-2){100}} 
                                                                 \put(27,248){\line(-1,-2){100}}

                                        \put(-70,47){\makebox(2,2){\circle*{11}}}
                                             \put(-91,28){\makebox(30,10){{\bf BCI} }}  
                                              \put(-131,50){\makebox(30,10){(Ex) $\Leftrightarrow$ (BB) }}    
                                              \put(-105,82){\line(1,-1){32}}                                                                                                 
                                              \put(-105,80){\line(1,-1){32}} 

                                               \put(47,24){\line(1,-1){34}} 
                                               \put(47,22){\line(1,-1){34}} 
                                               \put(87,-14){\makebox(2,2){\circle*{11}}}
                                                 \put(85,-35){\makebox(30,10){{\bf BCK} }}  
                                                 \put(23,-8){\makebox(30,10){(Ex) $\Leftrightarrow$ (BB) }}   
                                                 \put(152,65){\makebox(2,2){(An)}}
\put(148,315){\makebox(10,10){{\it Hierarchy 2.1 } }}
                            \put(-18,370){\makebox(10,10){{\it Hierarchy 3.1} }}

\end{picture}
 \end{center}
\caption{ Hierarchies 3.1 and 2.1,  determined by RM and by RML algebras, considering  (B), (*) and (Tr) properties   }
\label{fig:fig23sT}
\end{figure}
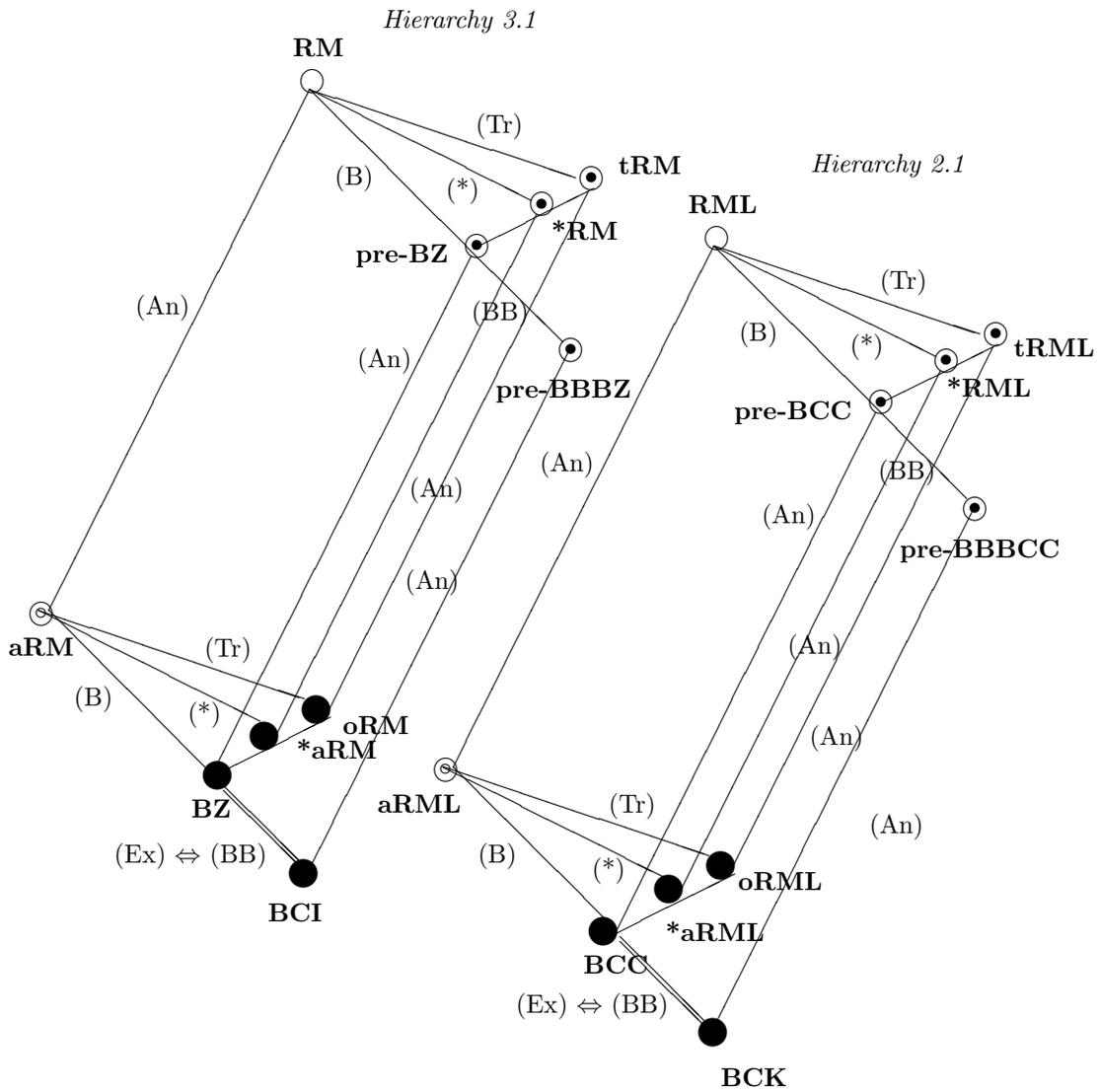

\begin{figure}[htbp]
\begin{center}
\begin{picture}(420,340)(-180,25)

       \put(-17,89){\line(1,-1){32}}    
                 
                      \put(15,55){\line(1,2){100}}                               
\put(-17,89){\line(1,2){100}}

      \put(83,289){\line(1,-1){32}}   
       
\put(83,289){\line(3,-1){65}}

\put(171,244){\line(-1,1){23}}
\put(171,244){\line(-4,1){56}}
\put(171,244){\makebox(2,2){$\bigcirc$}} \put(171,244){\makebox(2,2){$\bullet$}}\put(185,222){\makebox(10,10){{\bf *RML**} }}
\put(150,270){\makebox(2,2){$\bigcirc$}} \put(150,270){\makebox(2,2){$\bullet$}}\put(170,280){\makebox(10,10){{\bf RML**} }}
\put(99,280){\makebox(30,10){ (**)}}  \put(135,238){\makebox(30,10){ (**)}}
\put(-21,89){\line(3,-1){65}}   
    
\put(70,44){\line(-1,1){23}}
\put(70,44){\line(-4,1){56}}
 
             \put(55,70){\makebox(2,2){\circle*{11}}}
\put(70,41){\makebox(2,2){\circle*{11}}}
                             \put(50,70){\line(1,2){100}} 
 \put(70,41){\line(1,2){100}} 
                \put(75,80){\makebox(30,10){{\bf aRML**} }}
\put(70,21){\makebox(30,10){{\bf *aRML**} }}
\put(-5,80){\makebox(30,10){ (**)}}  \put(32,38){\makebox(30,10){ (**)}}
                            \put(143,300){\makebox(10,10){{\it Hierarchy 2.2} }}

\put(83,291){\makebox(2,2){$\bigcirc$}}  \put(83,300){\makebox(10,10){{\bf RML} }}
\put(115,252){\makebox(2,2){$\bigcirc$}} \put(115,252){\makebox(2,2){$\bullet$}}\put(88,240){\makebox(10,10){{\bf *RML} }}

\put(10,200){\makebox(30,10){ (An)}}
\put(105,130){\makebox(30,10){ (An)}}
                  \put(95,180){\makebox(30,10){ (An)}}
\put(35,125){\makebox(10,10){(An)}}

\put(95,260){\makebox(10,10){(*)}}   \put(165,252){\makebox(10,10){(*)}}

\put(20,51){\makebox(2,2){\circle*{11}}}
 \put(-5,60){\makebox(10,10){(*)}}  \put(55,52){\makebox(10,10){(*)}}

\put(-10,30){\makebox(30,10){{\bf *aRML} }}

\put(-21,87){\makebox(2,2){$\bigcirc$}}  \put(-21,87){\makebox(2,2){$\circ$}} 
\put(-45,95){\makebox(30,10){{\bf aRML}}}
 \put(-172,149){\line(1,-1){32}}    
     
                       \put(-108,132){\line(1,2){100}}                               
\put(-172,149){\line(1,2){100}}

      \put(-72,349){\line(1,-1){32}}   
      
\put(-72,349){\line(3,-1){65}}   
  
\put(16,304){\line(-1,1){23}}
\put(16,304){\line(-4,1){56}}

\put(-40,313){\makebox(2,2){$\bigcirc$}} \put(-40,313){\makebox(2,2){$\bullet$}}\put(30,282){\makebox(10,10){{\bf *RM**} }}
\put(20,303){\makebox(2,2){$\bigcirc$}} \put(20,303){\makebox(2,2){$\bullet$}}
\put(15,330){\makebox(10,10){{\bf RM**} }}
\put(-56,340){\makebox(30,10){ (**)}}   \put(-16,295){\makebox(30,10){ (**)}}
\put(-176,149){\line(3,-1){65}}   
   
\put(-85,104){\line(-1,1){23}}
\put(-85,104){\line(-4,1){56}}

\put(-105,132){\makebox(2,2){\circle*{11}}}
\put(-77,100){\makebox(2,2){\circle*{11}}}
                          \put(-140,110){\line(1,2){100}} 
 \put(-85,101){\line(1,2){100}} 
\put(-95,140){\makebox(30,10){{\bf aRM**} }}
\put(-80,80){\makebox(30,10){{\bf *aRM**} }}
\put(-155,140){\makebox(30,10){ (**)}}   \put(-115,140){\makebox(30,-80){ (**)}}
                            \put(-22,360){\makebox(10,10){{\it Hierarchy 3.2} }}

\put(-72,351){\makebox(2,2){$\bigcirc$}}  \put(-72,360){\makebox(10,10){{\bf RM} }}
\put(-8,332){\makebox(2,2){$\bigcirc$}} \put(-8,332){\makebox(2,2){$\bullet$}}
\put(-65,300){\makebox(10,10){{\bf *RM} }}

\put(-145,260){\makebox(30,10){ (An)}}
\put(-50,190){\makebox(30,10){ (An)}}
                  \put(-60,240){\makebox(30,10){ (An)}}
\put(-120,175){\makebox(10,10){(An)}}

\put(-60,320){\makebox(10,10){(*)}}  \put(10,312){\makebox(10,10){(*)}}

\put(-135,113){\makebox(2,2){\circle*{11}}}
 \put(-160,120){\makebox(10,10){(*)}}  \put(-90,112){\makebox(10,10){(*)}}

\put(-170,90){\makebox(30,10){{\bf *aRM} }}

\put(-176,147){\makebox(2,2){$\bigcirc$}}  \put(-176,147){\makebox(2,2){$\circ$}} 
\put(-195,158){\makebox(30,10){{\bf aRM}}}

\end{picture}
 \end{center}
\caption{ Hierarchies 3.2 and 2.2,  determined by RM and by RML algebras, considering  (*) and (**) properties   }
\label{fig:figrml***}
\end{figure}
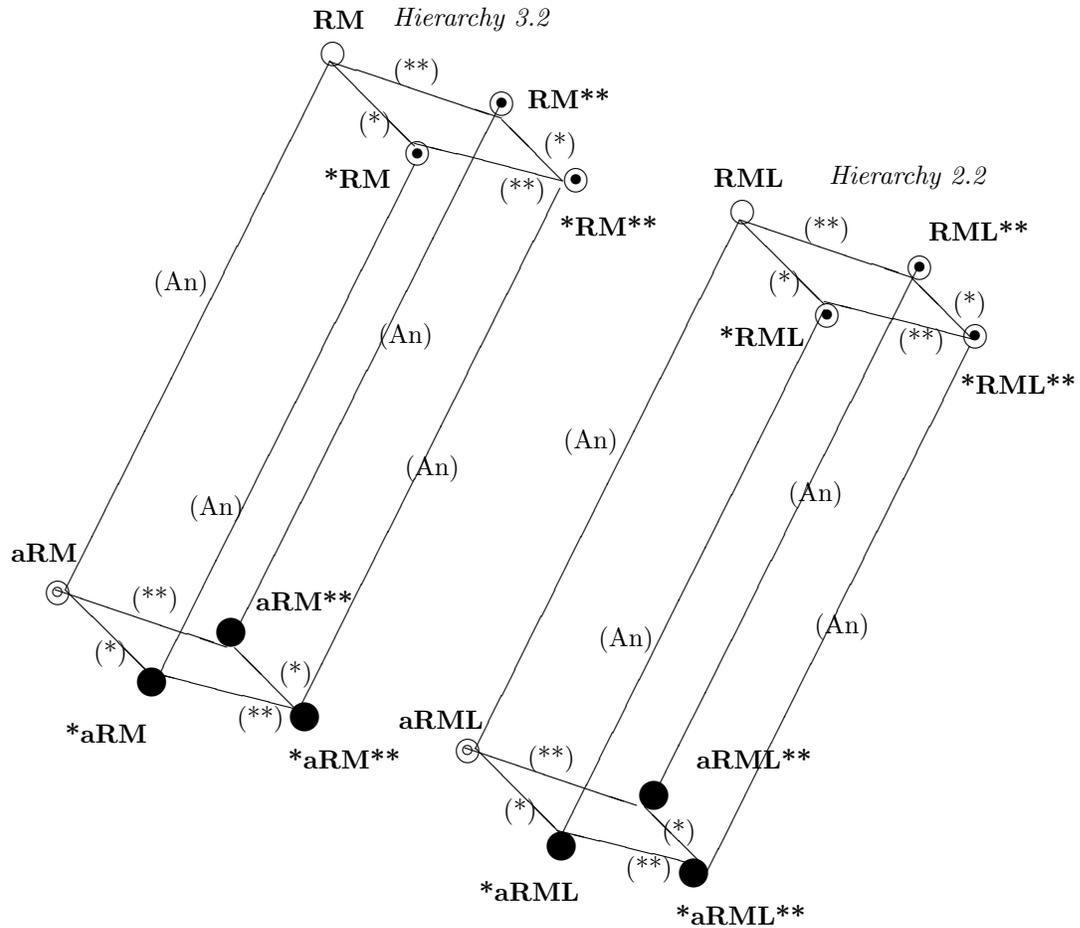

\begin{rem}\em J\" anis Cirulis introduced in 2006 the class of weak BCK-algebras.
 They were studied in \cite{Cirulis-1}
and as algebras with subtraction in \cite{Cirulis-2}, \cite{Cirulis-3}. We just noticed that weak BCK-algebras are in fact 
our aRML** algebras with property (D). 
\end{rem}

\begin{rem}\em Note that:\\
{\bf tRM} + (*) = {\bf *RM},\\
{\bf *RM} + (B) = {\bf pre-BZ},\\
{\bf RM**} + (B) = {\bf pre-BZ},\\
{\bf *RM**} + (B) = {\bf pre-BZ};\\
{\bf oRM} + (*) = {\bf *aRM},\\
{\bf *aRM} + (B) = {\bf BZ},\\
{\bf aRM**} + (B) = {\bf BZ},\\
{\bf *aRM**} + (B) = {\bf BZ};\\\\
{\bf tRML} + (*) = {\bf *RML},\\
{\bf *RML} + (B) = {\bf pre-BCC},\\
{\bf RML**} + (B) = {\bf pre-BCC},\\
{\bf *RML**} + (B) = {\bf pre-BCC};\\
{\bf oRML} + (*) = {\bf *aRML},\\
{\bf *aRML} + (B) = {\bf BCC},\\
{\bf aRML**} + (B) = {\bf BCC},\\
{\bf *aRML**} + (B) = {\bf BCC}.
\end{rem}

\begin{op}\em
Find other connections between the above RM algebras.
\end{op}
            \subsection{New RM and RML algebras, with property (Ex)}

By Theorem \ref{th1}, in RM and RML algebras satisfying condition (Ex),  we have:
$$ (B)  \; \Leftrightarrow \;  (BB) \; \Leftrightarrow \;  (*).$$
By Theorem \ref{th2}, in RM and RML algebras satisfying condition (Ex), we also have that:
$$ (**)  \; \Leftrightarrow \;  (Tr).$$
 
Hence, we have obviously the following:
\begin{cor}\label{m5}
In RM and RML algebras satisfying property (Ex), i.e. in RME, BCH and BE, aBE algebras, we have:
$$ (**) \;  \Leftrightarrow \; (Tr).$$
\end{cor}

Define now the following  last four new  algebras:
\begin{df}\em ${}$

28. A   {\it RME** algebra}   is a RME algebra verifying (**).

29. A {\it BCH** algebra} is a BCH algebra verifying (**).

30.  A {\it BE** algebra} is a BE algebra verifying (**).

31. An {\it aBE** algebra} is an  aBE algebra verifying (**).
\end{df}

Note that Corollary \ref{m5} gives us equivalent definitions for the four algebras RM**, BCH** and BE**, aBE**.

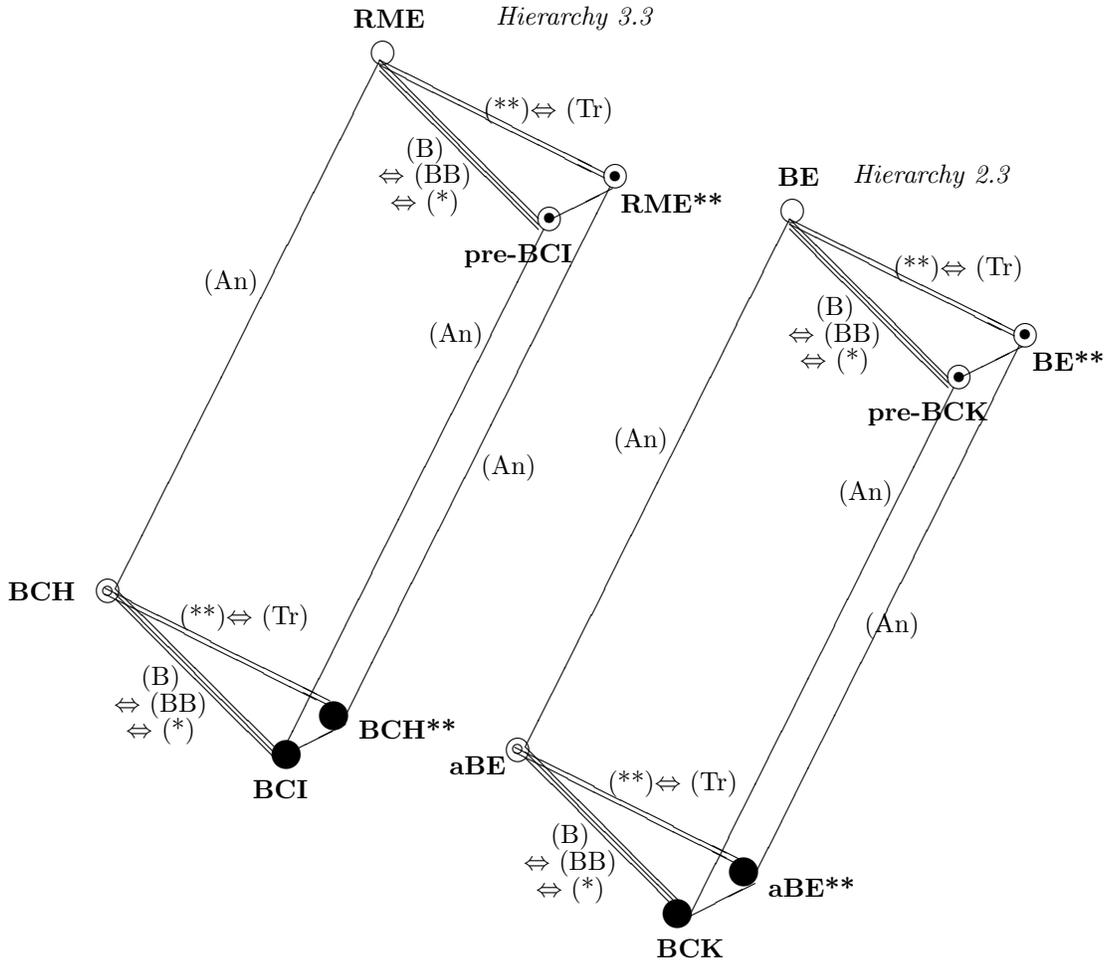
\begin{figure}[htbp]
\begin{center}
\begin{picture}(420,340)(-180,25)

       \put(-17,89){\line(1,-1){60}}    
       \put(-17,87){\line(1,-1){60}}    
       \put(-17,85){\line(1,-1){60}}

  \put(45,25){\line(1,2){100}}                               
  \put(-17,89){\line(1,2){100}}

      \put(83,289){\line(1,-1){60}}   
      \put(83,287){\line(1,-1){60}}   
       \put(83,285){\line(1,-1){60}}   
\put(83,289){\line(2,-1){85}} 
\put(83,287){\line(2,-1){85}} 

       \put(146,228){\line(2,1){25}}  
\put(171,244){\makebox(2,2){$\bigcirc$}} \put(171,244){\makebox(2,2){$\bullet$}}\put(185,230){\makebox(10,10){{\bf BE**} }}

\put(140,265){\makebox(10,10){ (**)$\Leftrightarrow$ (Tr)}}

\put(-21,89){\line(2,-1){85}} 
\put(-21,87){\line(2,-1){85}}   
  \put(45,25){\line(2,1){25}}    
\put(70,41){\makebox(2,2){\circle*{11}}}
                  
                    \put(70,41){\line(1,2){100}} 
\put(78,31){\makebox(30,10){{\bf aBE**} }}
\put(32,70){\makebox(10,10){ (**)$\Leftrightarrow$ (Tr)}}

                                    \put(133,300){\makebox(10,10){{\it Hierarchy 2.3} }}

\put(83,291){\makebox(2,2){$\bigcirc$}}  \put(83,300){\makebox(10,10){{\bf BE} }}
\put(146,228){\makebox(2,2){$\bigcirc$}} \put(146,228){\makebox(2,2){$\bullet$}}\put(132,210){\makebox(10,10){{\bf pre-BCK} }}

\put(10,200){\makebox(30,10){ (An)}}
\put(105,130){\makebox(30,10){ (An)}}
                  \put(95,180){\makebox(30,10){ (An)}}

\put(95,250){\makebox(10,10){(B)}}
\put(95,240){\makebox(10,10){$\Leftrightarrow$ (BB)}}
\put(95,230){\makebox(10,10){$\Leftrightarrow$ (*)}}

\put(45,25){\makebox(2,2){\circle*{11}}}
 \put(-5,50){\makebox(10,10){(B)}}
\put(-5,40){\makebox(10,10){$\Leftrightarrow$ (BB)}}
\put(-5,30){\makebox(10,10){$\Leftrightarrow$ (*)}}

\put(32,8){\makebox(30,10){{\bf BCK} }}

\put(-21,87){\makebox(2,2){$\bigcirc$}}  \put(-21,87){\makebox(2,2){$\circ$}} 
\put(-50,77){\makebox(30,10){{\bf aBE}}}
 \put(-172,149){\line(1,-1){60}}    
       \put(-172,147){\line(1,-1){60}}    
       \put(-172,145){\line(1,-1){60}}

  \put(-110,85){\line(1,2){100}}                               
\put(-172,149){\line(1,2){100}}

      \put(-72,349){\line(1,-1){60}}   
      \put(-72,347){\line(1,-1){60}}   
      \put(-72,345){\line(1,-1){60}}   
\put(-72,349){\line(2,-1){85}}
\put(-72,347){\line(2,-1){85}}   
   
  \put(-9,288){\line(2,1){25}}  
\put(16,304){\makebox(2,2){$\bigcirc$}} \put(16,304){\makebox(2,2){$\bullet$}}\put(35,290){\makebox(10,10){{\bf RME**} }}
 \put(-15,325){\makebox(10,10){ (**)$\Leftrightarrow$ (Tr)}}

\put(-176,149){\line(2,-1){85}}  
\put(-176,147){\line(2,-1){85}}  
 
  \put(-110,85){\line(2,1){25}}    
\put(-85,100){\makebox(2,2){\circle*{11}}}
 
 \put(-85,101){\line(1,2){100}} 
\put(-75,91){\makebox(30,10){{\bf BCH**} }}
\put(-130,133){\makebox(10,10){ (**)$\Leftrightarrow$ (Tr)}}

                                    \put(-2,360){\makebox(10,10){{\it Hierarchy 3.3} }}

\put(-72,351){\makebox(2,2){$\bigcirc$}}  \put(-72,360){\makebox(10,10){{\bf RME} }}
\put(-9,288){\makebox(2,2){$\bigcirc$}} \put(-9,288){\makebox(2,2){$\bullet$}}\put(-23,270){\makebox(10,10){{\bf pre-BCI} }}

\put(-145,260){\makebox(30,10){ (An)}}
\put(-40,190){\makebox(30,10){ (An)}}
                 \put(-60,240){\makebox(30,10){ (An)}}

\put(-60,310){\makebox(10,10){(B)}}
\put(-60,300){\makebox(10,10){$\Leftrightarrow$ (BB)}}
\put(-60,290){\makebox(10,10){$\Leftrightarrow$ (*)}}

\put(-103,85){\makebox(2,2){\circle*{11}}}
 \put(-160,110){\makebox(10,10){(B)}}
\put(-160,100){\makebox(10,10){$\Leftrightarrow$ (BB)}}
\put(-160,90){\makebox(10,10){$\Leftrightarrow$ (*)}}

\put(-123,68){\makebox(30,10){{\bf BCI} }}

\put(-176,147){\makebox(2,2){$\bigcirc$}}  \put(-176,147){\makebox(2,2){$\circ$}} 
\put(-215,143){\makebox(30,10){{\bf BCH}}}

\end{picture}
 \end{center}
\caption{ Hierarchies 3.3 and 2.3, determined by RME and  BE algebras   }
\label{fig:fig23ssT}
\end{figure}

It follows that we have, for example,  the Hierarchies 3.3 and 2.3 
 from Figure \ref{fig:fig23ssT}; in order to simplify the figure,
 the connection by (L) between Hierarchies 3.3 and 2.3 is not drawn.

\begin{rem}\em Note that:\\
{\bf RME**} + (B) = {\bf pre-BCI},\\
{\bf BCH**} + (B) = {\bf BCI},\\
{\bf BE**} + (B) = {\bf pre-BCK},\\
{\bf aBE**} + (B) = {\bf BCK}.
\end{rem}
              \section{Generalizations of Hilbert algebras}

\begin{df}\label{Hilberta}\em (see \cite{Diego})

 A {\it Hilbert algebra} is an algebra  ${\cal A}=(A, \ra,1)$ of type $(2,0)$, satisfying, 
for all $x,y,z \in A$:\\
(h1) $x \ra (y \ra x)=1$,\\
(h2) $(x \ra (y \ra z)) \ra ((x \ra y) \ra (x \ra z))=1,$\\
(h3) if $x \ra y=y \ra x=1$, then $x=y$.
\end{df}

Note that (h1) is (K) and (h3) is (An).

Let $(A,\ra,1)$ be an algebra of type $(2,0)$.
Consider the properties: for all $x,y,z \in A$\\
(pimpl) $x \ra (y \ra z) = (x \ra y) \ra (x \ra z)$,\\
(pi) $y \ra (y \ra x) = y \ra x,$  \\
 and also\\
(p-1) $[x \ra (y \ra z)] \ra [  (x \ra y) \ra (x \ra z)]=1$,\\
(p-2) $[(x \ra y) \ra (x \ra z)] \ra [x \ra (y \ra z)]=1 $.\\

Note that (p-1) is (h2).

We know \cite{Iseki-Tanaka} that in BCK algebras we have the equivalence:
$$ (pimpl) \; \Longleftrightarrow \;  (pi)$$
and that \cite{book}  Hilbert algebras are just those BCK algebras, called {\it positive implicative}, 
 verifying  the equivalent properties  (pimpl) and (pi).

We shall see in this section  what happens in the old and new generalizations 
of BCK algebras with these two properties (pimpl) and (pi).

\begin{prop}\label{prop-pi} Let $(A, \ra, 1)$ be an algebra of type $(2,0)$. Then,

(i) (Re) + (pimpl) $\; \Longrightarrow \; $ (L); 

(ii) (Re) + (pi) $\; \Longrightarrow \; $ (L); 

(iii) (Re) +(M) + (pimpl) $\; \Longrightarrow \; $ (pi).
\end{prop}
{\bf Proof.}

(i): Take $y=z=x$ in (pimpl); we obtain: $x \ra (x \ra x)=(x \ra x) \ra (x \ra x)$, hence, by (Re), we obtain:
$x \ra 1=1 \ra 1$, hence, by (Re) again, $x \ra 1=1$, i.e. (L) holds.

(ii): Take $y=x$ in (pi); we obtain: $x \ra (x \ra x)=x \ra x$, hence, by (Re), $x \ra 1=1$, i.e. (L) holds.

(iii): Take $x=y$ in (pimpl); we obtain:  $y \ra (y \ra z) = (y \ra y) \ra (y \ra z);$ then, by (Re), we obtain:
 $y \ra (y \ra z)=1 \ra (y \ra z)$; then, by (M), we obtain: $y \ra (y \ra z)=y \ra z$, i.e. (pi) holds.
\hfill $\Box$\\
 
We shall see now in which properties  (pi) implies  (pimpl). 
\begin{prop} Let $(A, \ra,1)$ be an algebra of type $(2,0)$.  Then,

 (a) (Re) +  (L) + (Ex) + (**) $\; \Longrightarrow \; $ (p-2);

 (b)  (Ex) + (B) + (*) + (pi)  $\; \Longrightarrow \; $ (p-1);

 (c)  (p-1) + (p-2) + (An) $\; \Longrightarrow \; $ (pimpl);

(d) (Re) + (Ex) + (B) + (**) + (*) + (L) + (An) + (pi) $\; \Longrightarrow \; $ (pimpl).

\end{prop}

{\bf Proof.} 

 (a):  By Proposition \ref{propp}, (8), (Re) + (Ex) + (L) imply (K). By (K), we have: $y \ra ( x \ra y)=1$, hence, by 
(**), we obtain: $$[(x \ra y) \ra (x \ra z)] \ra [ y \ra (x \ra z)]=1;$$ hence, by (Ex), we obtain:
$$[(x \ra y) \ra (x \ra z)] \ra [x \ra (y \ra z)]=1,$$
i.e. (p-2) holds.

 (b): We must prove that $[x \ra (y \ra z)] \ra [  (x \ra y) \ra (x \ra z)]=1$.\\
Denote $X=(x \ra y) \ra (x \ra z)$. We obtain: $$X \stackrel{(pi)}{=} (x \ra y) \ra [x \ra (x \ra z)] \stackrel{(Ex)}{=}
x \ra [(x \ra y) \ra (x \ra z)].$$
By (B), we have:  $(y \ra z) \ra [(x \ra y) \ra (x \ra z)]=1$; then, by (*), we obtain:
$$[x \ra (y \ra z)] \ra [x \ra [(x \ra y) \ra (x \ra z)]]=1= [x \ra (y \ra z)] \ra X,$$ 
hence $$[x \ra (y \ra z)] \ra [(x \ra y) \ra (x \ra z)]=1,$$
i.e. (p-1) holds.

 (c): Obviously.

(d): By above (a), (b), (p-2) and (p-1) hold. Then, by (c), (pimpl) holds.

\hfill $\Box$

We then immediately obtain:
\begin{cor}${}$

(1) Properties (pimpl) and (pi) may hold only in the generalizations of BCK algebras verifying (L), i.e. in the 
RML algebras.

(2) In  RML algebras, (pimpl) implies (pi).

(3) BE**, aBE** and pre-BCK, BCK algebras verify (p-2).

(4) In pre-BCK and BCK  algebras, (pi) implies (p-1).

(5) In BCK algebras only, (pi) implies (pimpl), hence (pimpl) $\Leftrightarrow$ (pi).

 (6) Hilbert algebras are the only RML algebras where (pimpl) and (pi) hold and  we have (pimpl) $\Leftrightarrow$ (pi). 

\end{cor}

We prove  now  the following  result:
\begin{prop}\label{pkt}

$$ (pimpl) \; + \; (K) \; \Longrightarrow \; (B).$$
\end{prop}
{\bf Proof.} $y \ra z \stackrel{(K')}{\leq} x \ra (y \ra z) \stackrel{(pimpl)}{=} (x \ra y) \ra (x \ra z)$.
Hence,  $y \ra z \leq (x \ra y) \ra (x \ra z)$, i.e. (B') holds.
\hfill $\Box$

\begin{cor}
In proper BE, aBE, BE**, aBE** algebras, property (pimpl) cannot hold.
\end{cor}
{\bf Proof.}
Since  by Proposition \ref{propp} (8.), (Re) +  (L) + (Ex)  imply (K), it follows that in BE algebras and in  aBE algebras we have (K), hence in
 BE** and aBE** algebras  we have property (K). In proper BE, aBE, BE**, aBE** algebras, property (B) does not hold.
  If (pimpl) holds in proper BE, aBE, BE**, aBE** algebras, then, by above Proposition \ref{pkt}, (B) holds: contradiction.
\hfill $\Box$

\begin{cor}
In proper RML**, *RML**, aRML**, *aRML** algebras, property (pimpl) cannot hold.
\end{cor}
{\bf Proof.}
Since  by Proposition \ref{propp} (9'.), (M) +  (L) + (**)  imply (K), it follows that in RML**, *RML**, aRML**, *aRML**
 algebras 
 we have (K). In proper RML**, *RML**, aRML**, *aRML** algebras, property (B) does not hold.
  If (pimpl) holds in proper RML**, *RML**, aRML**, *aRML** algebras, then, by above Proposition \ref{pkt}, (B) holds: contradiction.
\hfill $\Box$

\begin{rems}\label{m6}\em ${}$

(1) There are plenty of RML algebras verifying (pi) and not verifying (pimpl).

(2) The examples we found of RML algebras verifying (pimpl) (hence (pi)) are pre-BCK, BCK  and pre-BBBCC algebras. 

(3) We know that Hilbert algebras cannot be generalized to the non-commutative case. 
In  non-commutative RML algebras $(A, \ra, \ca, 1)$ (as for example the pseudo-BCK algebras, the pseudo-BE algebras), 
property (pi) would become a pair of the properties:\\\\
 (pi$\ca$) $y \ra (y \ca x)=y \ca x$,\\
  (pi$\ra$) $y \ca (y \ra x)=y \ra x$,\\\\
and  since in the non-commutative case, (Ex) becomes:\\
 (pEx) $x \ra (y \ca z)=y \ca (x \ra z)$\\
it follows that, for all $x,y$:
 $$y \ra x \stackrel{(pi\ra)}{=}y \ca (y \ra x) \stackrel{(pEx)}{=} y \ra (y \ca x)\stackrel{(pi\ca)}{=}y \ca x,$$
i.e. $\ra=\ca$, hence this non-commutative algebra $(A, \ra, \ca, 1)$ is in fact a commutative one. In other words,
the RML algebras verifying  (Ex)  and (pi) (i.e. the BE, pre-BCK, aBE, BCK  and BE**, aBE** algebras verifying (pi)) 
cannot be generalized to the non-commutative case.
\end{rems}

Let us  then introduce  the following definitions:          
\begin{df}\em ${}$

(1) A  {\it pi-RML algebra} is a RML algebra  verifying property  (pi). 

(2) A  {\it positive implicative RML algebra}, or  a {\it pimpl-RML algebra} for short, 
 is a RML algebra verifying property  (pimpl).

(3)  A {\it generalization of Hilbert algebra} is any pi-RML algebra.

(4) A {\it proper generalization  of Hilbert algebra} 
is a pi-RML algebra that cannot be generalized to the non-commutative case.
\end{df}

\begin{rems} \em ${}$

(1) Any pimpl-RML algebra is pi-RML algebra; the converse holds only for BCK algebras.

(2) The pi-RML algebras and the pimpl-RML algebras,  their hierarchies and 
their connections with the corresponding RML algebras,
are presented in the next  Figures \ref{fig:fig1234H}, \ref{fig:fig23sTH}, \ref{fig:figrml***H}
 and \ref{fig:fig23ssTH}. We shall present in other section 
 examples of all of these generalizations of Hilbert algebras.

(3) By above Remarks \ref{m6} (3), the proper generalizations of Hilbert algebras are the  following pi-RML algebras: 
pi-BE, pi-pre-BCK,pimpl-pre-BCK,  pi-aBE, pi-BCK = pimpl-BCK = Hilbert and   pi-BE**, pi-aBE**.
 \end{rems}

\begin{op}\em As Hilbert algebra  (see Definition  \ref{Hilberta}) is  equivalent 
with  pimpl-BCK algebra (= pi-BCK algebra), 
define a {\it pre-Hilbert algebra} as an algebra $(A, \ra,1)$ equivalent to pimpl-pre-BCK algebra. 
\end{op}

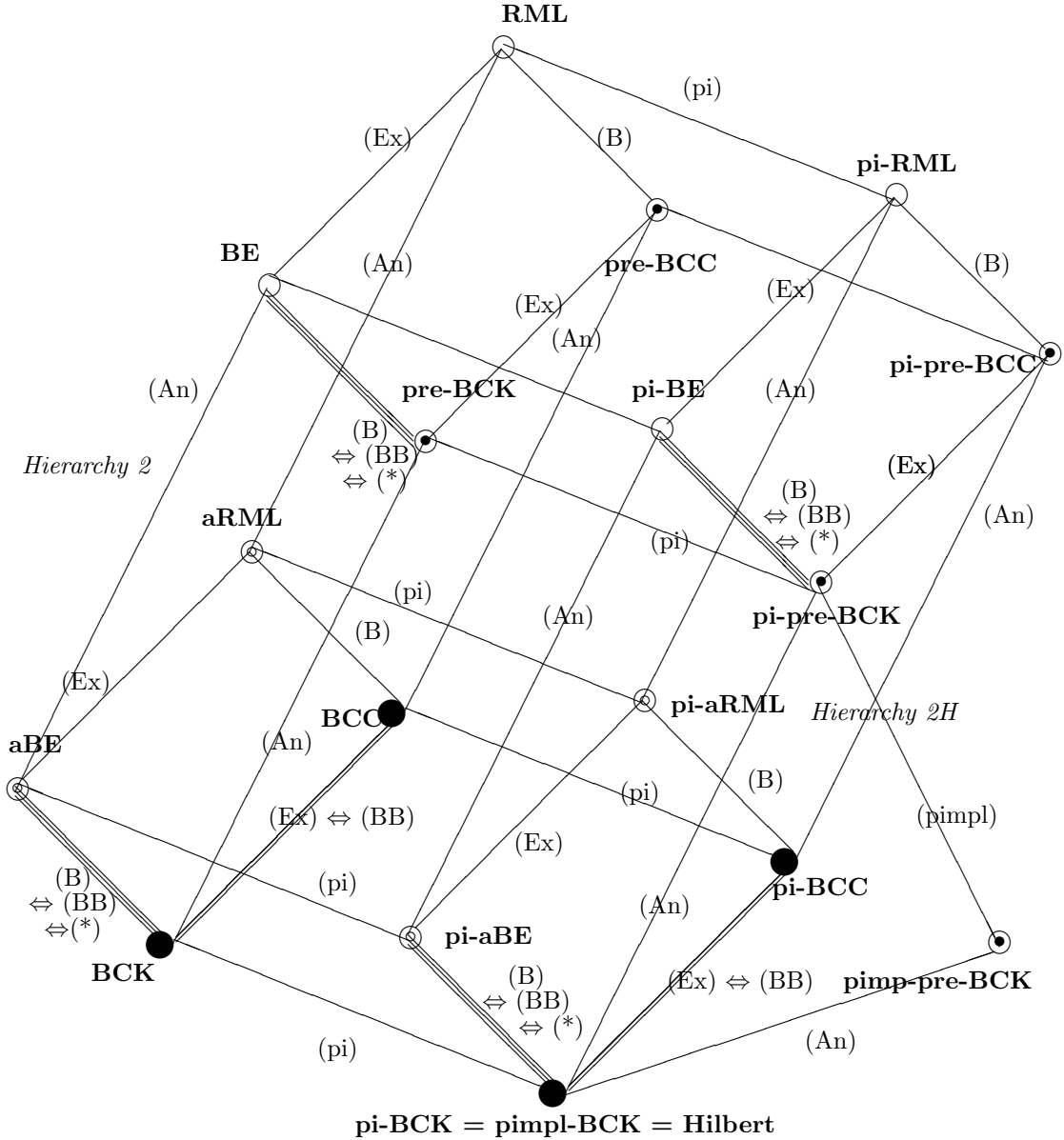
\begin{figure}[htbp]
\begin{center}
\begin{picture}(420,420)(-180,25) 
\put(155,172){\makebox(30,10){ {\it  Hierarchy 2H}}}
\put(-160,270){\makebox(30,10){{\it Hierarchy 2}}}

\put(76,182){\line(-1,-1){90}} 
   \put(137,119){\line(-1,-1){90}}    
     \put(76,182){\line(1,-1){60}}  
           \put(-17,89){\line(1,-1){60}}    
           \put(-17,87){\line(1,-1){60}}    
           \put(-17,85){\line(1,-1){60}}    

  \put(76,181){\line(-5,2){155}}
  \put(45,25){\line(-5,2){155}}
 \put(-17,87){\line(-5,2){155}}
  \put(137,117){\line(-5,2){155}}
  \put(45,25){\line(1,2){100}} 
\put(137,119){\line(1,2){100}} 
\put(-17,89){\line(1,2){100}} 
\put(76,182){\line(1,2){100}} 

\put(176,382){\line(-1,-1){90}} 
  \put(237,319){\line(-1,-1){90}}  
\put(176,382){\line(1,-1){60}}  

\put(83,289){\line(1,-1){59}}   
\put(83,287){\line(1,-1){59}}   
\put(83,285){\line(1,-1){59}}

\put(83,289){\makebox(2,2){$\bigcirc$}}  \put(83,300){\makebox(10,10){{\bf pi-BE} }}
\put(176,382){\makebox(2,2){$\bigcirc$}} \put(176,390){\makebox(10,10){{\bf pi-RML}}}
\put(237,319){\makebox(2,2){$\bigcirc$}} \put(237,319){\makebox(2,2){$\bullet$}} \put(200,310){\makebox(10,10){{\bf pi-pre-BCC} }}
\put(146,228){\makebox(2,2){$\bigcirc$}}  \put(146,228){\makebox(2,2){$\bullet$}} 
             \put(146,210){\makebox(10,10){{\bf pi-pre-BCK} }}
                                 \put(146,228){\line(1,-2){72}} 
                                \put(45,25){\line(3,1){171}}   
                                 \put(217,85){\makebox(2,2){$\bigcirc$}}
                                  \put(217,85){\makebox(2,2){$\bullet$}} 
                                  \put(190,65){\makebox(10,10){{\bf pimp-pre-BCK} }}
                                  \put(150,45){\makebox(2,2){(An)}} 
                                   \put(200,135){\makebox(2,2){(pimpl)}} 

\put(30,210){\makebox(30,10){ (An)}}
\put(205,250){\makebox(30,10){ (An)}}
\put(80,95){\makebox(10,10){(An)}}

\put(130,300){\makebox(10,10){(An)}}

\put(210,350){\makebox(10,10){(B)}}
\put(130,340){\makebox(10,10){(Ex)}}

\put(82,240){\makebox(10,10){(pi)}}

\put(135,260){\makebox(10,10){(B) }}
\put(135,250){\makebox(10,10){ $\Leftrightarrow$ (BB)}}
\put(135,240){\makebox(10,10){ $\Leftrightarrow$  (*)}}

\put(178,270){\makebox(10,10){(Ex)}}
 \put(-173,146){\line(1,2){100}} 
\put(-80,241){\line(1,2){100}} 
\put(-111,84){\line(1,2){100}} 
\put(-19,176){\line(1,2){100}} 

\put(20,441){\line(-1,-1){90}} 
  \put(81,376){\line(-1,-1){90}}  
\put(20,441){\line(1,-1){60}}  
         \put(-73,343){\line(1,-1){58}}   
         \put(-73,341){\line(1,-1){58}}   
         \put(-73,345){\line(1,-1){58}}

\put(-73,346){\makebox(2,2){$\bigcirc$}}  \put(-87,355){\makebox(10,10){{\bf BE} }}
\put(20,441){\makebox(2,2){$\bigcirc$}} \put(30,450){\makebox(10,10){{\bf RML} }}
\put(-11,284){\makebox(2,2){$\bigcirc$}} \put(-11,284){\makebox(2,2){$\bullet$}} \put(0,300){\makebox(10,10){{\bf pre-BCK} }}
\put(81,376){\makebox(2,2){$\bigcirc$}}  \put(81,376){\makebox(2,2){$\bullet$}} \put(80,350){\makebox(10,10){{\bf pre-BCC} }}

\put(60,400){\makebox(10,10){(B)}}
 \put(30,334){\makebox(10,10){(Ex)}}
\put(-30,400){\makebox(10,10){(Ex)}}  
\put(-30,350){\makebox(10,10){(An)}}
\put(-115,300){\makebox(10,10){(An)}}
\put(45,320){\makebox(10,10){(An)}} 

\put(-37,284){\makebox(10,10){(B)}}
\put(-35,274){\makebox(10,10){$\Leftrightarrow$ (BB)}}
\put(-35,264){\makebox(10,10){$\Leftrightarrow$ (*)}}

\put(178,270){\makebox(10,10){(Ex)}}
 \put(176,381){\line(-5,2){155}}
  \put(145,225){\line(-5,2){155}}
 \put(83,289){\line(-5,2){155}} 
  \put(237,317){\line(-5,2){155}}

\put(85,420){\makebox(30,10){ (pi) }}

\put(45,25){\makebox(2,2){\circle*{11}}}

 \put(25,67){\makebox(10,10){(B)}}
\put(25,57){\makebox(10,10){$\Leftrightarrow$ (BB)}}
\put(35,47){\makebox(10,10){$\Leftrightarrow$ (*)}}


\put(32,8){\makebox(30,10){{\bf pi-BCK = pimpl-BCK = Hilbert}  }}

\put(-17,87){\makebox(2,2){$\bigcirc$}}  \put(-17,87){\makebox(2,2){$\circ$}} 

\put(30,120){\makebox(10,10){(Ex)}}

\put(0,83){\makebox(30,10){{\bf pi-aBE}}}

\put(76,181){\makebox(2,2){$\bigcirc$}} \put(76,181){\makebox(2,2){$\circ$}}
 
\put(95,175){\makebox(30,10){ {\bf  pi-aRML} }}
\put(137,117){\makebox(2,2){\circle*{11}}} 
\put(120,145){\makebox(10,10){(B)}}
                                          
\put(130,102){\makebox(30,10){  {\bf pi-BCC}}}

\put(-80,241){\line(-1,-1){90}} 
   \put(-19,178){\line(-1,-1){90}}    
     \put(-80,241){\line(1,-1){60}} 
        \put(-173,148){\line(1,-1){60}}    
        \put(-173,146){\line(1,-1){60}}    
        \put(-173,144){\line(1,-1){60}}    

  
\put(-111,84){\makebox(2,2){\circle*{11}}} 

 \put(-155,106){\makebox(10,10){(B)}}
\put(-155,96){\makebox(10,10){$\Leftrightarrow$ (BB)}}
\put(-155,86){\makebox(10,10){$\Leftrightarrow$(*)}}

\put(-150,69){ \makebox(30,10){ {\bf BCK}}}
\put(-173,146){\makebox(2,2){$\bigcirc$}} \put(-173,146){\makebox(2,2){$\circ$}}

 \put(-150,184){\makebox(10,10){(Ex)}}

\put(-180,160){\makebox(30,10){{\bf aBE}}}
 
\put(-80,240){\makebox(2,2){$\bigcirc$}} \put(-80,240){\makebox(2,2){$\circ$}}

\put(-96,250){\makebox(30,10){{\bf aRML }}}

\put(-80,160){\makebox(30,10){ (An) }} 

\put(-30,220){\makebox(30,10){ (pi) }}
\put(-60,104){\makebox(30,10){(pi)}}
\put(60,140){\makebox(30,10){(pi)}}
\put(-60,38){\makebox(30,10){(pi)}}

\put(-19,176){\makebox(2,2){\circle*{11}}}
 \put(-36,204){\makebox(10,10){(B)}}
\put(-56,170){\makebox(30,10){ {\bf BCC}}}                                         
    \put(-19,178){\line(-1,-1){90}}
                                                 \put(-19,176){\line(-1,-1){90}} 
                                                \put(137,119){\line(-1,-1){90}} 
                                             \put(137,117){\line(-1,-1){90}}                                      
                                              \put(-60,130){\makebox(30,10){ (Ex) $\Leftrightarrow$  (BB)}}
                                              \put(110,65){\makebox(10,10){(Ex) $\Leftrightarrow $ (BB)}}

\end{picture}
 \end{center}
\caption{ The  hierarchy connecting Hierarchies 2 and 2H }
\label{fig:fig1234H}
\end{figure}

\begin{figure}[htbp]
\begin{center}
\begin{picture}(420,400)(-180,-25)

 \put(-17,89){\line(1,-1){60}}   
                                                      

  \put(45,25){\line(1,2){100}}                               
\put(-17,89){\line(1,2){100}}

      \put(83,289){\line(1,-1){97}}   
\put(83,289){\line(2,-1){85}}   
\put(83,289){\line(3,-1){102}}   
  \put(146,228){\line(2,1){46}}  
\put(171,244){\makebox(2,2){$\bigcirc$}} \put(171,244){\makebox(2,2){$\bullet$}}\put(185,230){\makebox(10,10){{\bf pi-*RML} }}
\put(190,254){\makebox(2,2){$\bigcirc$}} \put(190,254){\makebox(2,2){$\bullet$}}\put(183,265){\makebox(10,10){{\bf pi-tRML} }}
\put(139,270){\makebox(30,10){ (Tr)}}
\put(125,245){\makebox(30,10){ (*)}}
\put(-21,89){\line(2,-1){85}}   
\put(-21,89){\line(3,-1){102}}   
   \put(45,25){\line(2,1){46}}    
\put(90,50){\makebox(2,2){\circle*{11}}}
\put(70,41){\makebox(2,2){\circle*{11}}}
 \put(90,50){\line(1,2){100}} 
                                                                \put(182,187){\line(-1,-2){100}} 
                                                                 \put(27,248){\line(-1,-2){100}}

 \put(70,41){\line(1,2){100}} 
\put(105,40){\makebox(30,10){{\bf pi-oRML} }}
\put(70,21){\makebox(30,10){{\bf pi-*aRML} }}
\put(35,69){\makebox(30,10){ (Tr)}}
\put(26,44){\makebox(30,10){ (*)}}

\put(83,291){\makebox(2,2){$\bigcirc$}}  \put(83,300){\makebox(10,10){{\bf pi-RML} }}
\put(146,228){\makebox(2,2){$\bigcirc$}} \put(146,228){\makebox(2,2){$\bullet$}} 
                              \put(110,220){\makebox(10,10){{\bf pi-pre-BCC} }}
\put(182,187){\makebox(2,2){$\bigcirc$}}
\put(182,187){\makebox(2,2){$\bullet$}}
\put(170,167){\makebox(10,10){{\bf pi-pre-BBBCC}}}
\put(152,197){\makebox(10,10){(BB)}}
\put(182,187){\line(1,-3){26}} 
\put(206,107){\makebox(2,2){$\bigcirc$}}
\put(206,107){\makebox(2,2){$\bullet$}}
\put(160,85){\makebox(10,10){{\bf pimpl-pre-BBBCC}}}
\put(180,137){\makebox(10,10){ (pimpl)}}

\put(10,200){\makebox(30,10){ (An)}}
\put(105,130){\makebox(30,10){ (An)}}
                  \put(95,180){\makebox(30,10){ (An)}}
\put(134,120){\makebox(10,10){(An)}}

\put(95,250){\makebox(10,10){(B)}}

\put(49,25){\makebox(2,2){\circle*{11}}}
 \put(-5,50){\makebox(10,10){(B)}}


\put(0,18){\makebox(30,10){{\bf pi-BCC} }}

\put(-21,87){\makebox(2,2){$\bigcirc$}}  \put(-21,87){\makebox(2,2){$\circ$}} 
\put(-15,99){\makebox(30,10){{\bf pi-aRML}}}
 \put(-172,149){\line(1,-1){60}}    

  \put(-110,85){\line(1,2){100}}                               
\put(-172,149){\line(1,2){100}}

      \put(-72,349){\line(1,-1){97}}   

\put(-72,349){\line(2,-1){85}}   
\put(-72,349){\line(3,-1){102}}   
  \put(-9,288){\line(2,1){46}}  
\put(16,304){\makebox(2,2){$\bigcirc$}} \put(16,304){\makebox(2,2){$\bullet$}}\put(30,290){\makebox(10,10){{\bf *RML} }}
\put(35,314){\makebox(2,2){$\bigcirc$}} \put(35,314){\makebox(2,2){$\bullet$}}\put(55,315){\makebox(10,10){{\bf tRML} }}
\put(-16,330){\makebox(30,10){ (Tr)}}
\put(-30,305){\makebox(30,10){ (*)}}
\put(-176,149){\line(2,-1){85}}   
\put(-176,149){\line(3,-1){102}}   
  \put(-110,85){\line(2,1){46}}    
\put(-65,110){\makebox(2,2){\circle*{11}}}
\put(-85,100){\makebox(2,2){\circle*{11}}}
  \put(-65,110){\line(1,2){100}} 
 \put(-85,101){\line(1,2){100}} 
\put(-60,100){\makebox(30,10){{\bf oRML} }}
\put(-75,91){\makebox(30,10){{\bf *aRML} }}
\put(-120,129){\makebox(30,10){ (Tr)}}
\put(-129,104){\makebox(30,10){ (*)}}

\put(-72,351){\makebox(2,2){$\bigcirc$}}  \put(-72,360){\makebox(10,10){{\bf RML} }}
\put(-9,288){\makebox(2,2){$\bigcirc$}} \put(-9,288){\makebox(2,2){$\bullet$}}\put(-40,280){\makebox(10,10){{\bf pre-BCC} }}

\put(27,248){\makebox(2,2){$\bigcirc$}}
 \put(27,248){\makebox(2,2){$\bullet$}}
\put(22,228){\makebox(10,10){{\bf pre-BBBCC} }}
\put(-4,258){\makebox(10,10){ (BB) }}

\put(-145,260){\makebox(30,10){ (An)}}
\put(-40,190){\makebox(30,10){ (An)}}
                  \put(-60,240){\makebox(30,10){ (An)}}
\put(-30,155){\makebox(10,10){(An)}}

\put(-60,310){\makebox(10,10){(B)}}

 \put(-103,85){\makebox(2,2){\circle*{11}}}
 \put(-160,110){\makebox(10,10){(B)}}

\put(-141,75){\makebox(30,10){{\bf BCC} }}   
                                              \put(-70,47){\makebox(2,2){\circle*{11}}}
                                             \put(-91,28){\makebox(30,10){{\bf BCK} }}  
                                              \put(-131,55){\makebox(30,10){(Ex) $\Leftrightarrow$ (BB) }}    
                                              \put(-105,82){\line(1,-1){32}}                                                                                                 
                                              \put(-105,80){\line(1,-1){32}} 

                                               \put(47,24){\line(1,-1){34}} 
                                               \put(47,22){\line(1,-1){34}} 
                                               \put(87,-14){\makebox(2,2){\circle*{11}}}
                                                 \put(85,-35){\makebox(30,10){{\bf pi-BCK = pimpl-BCK = Hilbert} }}  
                                                 \put(23,-3){\makebox(30,10){(Ex) $\Leftrightarrow$ (BB) }}   
                                                 \put(85,-12){\line(1,1){120}} 
                                                 \put(152,65){\makebox(2,2){(An)}}

\put(-176,147){\makebox(2,2){$\bigcirc$}}  \put(-176,147){\makebox(2,2){$\circ$}} 
\put(-175,160){\makebox(30,10){{\bf aRML}}}
                          \put(148,315){\makebox(10,10){{\it Hierarchy 2.1 H} }}
                            \put(-18,370){\makebox(10,10){{\it Hierarchy 2.1} }}

\end{picture}
 \end{center}
\caption{ Hierarchies 2.1 and 2.1 H, determined by RML and by pi-RML algebras, considering  (B), (*) and (Tr) properties   }
\label{fig:fig23sTH}
\end{figure}
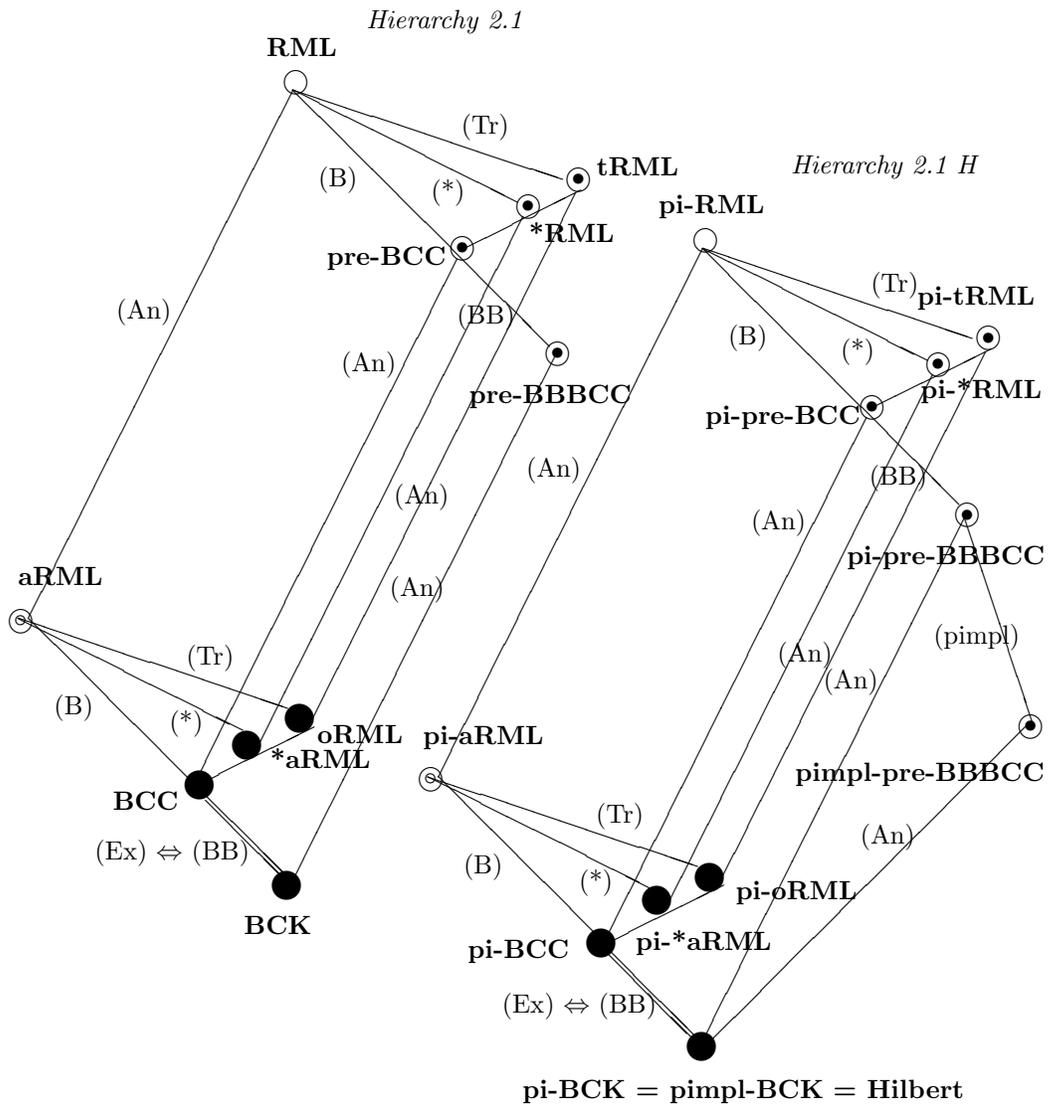


\begin{figure}[htbp]
\begin{center}
\begin{picture}(420,340)(-180,25)

       \put(-17,89){\line(1,-1){32}}    
                 
                      \put(15,55){\line(1,2){100}}                               
\put(-17,89){\line(1,2){100}}

      \put(83,289){\line(1,-1){32}}   
       
\put(83,289){\line(3,-1){65}}

\put(171,244){\line(-1,1){23}}
\put(171,244){\line(-4,1){56}}
\put(171,244){\makebox(2,2){$\bigcirc$}} \put(171,244){\makebox(2,2){$\bullet$}}\put(185,222){\makebox(10,10){{\bf pi-*RML**} (?) }}
\put(150,270){\makebox(2,2){$\bigcirc$}} \put(150,270){\makebox(2,2){$\bullet$}}\put(170,280){\makebox(10,10){{\bf pi-RML**} }}
\put(99,280){\makebox(30,10){ (**)}}  \put(135,238){\makebox(30,10){ (**)}}
\put(-21,89){\line(3,-1){65}}   
    
\put(70,44){\line(-1,1){23}}
\put(70,44){\line(-4,1){56}}
 
             \put(55,70){\makebox(2,2){\circle*{11}}}
\put(70,41){\makebox(2,2){\circle*{11}}}
                             \put(50,70){\line(1,2){100}} 
 \put(70,41){\line(1,2){100}} 
                \put(75,80){\makebox(30,10){{\bf pi-aRML**} }}
\put(70,21){\makebox(30,10){{\bf pi-*aRML**} }}
\put(-5,80){\makebox(30,10){ (**)}}  \put(32,38){\makebox(30,10){ (**)}}
\put(83,291){\makebox(2,2){$\bigcirc$}}  \put(83,300){\makebox(10,10){{\bf pi-RML} }}
\put(115,252){\makebox(2,2){$\bigcirc$}} \put(115,252){\makebox(2,2){$\bullet$}}\put(88,240){\makebox(10,10){{\bf pi-*RML} }}

\put(10,200){\makebox(30,10){ (An)}}
\put(105,130){\makebox(30,10){ (An)}}
                  \put(95,180){\makebox(30,10){ (An)}}
\put(35,125){\makebox(10,10){(An)}}

\put(95,260){\makebox(10,10){(*)}}   \put(165,252){\makebox(10,10){(*)}}

\put(20,51){\makebox(2,2){\circle*{11}}}
 \put(-5,60){\makebox(10,10){(*)}}  \put(55,52){\makebox(10,10){(*)}}

\put(-10,30){\makebox(30,10){{\bf pi-*aRML} }}

\put(-21,87){\makebox(2,2){$\bigcirc$}}  \put(-21,87){\makebox(2,2){$\circ$}} 
\put(-45,95){\makebox(30,10){{\bf pi-aRML}}}
 \put(-172,149){\line(1,-1){32}}    
     
                       \put(-108,132){\line(1,2){100}}                               
\put(-172,149){\line(1,2){100}}

      \put(-72,349){\line(1,-1){32}}   
      
\put(-72,349){\line(3,-1){65}}   
  
\put(16,304){\line(-1,1){23}}
\put(16,304){\line(-4,1){56}}

\put(-40,313){\makebox(2,2){$\bigcirc$}} \put(-40,313){\makebox(2,2){$\bullet$}}\put(30,282){\makebox(10,10){{\bf *RML**} }}
\put(20,303){\makebox(2,2){$\bigcirc$}} \put(20,303){\makebox(2,2){$\bullet$}}
\put(15,330){\makebox(10,10){{\bf RML**} }}
\put(-56,340){\makebox(30,10){ (**)}}   \put(-16,295){\makebox(30,10){ (**)}}
\put(-176,149){\line(3,-1){65}}   
   
\put(-85,104){\line(-1,1){23}}
\put(-85,104){\line(-4,1){56}}

\put(-105,132){\makebox(2,2){\circle*{11}}}
\put(-77,100){\makebox(2,2){\circle*{11}}}
                          \put(-140,110){\line(1,2){100}} 
 \put(-85,101){\line(1,2){100}} 
\put(-95,140){\makebox(30,10){{\bf aRML**} }}
\put(-80,80){\makebox(30,10){{\bf *aRML**} }}
\put(-155,140){\makebox(30,10){ (**)}}   \put(-115,140){\makebox(30,-80){ (**)}}
\put(-72,351){\makebox(2,2){$\bigcirc$}}  \put(-72,360){\makebox(10,10){{\bf RML} }}
\put(-8,332){\makebox(2,2){$\bigcirc$}} \put(-8,332){\makebox(2,2){$\bullet$}}
\put(-65,300){\makebox(10,10){{\bf *RML} }}

\put(-145,260){\makebox(30,10){ (An)}}
\put(-50,190){\makebox(30,10){ (An)}}
                  \put(-60,240){\makebox(30,10){ (An)}}
\put(-120,175){\makebox(10,10){(An)}}

\put(-60,320){\makebox(10,10){(*)}}  \put(10,312){\makebox(10,10){(*)}}

\put(-135,113){\makebox(2,2){\circle*{11}}}
 \put(-160,120){\makebox(10,10){(*)}}  \put(-90,112){\makebox(10,10){(*)}}

\put(-170,90){\makebox(30,10){{\bf *aRML} }}

\put(-176,147){\makebox(2,2){$\bigcirc$}}  \put(-176,147){\makebox(2,2){$\circ$}} 
\put(-195,158){\makebox(30,10){{\bf aRML}}}
                  \put(143,300){\makebox(10,10){{\it Hierarchy 2.2 H} }}
                            \put(-22,360){\makebox(10,10){{\it Hierarchy 2.2} }}

\end{picture}
 \end{center}
\caption{ Hierarchies 2.2 and 2.2 H, determined by RML and by pi-RML algebras, considering  (*) and (**) properties   }
\label{fig:figrml***H}
\end{figure}
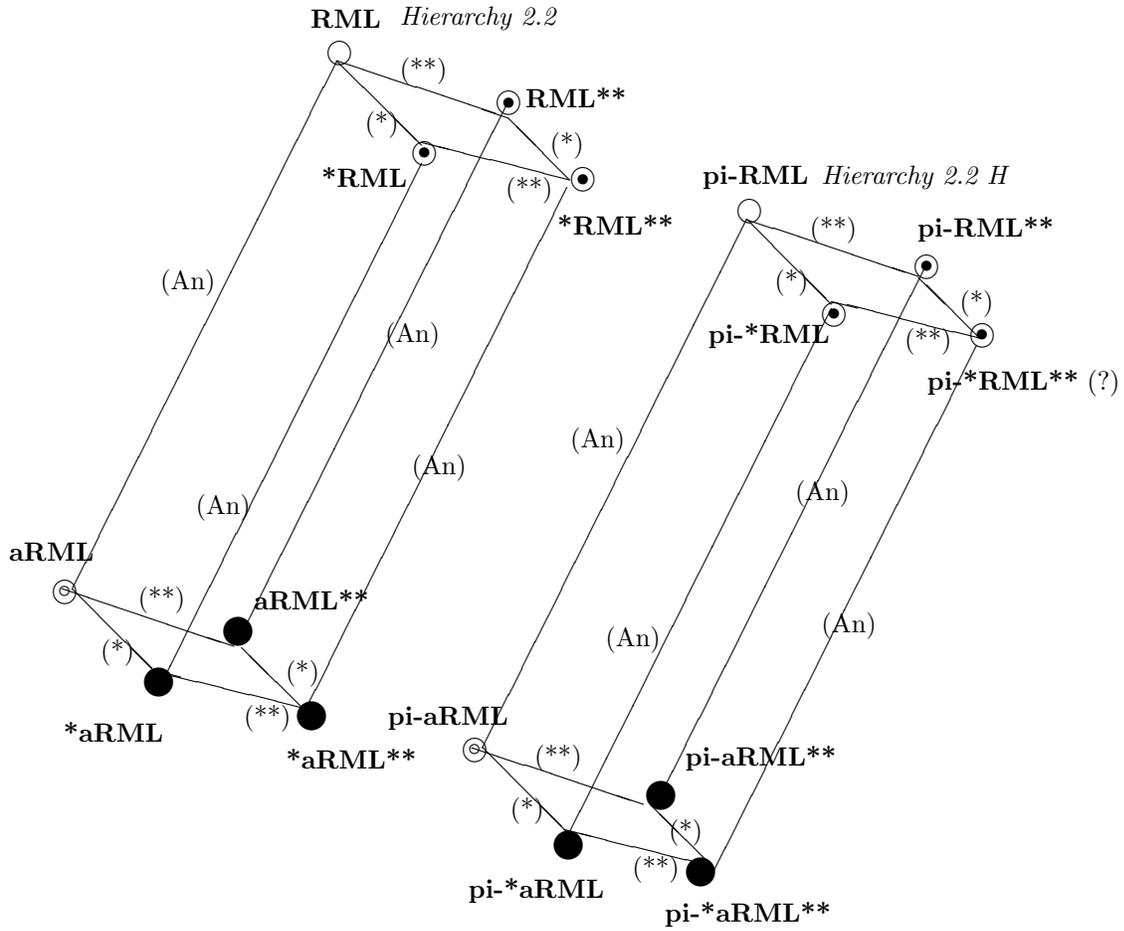

\begin{figure}[htbp]
\begin{center}
\begin{picture}(420,340)(-180,25)

       \put(-17,89){\line(1,-1){60}}    
       \put(-17,87){\line(1,-1){60}}    
       \put(-17,85){\line(1,-1){60}}

  \put(45,25){\line(1,2){100}}                               
  \put(-17,89){\line(1,2){100}}

      \put(83,289){\line(1,-1){60}}   
      \put(83,287){\line(1,-1){60}}   
       \put(83,285){\line(1,-1){60}}   
 \put(146,228){\line(1,-2){72}} 
                                \put(45,22){\line(3,1){171}}   
                                 \put(217,82){\makebox(2,2){$\bigcirc$}}
                                  \put(217,82){\makebox(2,2){$\bullet$}} 
                                  \put(190,65){\makebox(10,10){{\bf pimpl-pre-BCK} }}
                                  \put(150,45){\makebox(2,2){(An)}} 
                                   \put(200,135){\makebox(2,2){(pimpl)}}

\put(83,289){\line(2,-1){85}} 
\put(83,287){\line(2,-1){85}} 

       \put(146,228){\line(2,1){25}}  
\put(171,244){\makebox(2,2){$\bigcirc$}} \put(171,244){\makebox(2,2){$\bullet$}}\put(185,230){\makebox(10,10){{\bf pi-BE**} }}

\put(140,265){\makebox(10,10){ (**)$\Leftrightarrow$ (Tr)}}

\put(-21,89){\line(2,-1){85}} 
\put(-21,87){\line(2,-1){85}}   
  \put(45,25){\line(2,1){25}}    
\put(70,41){\makebox(2,2){\circle*{11}}}
                  
                    \put(70,41){\line(1,2){100}} 
\put(78,51){\makebox(30,10){{\bf pi-aBE**} }}
\put(32,70){\makebox(10,10){ (**)$\Leftrightarrow$ (Tr)}}

\put(83,291){\makebox(2,2){$\bigcirc$}}  \put(83,300){\makebox(10,10){{\bf pi-BE} }}
\put(141,226){\makebox(2,2){$\bigcirc$}} \put(141,226){\makebox(2,2){$\bullet$}}\put(132,210){\makebox(10,10){{\bf pi-pre-BCK} }}

\put(10,200){\makebox(30,10){ (An)}}
\put(105,130){\makebox(30,10){ (An)}}
                  \put(95,180){\makebox(30,10){ (An)}}

\put(95,250){\makebox(10,10){(B)}}
\put(95,240){\makebox(10,10){$\Leftrightarrow$ (BB)}}
\put(95,230){\makebox(10,10){$\Leftrightarrow$ (*)}}

\put(45,25){\makebox(2,2){\circle*{11}}}
 \put(-5,50){\makebox(10,10){(B)}}
\put(-5,40){\makebox(10,10){$\Leftrightarrow$ (BB)}}
\put(-5,30){\makebox(10,10){$\Leftrightarrow$ (*)}}

\put(32,8){\makebox(30,10){{\bf pi-BCK = pimpl-BCK = Hilbert} }}

\put(-21,87){\makebox(2,2){$\bigcirc$}}  \put(-21,87){\makebox(2,2){$\circ$}} 
\put(-50,65){\makebox(30,10){{\bf pi-aBE}}}
 \put(-172,149){\line(1,-1){60}}    
       \put(-172,147){\line(1,-1){60}}    
       \put(-172,145){\line(1,-1){60}}

  \put(-110,85){\line(1,2){100}}                               
\put(-172,149){\line(1,2){100}}

      \put(-72,349){\line(1,-1){60}}   
      \put(-72,347){\line(1,-1){60}}   
      \put(-72,345){\line(1,-1){60}}   
\put(-72,349){\line(2,-1){85}}
\put(-72,347){\line(2,-1){85}}   
   
  \put(-9,288){\line(2,1){25}}  
\put(16,304){\makebox(2,2){$\bigcirc$}} \put(16,304){\makebox(2,2){$\bullet$}}\put(35,290){\makebox(10,10){{\bf BE**} }}
 \put(-15,325){\makebox(10,10){ (**)$\Leftrightarrow$ (Tr)}}

\put(-176,149){\line(2,-1){85}}  
\put(-176,147){\line(2,-1){85}}  
 
  \put(-110,85){\line(2,1){25}}    
\put(-85,100){\makebox(2,2){\circle*{11}}}
 
 \put(-85,101){\line(1,2){100}} 
\put(-75,91){\makebox(30,10){{\bf aBE**} }}
\put(-130,133){\makebox(10,10){ (**)$\Leftrightarrow$ (Tr)}}

\put(-72,351){\makebox(2,2){$\bigcirc$}}  \put(-72,360){\makebox(10,10){{\bf BE} }}
\put(-9,288){\makebox(2,2){$\bigcirc$}} \put(-9,288){\makebox(2,2){$\bullet$}}\put(-23,270){\makebox(10,10){{\bf pre-BCK} }}

\put(-145,260){\makebox(30,10){ (An)}}
\put(-40,190){\makebox(30,10){ (An)}}
                 \put(-60,240){\makebox(30,10){ (An)}}

\put(-60,310){\makebox(10,10){(B)}}
\put(-60,300){\makebox(10,10){$\Leftrightarrow$ (BB)}}
\put(-60,290){\makebox(10,10){$\Leftrightarrow$ (*)}}

\put(-103,85){\makebox(2,2){\circle*{11}}}
 \put(-160,110){\makebox(10,10){(B)}}
\put(-160,100){\makebox(10,10){$\Leftrightarrow$ (BB)}}
\put(-160,90){\makebox(10,10){$\Leftrightarrow$ (*)}}

\put(-123,68){\makebox(30,10){{\bf BCK} }}

\put(-176,147){\makebox(2,2){$\bigcirc$}}  \put(-176,147){\makebox(2,2){$\circ$}} 
\put(-215,143){\makebox(30,10){{\bf aBE}}}
 \put(143,300){\makebox(10,10){{\it Hierarchy 2.3 H} }}
                            \put(-22,360){\makebox(10,10){{\it Hierarchy 2.3} }}

\end{picture}
 \end{center}
\caption{ Hierarchies 2.3 and 2.3 H,  determined by BE  and  pi-BE algebras   }
\label{fig:fig23ssTH}
\end{figure}
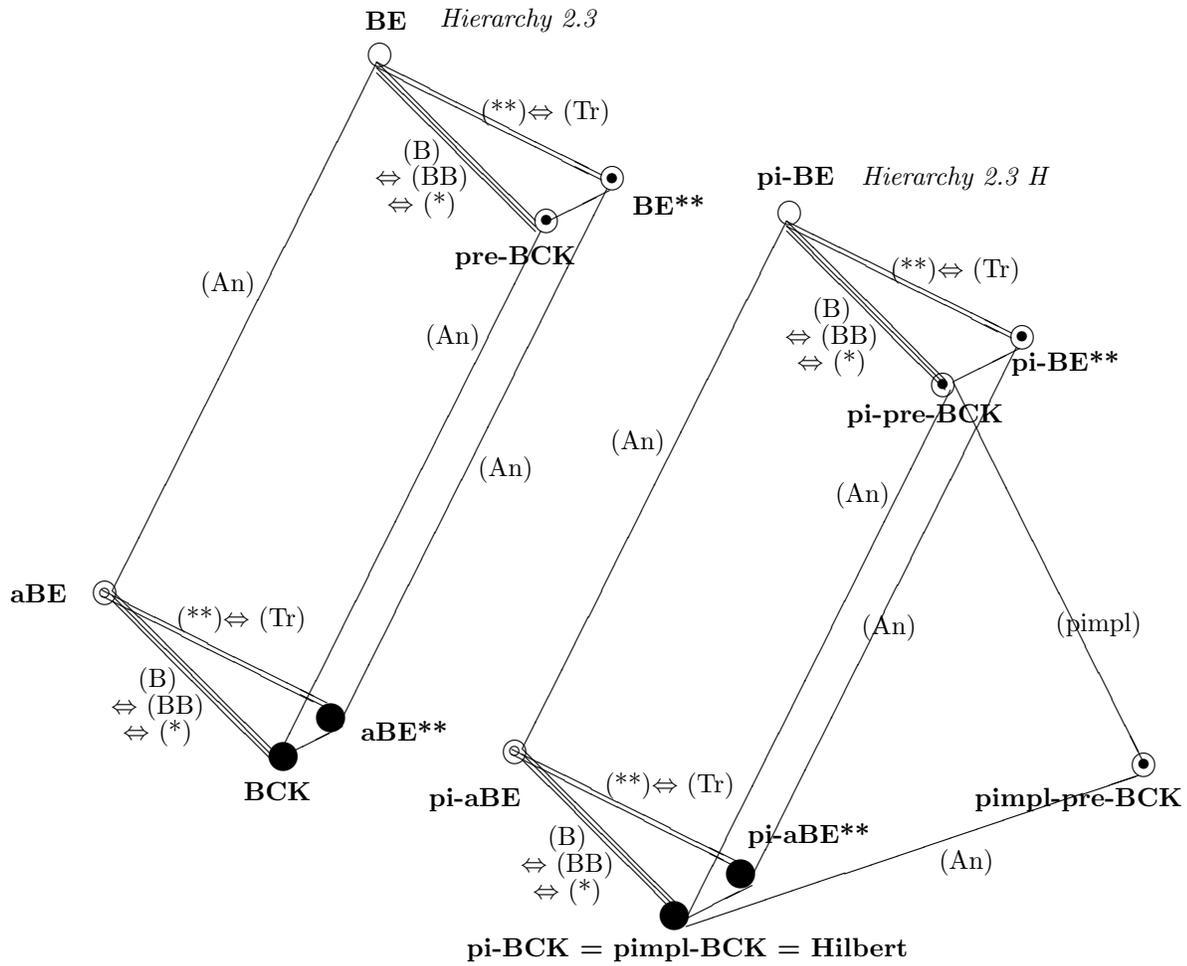

                           \section{Proper algebras}

In this section, we define the proper algebras discussed in the paper; the  definitions  are needed in the next sections,
where we present examples
of proper algebras. The definitions are denoted by: PO1-PO7 (following the corresponding  notations O1 - O7) and  P1 - P31
 (following
the corresponding notations 1 - 31), for the generalizations of BCI and BCK algebras  
and  are denoted by Pnm-pi or Pnm-pimpl (following
the corresponding notation Pnm), for the corresponding generalizations of  Hilbert algebras verifying (pi) or (pimpl).

 Hence, we have the following obvious definitions.

\begin{df}\em ${}$

PO1.  A  {\it proper  BCI algebra}  is a BCI algebra (i.e. verifying  
(BB), (D), (Re), (N), (An) or, equivalently, (BB), (D), (Re), (An), or, equivalently, (BB), (M), (An)) not verifying (L).

PO2. A {\it proper  BCK algebra} is a BCK algebra (i.e. verifying  (BB), (D),  (Re), (L), (An), 
or, equivalently,  (BB), (M), (L), (An)) not verifying (pi) (hence (pimpl)).

PO3.  A {\it proper BCH algebra} is a BCH algebra (i.e. verifying (Re), (Ex), (An)) 
not verifying (B), (BB), (*), (**), (Tr), (L).  

PO4. A {\it proper BCC algebra}  is a BCC algebra (i.e. verifying (Re), (M), (L), (B), (An)) 
not verifying  (Ex), (BB), (pi).\\
PO4-pi. A {\it proper pi-BCC algebra}  is a pi-BCC algebra (i.e. verifying (Re), (M), (L), (B), (An), (pi)) 
not verifying  (Ex), (pimpl).

PO5.  A {\it proper  BZ algebra} is a BZ algebra (i.e. verifying (Re), (M), (B), (An))
not verifying (L), (Ex), (BB).

PO6.  A {\it proper BE algebra} is a BE algebra (i.e. verifying  (Re), (M), (L), (Ex) (hence (D))) not
verifying  (An), (Tr), (*) (hence (B)), (**) (hence (BB)),  (pi).\\
PO6-pi.  A {\it proper pi-BE algebra} is a pi-BE algebra (i.e. verifying  (Re), (M), (L), (Ex) (hence (D)), (pi)) not
verifying  (An), (B), (BB), (*), (**), (Tr),  (pimpl).

PO7. A {\it proper pre-BCK algebra} is a pre-BCK algebra (i.e. verifying  (Re), (M), (L), (Ex) (hence (D)), (*)
(hence (B), (BB), (**), (Tr))) not
verifying (An), (pi). \\
PO7-pi. A {\it proper pi-pre-BCK algebra} is a pi-pre-BCK algebra (i.e. verifying  (Re), (M), (L), (Ex), (*), (pi)) not
verifying (An), (pimpl).  \\
PO7-pimpl. A {\it proper pimpl-pre-BCK algebra} is a pimpl-pre-BCK algebra 
(i.e. verifying  (Re), (M), (L), (Ex), (*), (pimpl)) not
verifying (An).  
\end{df}

\begin{df}\em ${}$

P1. A {\it proper pre-BCC algebra} is a pre-BCC algebra (i.e. verifying (Re), (M), (L), (B) (hence (*), (Tr))) not
verifying (An), (Ex), (BB),  (pi).\\
P1-pi. A {\it proper pi-pre-BCC algebra} is a pi-pre-BCC algebra
 (i.e. verifying (Re), (M), (L), (B) (hence (*), (Tr)), (pi)) not
verifying (An), (Ex), (BB),   (pimpl).

P2. A {\it proper aBE algebra} is an aBE algebra (i.e. verifying (Re), (M), (L), (Ex) (hence (D)), (An))
not verifying (B), (BB), (*), (**), (Tr), (pi). \\
P2-pi. A {\it proper pi-aBE algebra} is a pi-aBE algebra (i.e. verifying (Re), (M), (L), (Ex) (hence (D)), (An), (pi))
not verifying (B), (BB), (*), (**), (Tr), (pimpl).

P3. A {\it proper RM algebra} is an  RM algebra (i.e. verifying (Re), (M)) 
not  verifying  (Ex), (An), (L), (B), (BB), (*), (**), (Tr).

P4. A {\it proper pre-BZ algebra} is a pre-BZ algebra (i.e. verifying  (Re), (M), (B) (hence (*), (Tr)))
not  verifying  (Ex), (An), (L), (BB).

P5. A {\it proper aRM algebra} is an aRM  algebra (i.e. verifying (Re), (M), (An))
not  verifying  (L), (Ex), (B), (*), (BB), (**), (Tr).

P6. A  {\it proper RME algebra} is a RME algebra  (i.e. verifying (Re), (M), (Ex) (hence (D)))
not  verifying  (An), (L), (B), (BB), (*), (**), (Tr).

P7. A {\it proper pre-BCI algebra} is a pre-BCI algebra 
(i.e. verifying  (Re), (M), (Ex), (B) (hence (BB), (*), (**), (Tr)))
not  verifying  (An), (L).

P8. A {\it proper RML algebra} is a RML algebra  (i.e. verifying  (Re), (M), (L))
not  verifying  (Ex), (An), (B), (*), (BB), (**), (Tr), (pi).\\
P8-pi. A {\it proper pi-RML algebra} is a pi-RML algebra  (i.e. verifying  (Re), (M), (L), (pi))
not  verifying  (Ex), (An), (B), (*), (BB), (**), (Tr), (pimpl).

P9. A {\it proper  aRML algebra} is an aRML algebra (i.e.  verifying  (Re), (M), (L), (An))
   not  verifying  (Ex), (B), (*), (BB), (**), (Tr), (pi).\\
P9-pi. A {\it proper  pi-aRML algebra} is a pi-aRML algebra (i.e.  verifying  (Re), (M), (L), (An), (pi))
   not  verifying  (Ex), (B), (*), (BB), (**), (Tr), (pimpl).
\end{df}
\begin{df}${}$

P10. A   {\it proper tRM algebra}   is a tRM algebra (i.e.  verifying (Re), (M), (Tr))
  not verifying (Ex), (An), (L), (*) (hence (B)), (**) (hence (BB)).

P11. A   {\it proper *RM algebra}   is a *RM algebra (i.e. verifying (Re), (M), (*) (hence (Tr)))
  not verifying  (Ex), (An), (L), (B), (**) (hence (BB)).

P12. A   {\it proper RM** algebra}   is a RM**  algebra (i.e. verifying  (Re), (M), (**) (hence (Tr)))
  not verifying  (Ex), (An), (L), (BB), (*) (hence (B)).

P13. A   {\it proper *RM** algebra}   is a *RM** algebra (i.e. verifying (Re), (M), (*), (**) (hence (Tr)))
   not verifying (Ex), (An), (L), (B), (BB).

P14. A   {\it proper pre-BBBZ algebra}   is a pre-BBBZ algebra  (i.e. verifying  (Re), (M), (B),
 (BB) (hence (*), (**), (Tr), (D)))
  not verifying (Ex), (An), (L).

P15. A   {\it proper oRM algebra}   is an  oRM algebra (i.e. verifying (Re), (M), (An), (Tr))
  not verifying (*), (**), (L).

P16. A   {\it proper *aRM algebra}   is an  *aRM algebra (i.e. verifying  (Re), (M), (An), (*) (hence (Tr)))
  not verifying  (Ex), (L), (B),  (**) (hence  not (BB)).

P17. A   {\it proper aRM** algebra}   is an  aRM** algebra  (i.e. verifying (Re), (M), (An), (**) (hence (Tr)))
 not verifying (Ex), (L),  (BB), (*)  (hence not (B)).

P18. A   {\it proper *aRM** algebra}   is a  *aRM** algebra  (i.e. verifying (Re), (M), (An), (*), (**) (hence (Tr)))
 not verifying (Ex), (L), (B), (BB).

P19. A   {\it proper tRML algebra}   is a tRML algebra  (i.e. verifying (Re), (M),  (L), (Tr))
 not verifying (An),  (Ex), (*) (hence (B)), (**) (hence (BB)),  (pi).\\
P19-pi. A   {\it proper pi-tRML algebra}   is a pi-tRML algebra  (i.e. verifying (Re), (M),  (L), (Tr), (pi))
 not verifying  (An), (Ex), (*), (**),  (pimpl).

P20. A   {\it proper *RML algebra}   is a *RML algebra  (i.e. verifying (Re), (M), (L),  (*) (hence (Tr)))
 not verifying  (An), (Ex), (B), (**) (hence (BB)),  (pi).\\
P20-pi. A   {\it proper pi-*RML algebra}   is a pi-*RML algebra  (i.e. verifying (Re), (M), (L),  (*) (hence (Tr)), (pi))
 not verifying (An), (Ex), (B), (**), (BB),  (pimpl).

P21. A   {\it proper RML** algebra}   is a RML** algebra (i.e. verifying  (Re), (M), (L), (**) (hence (Tr)))
 not verifying  (An), (Ex), (BB), (*) (hence (B)),  (pi).\\
P21-pi. A   {\it proper pi-RML** algebra}   is a pi-RML** algebra (i.e. verifying  (Re), (M), (L), (**) (hence (Tr)), (pi))
 not verifying (An), (Ex), (BB), (*), (B), (pimpl).

P22. A   {\it proper *RML** algebra}   is a *RML** algebra  (i.e. verifying (Re), (M), (L), (*), (**) (hence (Tr)))
 not verifying (An), (Ex), (B), (BB),  (pi).\\
P22-pi. A   {\it proper pi-*RML** algebra}   is a pi-*RML** algebra  (i.e. verifying (Re), (M), (L), (*), (**) (hence (Tr)), (pi))
 not verifying  (An), (Ex), (B), (BB),  (pimpl).

P23. A  {\it proper pre-BBBCC algebra} is a pre-BBBCC algebra  
(i.e. verifying (Re), (M), (L), (B), (BB) (hence (*), (**), (Tr), (D))) not verifying (An), (Ex),(pi).\\
P23-pi. A  {\it proper pi-pre-BBBCC algebra} is a pi-pre-BBBCC algebra  
(i.e. verifying (Re), (M), (L), (B), (BB) (hence (*), (**), (Tr), (D)), (pi)) not verifying (An), (Ex), (pimpl).\\
P23-pimpl. A  {\it proper pimpl-pre-BBBCC algebra} is a pimpl-pre-BBBCC algebra  
(i.e. verifying (Re), (M), (L), (B), (BB) (hence (*), (**), (Tr), (D)), (pimpl)) not verifying (An), (Ex).

P24. A   {\it proper oRML algebra}   is an oRML algebra  (i.e. verifying  (Re), (M), (L), (An), (Tr))
 not verifying  (Ex), (*) (hence (B)), (**) (hence (BB)), (pi).\\
P24-pi. A   {\it proper pi-oRML algebra}   is a pi-oRML algebra  (i.e. verifying  (Re), (M), (L), (An), (Tr), (pi))
 not verifying (Ex), (*) (hence (B)), (**) (hence (BB)), (pimpl).

P25. A   {\it proper  *aRML algebra}   is a  *aRML algebra  (i.e. verifying  (Re), (M), (L), (An), (*) (hence (Tr)))
not verifying (B), (**) (hence (BB)), (pi).\\
P25-pi. A   {\it proper  pi-*aRML algebra}   is a  pi-*aRML algebra 
 (i.e. verifying  (Re), (M), (L), (An), (*) (hence (Tr)), (pi))
not verifying (B),  (**) (hence (BB)), (pimpl).

P26. A   {\it proper aRML** algebra}   is an aRML** algebra (i.e.  verifying (Re), (M), (L), (An), (**) (hence (Tr)))
 not verifying (Ex), (BB), (*) (hence (B)), (pi).\\
P26-pi. A   {\it proper pi-aRML** algebra}   is a pi-aRML** algebra 
(i.e.  verifying (Re), (M), (L), (An), (**) (hence (Tr)), (pi))
 not verifying (Ex), (BB), (*) (hence  (B)), (pimpl).

P27. A   {\it proper *aRML** algebra}   is a *aRML** algebra  
(i.e. verifying (Re), (M), (L), (An), (*), (**) (hence (Tr)))
 not verifying (Ex), (B), (BB), (pi).\\
P27-pi. A   {\it proper pi-*aRML** algebra}   is a pi-*aRML** algebra 
 (i.e. verifying (Re), (M), (L), (An), (*), (**) (hence (Tr)), (pi))
 not verifying (Ex), (B), (BB), (pimpl).

P28. A   {\it proper RME** algebra}   is a RME** algebra  
(i.e. verifying (Re), (M), (Ex) (hence (D)), (**) (hence (Tr)), 
 not verifying (An),  (BB), (*) (hence (B)).

P29. A   {\it proper BCH** algebra}   is a BCH** algebra  
(i.e. verifying (Re), (M), (Ex) (hence (D)), (**) (hence (Tr)), (An)), 
 not verifying   (BB), (*) (hence (B)).

P30. A   {\it proper BE** algebra}   is a BE** algebra  
(i.e. verifying (Re), (M), (L), (Ex) (hence (D)), (**) (hence (Tr)), 
 not verifying (An),  (BB), (*) (hence (B)), (pi).\\
P30-pi. A   {\it proper pi-BE** algebra}   is a pi-BE** algebra 
 (i.e. verifying (Re), (M), (L), (Ex) (hence (D)), (**) (hence (Tr)), (pi))
 not verifying (An),  (BB), (*) (hence (B)), (pimpl).

P31. A   {\it proper aBE** algebra}   is an aBE** algebra  
(i.e. verifying (Re), (M), (L), (Ex) (hence (D)), (An), (**) (hence (Tr)), 
 not verifying   (BB), (*) (hence (B)), (pi).\\
P31-pi. A   {\it proper pi-aBE** algebra}   is a pi-aBE** algebra 
 (i.e. verifying (Re), (M), (L), (Ex) (hence (D)), (**) (hence (Tr)), (pi))
 not verifying  (BB), (*) (hence (B)), (pimpl).
\end{df}

In the remaining sections we shall give proper examples of all the algebras, old and new,  
discussed in the paper.
 The examples  will be presented mainly 
in the same order as the corresponding algebras  were introduced: \\
- first, the sixteen old and new algebras from Hierarchy 3 (the eight  RM algebras, generalizations of BCI algebras)
 and  Hierarchy 2 
(the eight RML algebras, 
proper generalizations of BCK algebras), \\
- then, the twenty two new algebras from Hierarchies 3.1 - 3.3  (other eleven RM algebras) and  Hierarchies 2.1 - 2.3
 (other eleven RML algebras),\\
- finally, the (proper) generalizations of Hilbert algebras (pi-RML algebras and pimpl-RML algebras).

                                    \section{Examples of the first sixteen  old and  new   RM and RML algebras
 PO1 - PO7, P1 - P9}  

For a coeherent presentation of the examples, we shall provide  also  examples of the seven old algebras.

The examples will be presented in the following order: first the RM algebras from Hierarchy 3, 
then the RML algebras (Hierarchy 2).

Note that even if the algebras were obtained from the bottom to the top, 
the examples will be presented from the  top to the bottom.

                  \subsection{Examples of the RM algebras from Hierarchy 3}

Recall that we have proved that all the algebras from Hierarchy 3 (i.e. not satisfying condition (L)) cannot have
the properties (pimpl) and (pi).  

In the first subsubsection, we present proper examples of RM algebras not satisfying  property (Ex) and in the second
subsubsection, we present examples of RM algebras satisfying property (Ex).
                  
                          \subsubsection{ Examples of RM algebras  without condition (Ex):\\ proper RM, pre-BZ, aRM, BZ algberas }
${}$

{\bf $\bullet$  Proper RM algebras} (P3)

\begin{ex}\em
Let $(A=\{a,b,1\}, \ra,1)$ with the following matrix (table)  of implication:\\

\noindent  $\qquad \qquad \qquad \qquad$
\begin{tabular}{*{1}{c|}ccc}
$\ra$ & a & b& 1  \\
\hline
a& 1& 1 & a   \\
b& 1 & 1&  1\\
1& a& b& 1 
\end{tabular}
 
\noindent\\
Properties  (Re), (M), (D) are satisfied. (An) is not satisfied for $(x,y)=(a,b)$. (L) is not satisfied for $x=a$.
(Ex) is not satisfied for $(x,y,z)=(a,b,a)$,  
   (BB) is not satisfied for $(x,y,z)=(a,b,1)$,  
(**) is not satisfied for $(x,y,z)=( a,b,1)$, (B) is not satisfied for $(x,y,z)=(a,b,1)$,  
 (*) is not satisfied for $(x,y,z)=(a,b,1)$, (Tr) is not satisfied for $(x,y,z)=(a,b,1)$.\\
 Hence, ${\cal A}$ is a  {\bf proper RM algebra} with the minimum number of elements, three, with (D).
\end{ex}

\begin{ex}\em
Let $(A=\{a,b,1\}, \ra,1)$ with the following matrix   of implication:\\

\noindent $\qquad \qquad \qquad \qquad$
\begin{tabular}{*{1}{c|}ccc}
$\ra$ & a & b& 1  \\
\hline
a& 1& 1 & b   \\
b& 1 & 1&  1\\
1& a& b& 1 
\end{tabular} 

\noindent\\
Properties (Re), (M) are satisfied.
(Ex) is not satisfied for a,b,a; (An) is not satisfied by a,b;
 (BB) is not satisfied for a,b,1; (**) is not satisfied for a,b,1; (B) is not satisfied for a,b,1;  
 (*) is not satisfied for a,b,1; (Tr) is not satisfied for a,b,1; (D)
 is not satisfied for 1,a.\\
Hence, ${\cal A}$ is a {\bf proper RM algebra} with three elements too, without (D).
\end{ex}

\begin{ex}\em
Consider the set $A=\{0, a,b, 1\}$ with the following table of $\ra$: \\ 

\noindent  $\qquad \qquad \qquad \qquad$
\begin{tabular}{*{1}{c|}cccc} 
$\ra$ & 0 & a & b & 1  \\
\hline
0 & 1 & 1 & 1 & 1  \\
a & 1 & 1 & 1 & a \\
b & 1 & 1 & 1 & a  \\
1 & 0 & a & b & 1 
\end{tabular} 

\noindent\\
 Then the algebra $ {\cal A}=(A, \ra, 1)$ verifies properties (Re), (M) and (D).
 It does not verify properties: (Ex) for $x=a, \; y=0, \; z=b$, (BB) for 
$x=0, \; y=a, \; z=1$, (B) for $x=0, \; y=a, \; z=1$,   (**) for $x=a, \; y=0, \; z=1$, (*) for 
$x=a, \; y=0, \; z=1$, (Tr) for $x=a, \; y=0, \; z=1$, (An) for $x=a, \; y=b$.\\
  Hence, {\bf ${\cal A}$ is a proper RM algebra}, with four elements, with (D).
\end{ex} 

\begin{ex}\em
Consider the set $A=\{0, a,b, 1\}$ with the following table of $\ra$: \\ 

\noindent  $\qquad \qquad \qquad \qquad$
\begin{tabular}{*{1}{c|}cccc}
$\ra$ & 0 & a & b & 1  \\
\hline
0 & 1 & 1 & 1 & b  \\
a & 1 & 1 & a & b \\
b & 0 & a & 1 & 1  \\
1 & 0 & a & b & 1 
\end{tabular} 

\noindent\\
 Then the algebra $ {\cal A}=(A, \ra, 1)$ verifies properties (Re), (M).
 It does not verify properties (Ex), (BB), (B), (**), (*), (Tr), (An), (D). \\
  Hence,  ${\cal A}$ is a  {\bf proper RM algebra},  with four elements, without (D).
\end{ex}

{\bf $\bullet$  Proper pre-BZ algebras} (P4)

\begin{ex}\em
Consider the set $A=\{a, b,c,d, 1\}$ with the following table of $\ra$: \\ 

\noindent  $\qquad \qquad \qquad \qquad$ 
\begin{tabular}{*{1}{c|}ccccc}
$\ra$ & a & b & c & d & 1 \\
\hline
a  & 1& a & a & b & a \\
b & a & 1 & 1 & d & 1\\ 
c & a & 1 & 1 & d & 1\\ 
d & 1 & a & a & 1 & a\\ 
1 & a & b & c & d & 1 
\end{tabular} 

\noindent\\
Properties  (Re), (M), (B) (hence (*), (Tr), (**))  are satisfied. Properties (An) and (L) do not hold obviously.
Property  (Ex) does not hold   for $a,b,d$; 
  (BB) does not hold  for $d,b,a$; (D) is not satisfied for $d,a$.\\
Hence, ${\cal A}$ is a {\bf proper pre-BZ algebra},  without (D).
\end{ex} 
{\bf $\bullet$  Proper aRM algebras} (P5)

\begin{ex}\em
 Consider the set $A=\{a, b,c, 1\}$ with the following table of $\ra$: \\ 

\noindent  $\qquad \qquad \qquad \qquad$
\begin{tabular}{*{1}{c|}cccc}
$\ra$ & a & b & c & 1  \\
\hline
a&  1& a& a& a\\ 
b&  a &1& a& 1\\ 
c&  a& 1& 1& a\\ 
1&  a& b& c& 1 
\end{tabular} 

\noindent\\
${\cal A}$ satisfies (Re), (M), (An).
  (Ex) is not satisfied for $a,b,c$;  (L) is not satisfied for $a$; (BB) is not satisfied for $a,b,c$;  (**) does  not hold for $c,b,1$;
  (B) does not hold for $a,b,c$;
 (*) is not satisfied for $b,c,b$;  (Tr) is  not satisfied for $c,b,1$;   (D) is not satisfied for $a,c$.\\
Hence, ${\cal A}$ is a {\bf proper aRM algebra}, without (D).
\end{ex}

{\bf $\bullet$  Proper BZ algebras} (PO5)

\begin{ex}\em
Consider the set $A=\{a,b,c,d,1\}$ with the following table of $\ra$: \\ 

\noindent   $\qquad \qquad \qquad \qquad$
\begin{tabular}{*{1}{c|}ccccc}
$\ra$ & a & b & c & d & 1 \\
\hline
a & 1 &a& a &a& a \\
b & a &1& b& c &1 \\
c & a &1& 1& c &1 \\
d & a& 1& 1& 1& 1 \\
1 & a& b& c &d &1 
\end{tabular} 

\noindent\\
Properties  (Re), (M), (An), (B) (hence (*), (Tr), (**)) are satisfied.
Property (Ex) is not satisfied for $b,c,d$,   (BB) for $d,1,b$,    (D) for $d,b$.\\
Hence, ${\cal A}$ is a {\bf proper BZ algebra},   without (D).
\end{ex}

       \subsubsection{ Examples of RM algebras  with condition (Ex): \\RME, pre-BCI, BCH, BCI }

These algebras verify all property (D), by Proposition \ref{propp} (4.).\\

{\bf $\bullet$  Proper RME algebras} (P6)
 
\begin{ex}\em
 Consider the set $A=\{a,b,c,d,1\}$ with the following table of $\ra$: \\ 

\noindent  $\qquad \qquad \qquad \qquad$ 
\begin{tabular}{*{1}{c|}ccccc}
$\ra$ & a & b & c & d & 1 \\
\hline
a & 1 & a & a & a & a\\ 
b & a & 1 & c & d & 1 \\
c & a & b & 1 & 1 & 1 \\
d & a & 1 & 1 & 1 & 1 \\
1 & a & b & c&  d & 1 
\end{tabular} 

\noindent\\
${\cal A}$ satisfies (Re), (M), (Ex). It does not satisfies:  (L) for $a$,  (An) for $c,d$,   (BB) for $b,c,d$,
  (**) for $b,c,d$,   (B) and  (*) for $c,d,b$,  (Tr) for $c,d,b$.\\
Hence, ${\cal A}$ is a {\bf proper RME algebra}.
\end{ex}

{\bf $\bullet$  Proper pre-BCI algebras} (P7)

\begin{ex}\em
 Consider the set $A=\{a,b,1\}$ with the following table of $\ra$: \\ 

\noindent  $\qquad \qquad \qquad \qquad$
\begin{tabular}{*{1}{c|}ccc}
$\ra$ & a & b& 1  \\
\hline
a & 1 & 1 & a \\
b & 1 & 1 & a \\
1 & a & b & 1 
\end{tabular} 

\noindent\\ 
Properties (Re), (M), (Ex), (B) (hence (BB), (*), (**), (Tr)) are  satisfied. 
(L) and (An) are not satisfied obviously.\\
Hence, ${\cal A}$ is a {\bf proper pre-BCI algebra}.
\end{ex}

{\bf $\bullet$ Proper BCH algebras} (PO3)

\begin{ex}\em
 Consider the set $A=\{a,b,c,d,1\}$ with the following table of $\ra$: \\ 

\noindent   $\qquad \qquad \qquad \qquad$
\begin{tabular}{*{1}{c|}ccccc}
$\ra$ & a & b & c & d & 1 \\
\hline
a & 1 & a & a & a & a \\
b & a & 1 & b & 1 & 1 \\
c & a & 1 & 1 & d & 1 \\
d & a & b & c & 1 & 1 \\
1 & a & b & c & d & 1 
\end{tabular} 

\noindent\\
Properties (Re), (M), (Ex), (An) are satisfied. ${\cal A}$ does not satisfy  (L) for $a$,
  (BB) and (**)  for $d,c,b$,
(B), (*), (Tr) for $c,b,d$. \\
Hence, ${\cal A}$ is a {\bf proper BCH algebra}.
\end{ex}

{\bf $\bullet$  Proper BCI algebras} (PO1)

\begin{ex}\em
 
Consider the set $A=\{a,b,1\}$ with the following table of $\ra$: \\ 

\noindent  $\qquad \qquad \qquad \qquad$
\begin{tabular}{*{1}{c|}ccc}
$\ra$ & a & b& 1  \\
\hline
a & 1 & a & a \\
b & a & 1 & 1 \\
1 & a & b & 1 
\end{tabular} 

\noindent\\
Properties (Re), (An), (BB), (D), (N)  are satisfied. (L)  is not satisfied for $a$.\\
Hence, ${\cal A}$ is a {\bf proper BCI algebra}.
\end{ex}


                          \subsection{Examples of the RML algebras from Hierarchy 2}

Recall that we have proved that only these algebras can have the properties (pimpl) and (pi).

In the first subsubsection we present proper examples of RML algebras not satisfying  property (Ex) and in the second
subsubsection we present examples of RML algebras satisfying property (Ex).

              \subsubsection{ Examples of RML algebras  without condition (Ex):\\ RML, pre-BCC, aRML, BCC }
${}$

{\bf $\bullet$  Proper RML algebras} (P8)

\begin{ex}\em

Consider the set $A=\{a,b,c,1\}$ with the following table of $\ra$: \\ 

\noindent  $\qquad \qquad \qquad \qquad$
\begin{tabular}{*{1}{c|}cccc}
$\ra$ & a & b & c & 1  \\
\hline
a & 1 & a & a & 1 \\
b & a & 1 & 1 & 1 \\
c & 1 & 1 & 1 & 1 \\
1 & a & b & c & 1 
\end{tabular} 

\noindent\\
${\cal A}$ satisfies (Re), (M), (L), (D). It does not satisfy   (Ex) for $a,b,b$,
 (pi) for $b,a$,
 (BB), (**) for $a,b,c$;  (B), (*), (Tr) for $b,c,a$. \\
Hence, ${\cal A}$ is a {\bf proper RML algebra}, with (D).
\end{ex}

\begin{ex}\em
Consider the set $A=\{0, a,b,c, 1\}$ with the following table of $\ra$: \\ 

\noindent   $\qquad \qquad \qquad \qquad$
\begin{tabular}{*{1}{c|}ccccc}
$\ra$ & 0 & a & b & c & 1 \\
\hline
0 & 1 & 1 & 1 & 1 & 1  \\
a & c & 1 & 1 & 1 & 1 \\
b & b & {\bf a} & 1 & 1 & 1  \\
c &  a & {\bf 1} & {\bf a} & 1 &  1  \\
1 & 0 & a & b & c & 1
\end{tabular} 

\noindent\\
Then the bounded algebra $ {\cal A}=(A, \ra,0, 1)$ verifies properties (Re), (M), (L),(D), and (DN).\\
 It does not verify properties (Ex), (An), (BB), (B), (*), (**), (Tr),  (pi).\\
The relation $\leq$ is only  reflexive; it is not antisymmetrique (not (An)): $a \leq c$ and $c \leq a$, but $a \not = c$
  and it is not  tranzitive (not (Tr)).\\
Hence,  ${\cal A}$ is a {\bf proper RML algebra},  with (D), (DN).
\end{ex}

{\bf $\bullet$  Proper pre-BCC algebras} (P1)

\begin{ex}\em
Consider the set $A=\{a,b,c,d,1\}$ with the following table of $\ra$: \\ 

\noindent  $\qquad \qquad \qquad \qquad$ 
\begin{tabular}{*{1}{c|}ccccc}
$\ra$ & a & b & c & d & 1 \\
\hline
a & 1 & a & a & b & 1 \\
b & 1 & 1 & 1 & b & 1 \\
c & 1 & 1 & 1 & b & 1 \\
d & 1 & 1 & 1 & 1 & 1 \\
1 & a & b & c & d & 1 
\end{tabular} 

\noindent\\
${\cal A}$ satisfies properties (Re), (M), (L), (B) (hence (*), (Tr)) and (**). It does not satisfy 
 (Ex) for $a,b,d$, (An) for $b,c$,
    (BB) for $d,1,a$, (pi) for $b,a$,
 (D) for $d,a$.\\
Hence, ${\cal A}$ is a {\bf proper pre-BCC algebra},  with (**), without (D).
\end{ex}

 {\bf $\bullet$  Proper aRML algebras} (P9)

\begin{ex}\em
Consider the set $A=\{a,b,c,1\}$ with the following table of $\ra$: \\ 

\noindent  $\qquad \qquad \qquad \qquad$
\begin{tabular}{*{1}{c|}cccc}
$\ra$ & a & b & c & 1  \\
\hline
a & 1 & a & a & 1 \\
b & a & 1 & 1 & 1 \\
c & 1 & a & 1 & 1 \\
1 & a & b & c & 1 
\end{tabular} 

\noindent\\
It verifies (Re), (M), (L), (An), (D). It does not verify   (Ex) for $a,b,b$,
  (pi) for $b,a$,
 (BB), (**) for $a,b,c$,
  (B), (*), (Tr) for $b,c,a$.\\
Hence, ${\cal A}$ is a {\bf proper aRML algebra},  with (D).
\end{ex}

{\bf $\bullet$  Proper BCC algebras} (PO4)

\begin{ex}\em
Consider the set $A=\{a,b,c,d,1\}$ with the following table of $\ra$: \\ 

\noindent  $\qquad \qquad \qquad \qquad$ 
\begin{tabular}{*{1}{c|}ccccc}
$\ra$ & a & b & c & d & 1 \\
\hline
a & 1 & a & a & b & 1 \\
b & 1 & 1 & a & a & 1 \\
c & 1 & a & 1 & a & 1 \\
d & 1 & 1 & a & 1 & 1 \\
1 & a & b & c & d & 1 
\end{tabular} 

\noindent\\
${\cal A}$ satisfies properties (Re), (M), (L), (An), (B) (hence  (*), (Tr)) and (**). It does not satisfy 
 (Ex) for $a,c,d$,
  (pi) for $b,a$,
  (BB) for $d,1,c$,
 (D) for $d,c$.\\
Hence, ${\cal A}$ is a {\bf proper BCC algebra},  with (**), without (D).
\end{ex}

             \subsubsection{ Examples of RML algebras  with condition (Ex): \\BE, pre-BCK, aBE, BCK }
${}$

{\bf $\bullet$  Proper BE algebras} (PO6)

\begin{ex}\em

Consider the set $A=\{a,b,c,1\}$ with the following table of $\ra$: \\ 

\noindent $\qquad \qquad \qquad \qquad$ 
\begin{tabular}{*{1}{c|}cccc}
$\ra$ & a & b & c & 1  \\
\hline
a & 1 & a & 1 & 1 \\
b & 1 & 1 & a & 1 \\
c & 1 & a & 1 & 1 \\
1 & a & b & c & 1 
\end{tabular} 

\noindent\\
${\cal A}$ satisfies (Re), (M), (L), (Ex) (hence (D)). It does not satisfy:
 (pi)  for $b,a$,    (BB), (**)  for $c,b,a$,   (B), (*), (Tr)  for $b,a,c$, (An) for $a,c$.\\
Hence, ${\cal A}$ is a {\bf proper BE algebra}.
\end{ex}

\begin{ex}\em
Consider the set $A=\{0, a,b,c, 1\}$ with the following table of $\ra$: \\ 

\noindent  $\qquad \qquad \qquad \qquad$ 
\begin{tabular}{*{1}{c|}ccccc}
$\ra$ & 0 & a & b & c & 1 \\
\hline
0 & 1 & 1 & 1 & 1 & 1  \\
a & c & 1 & 1 & 1 & 1 \\
b & b & {\bf 1} & 1 & 1 & 1  \\
c &  a & {\bf a} & {\bf 1} & 1 &  1  \\
1 & 0 & a & b & c & 1
\end{tabular} 

\noindent\\
Then the bounded algebra $ {\cal A}=(A, \ra, 0,1)$ verifies properties (Re), (M), (L), (Ex) (hence (D)) and (DN).
 It does not verify   (An), (BB), (B), (*), (**), (Tr),  (pi).\\
Hence,  ${\cal A}$ is a {\bf proper BE algebra}, with (DN).
\end{ex}

{\bf $\bullet$  Proper pre-BCK algebras} (PO7)

\begin{ex}\em

Consider the set $A=\{ a,b,c, 1\}$ with the following table of $\ra$: \\ 

\noindent  $\qquad \qquad \qquad \qquad$
\begin{tabular}{*{1}{c|}cccc}
$\ra$ & a & b & c & 1  \\
\hline
a & 1 & a & a & 1 \\
b & 1 & 1 & 1 & 1 \\
c & 1 & 1 & 1 & 1 \\
1 & a & b & c & 1 
\end{tabular} 

\noindent\\
Properties (Re), (M), (L), (Ex) (hence (D)), (*) (hence (B), (BB),  (**), (Tr)) are satisfied.
 It does not satisfy  (An) for $b,c$;   (pi) for $b,a$.\\
 Hence, ${\cal A}$ is a {\bf proper pre-BCK algebra}.
\end{ex}
                         
\begin{ex}\em

Consider the set $A=\{0,a, b, c, 1\}$ with the following table of $\ra$: \\ 

\noindent  $\qquad \qquad \qquad \qquad$
\begin{tabular}{*{1}{c|}ccccc}
$\ra$ & 0 & a & b & c & 1 \\
\hline
0 & 1 & 1 & 1 & 1 & 1  \\
a & c & 1 & 1 & 1 & 1 \\
b & b & {\bf 1} & 1 & 1 & 1  \\
c &  a & {\bf 1} & {\bf 1} & 1 &  1  \\
1 & 0 & a & b & c & 1
\end{tabular} 

\noindent\\ 
Then the bounded algebra $ {\cal A}=(A, \ra,0, 1)$ verifies properties (Re), (M), (L), (Ex) (hence (D)),
 (*)  (hence (B), (BB), (**), (Tr)) and (DN).  
It does not verify properties  (An),  (pi).\\
 Hence,  ${\cal A}$ is a {\bf proper  pre-BCK algebra},  with (DN).
\end{ex}

{\bf $\bullet$  Proper aBE algebras} (P2)

\begin{ex}\em
 Consider the set $A=\{ a,b,c, 1\}$ with the following table of $\ra$: \\ 
  
\noindent  $\qquad \qquad \qquad \qquad$
\begin{tabular}{*{1}{c|}cccc}
$\ra$ & a & b & c & 1  \\
\hline
a & 1 & a & 1 & 1 \\
b & 1 & 1 & c & 1 \\
c & a & b & 1 & 1 \\
1 & a & b & c & 1 
\end{tabular} 

\noindent\\
Properties (Re), (M), (L), (An)  and (Ex) (hence (D)) are satisfied. 
It does not satisfy:   (pi)  for $b,a$,
  (BB) and  (**) for $c,b,a$,    (B), (*), (Tr) for $b,a,c$.\\
Hence, ${\cal A}$ is a {\bf proper aBE algebra}.
\end{ex}

{\bf $\bullet$  Proper BCK algebras } (PO2)

\begin{ex}\em
Consider the set $A=\{a,b,c, 1\}$ with the following table of $\ra$: \\ 

\noindent  $\qquad \qquad \qquad \qquad$
\begin{tabular}{*{1}{c|}cccc}
$\ra$  & a & b & c & 1 \\
\hline
a & 1 & a & a & 1 \\
b & 1 & 1 & a & 1 \\
c & 1 & a & 1 & 1 \\
1 & a & b & c & 1 
\end{tabular} 

\noindent\\
Properties  (Re), (M), (L), (Ex) (hence (D)), (An), (BB) (hence (B), (*), (**), (Tr)) are satisfied.
 It does not satisfy
   (pi)  for $b,a$.\\
Hence, ${\cal A}$ is a {\bf proper BCK algebra}.
\end{ex}

One can find  many examples of proper BCK algebras in \cite{book}.

                                  \section{Examples of  the other  twenty two   new RM and RML algebras (P10) - (P31)}

We present examples of proper algebras.

              \subsection{Examples of RM and RML algebras without (Ex)}

\subsubsection{RM algebras without (Ex):\\ tRM, *RM, RM**, *RM**, pre-BBBZ,\\ oRM, *aRM, aRM**, *aRM**}
${}$

{\bf $\bullet$  Proper tRM algebras} (P10)

\begin{ex}\em\label{tRM2}
Consider the set $A=\{ a,b, c,1\}$ with the following table of $\ra$: \\ 

\noindent  $\qquad \qquad \qquad \qquad$
\begin{tabular}{*{1}{c|}cccc}
$\ra$ & a & b & c & 1  \\
\hline
a & 1 & b & 1 & b  \\
b & b & 1 & a & c \\
c & 1 & c & 1 & c  \\
1 & a & b & c & 1\
\end{tabular} 

\noindent\\ 
Then the algebra $ {\cal A}=(A, \ra, 1)$ verifies properties (Re), (M), (Tr).
 It does not verify: (Ex) for $a,b,a$, (An) for $a,c$, (L) for $x=b$, 
(**), (BB) for $b,a,c$, (*), (B) for $b,a,c$, (D) for $b,a$. \\
  Hence,  ${\cal A}$ is a {\bf proper tRM algebra}, without (D).
\end{ex}

{\bf $\bullet$  Proper *RM algebras} (P11)

\begin{ex}\em
Consider the set $A=\{0, a,b, 1\}$ with the following table of $\ra$: \\ 

\noindent  $\qquad \qquad \qquad \qquad$
\begin{tabular}{*{1}{c|}cccc}
$\ra$ & 0 & a & b & 1  \\
\hline
0 & 1 & 1 & 1 & 1  \\
a & 1 & 1 & 1 & 1 \\
b & 0 & a & 1 & a  \\
1 & 0 & a & b & 1 
\end{tabular} 

\noindent\\ 
Then the algebra $ {\cal A}=(A, \ra, 1)$ verifies properties  (Re), (M), (*) (hence  (Tr)).
 It does not verify properties  (Ex), (An), (L),  (B), (**) (hence (BB)),  (D). \\
  Hence, ${\cal A}$ is a {\bf proper *RM algebra}, without (D).
\end{ex}

{\bf $\bullet$  Proper RM** algebras} (P12)

\begin{ex}\em
Consider the set $A=\{a, b,c, 1\}$ with the following table of $\ra$: \\ 

\noindent  $\qquad \qquad \qquad \qquad$
\begin{tabular}{*{1}{c|}cccc}
$\ra$ & a & b & c & 1  \\
\hline
a & 1 & a & b & a  \\
b & b & 1 & 1 & a \\
c & c & 1 & 1 & a  \\
1 & a & b & c & 1 
\end{tabular} 

\noindent\\ 
Then the algebra $ {\cal A}=(A, \ra, 1)$ verifies properties (Re), (M), (**) (hence (Tr)).
 It does not verify: (Ex) for $a,b,b$, (An) for $b,c$, (L) for $x=a$,  (BB) for $a,a,c$, (B) for $a,b,a$, (*) for $a,b,c$, (D) for $c,a$.\\
  Hence, ${\cal A}$  is a  {\bf proper RM** algebra}, without (D).
\end{ex}

{\bf $\bullet$  Proper *RM** algebras} (P13)

\begin{ex}\em
Consider the set $A=\{a, b,c, 1\}$ with the following table of $\ra$: \\ 

\noindent  $\qquad \qquad \qquad \qquad$
\begin{tabular}{*{1}{c|}cccc}
$\ra$ & a & b & c & 1  \\
\hline
a & 1 & 1 & c & c  \\
b & 1 & 1 & c & c \\
c & a & b & 1 & 1  \\
1 & a & b & c & 1 
\end{tabular} 

\noindent\\ 
Then the algebra $ {\cal A}=(A, \ra, 1)$ verifies properties (Re), (M), (*), (**) (hence (Tr)).
 It does not verify properties (Ex), (An), (L), (BB), (B),   (D). \\
  Hence,  ${\cal A}$ is  a  {\bf proper *RM** algebra}, without (D).
\end{ex}

{\bf $\bullet$  Proper pre-BBBZ algebras} (P14)

\begin{ex}\em
Consider the set $A=\{a, b, 1\}$ with the following table of $\ra$: \\ 

\noindent  $\qquad \qquad \qquad \qquad$
\begin{tabular}{*{1}{c|}ccc}
$\ra$ & a & b& 1  \\
\hline
a& 1& 1 & a   \\
b& 1 & 1& b \\
1& a& b& 1 
\end{tabular} 

\noindent\\
Properties (Re), (M),  (B), (BB) (hence  (*), (**), (Tr), (D))  are satisfied.
(Ex) is not satisfied for $a,b,a$,  (An) for $a,b$, (L) for $x=a$.\\
Hence, ${\cal A}$ is a {\bf proper pre-BBBZ algebra}, with minimum number of elements, three.
\end{ex}

\begin{ex}\em
Consider the set $A=\{a, b,c, 1\}$ with the following table of $\ra$: \\ 

\noindent  $\qquad \qquad \qquad \qquad$
\begin{tabular}{*{1}{c|}cccc}
$\ra$ & a & b & c & 1  \\
\hline
a & 1 & a & a & b  \\
b & b & 1 & 1 & a \\
c & c & 1 & 1 & a  \\
1 & a & b & c & 1 
\end{tabular} 

\noindent\\
Properties (Re), (M),  (B), (BB) (hence  (*), (**), (Tr), (D)) are satisfied.
 (Ex) is not satisfied for $a,c,b$, (An) for $b,c$, (L) for $x=b$.\\
Hence, ${\cal A}$ is a {\bf proper pre-BBBZ algebra},  one of the very many pre-BBBZ algebras with four elements.
\end{ex}


{\bf $\bullet$  Proper oRM algebras} (P15)

\begin{ex}\em
Consider the set $A=\{a, b, 1\}$ with the following table of $\ra$: \\ 

\noindent  $\qquad \qquad \qquad \qquad$
\begin{tabular}{*{1}{c|}ccc}
$\ra$ & a & b& 1  \\
\hline
a& 1& a & a   \\
b& 1 & 1& b \\
1& a& b& 1 
\end{tabular} 

\noindent\\
Properties (Re), (M), (An), (Tr), (D) are satisfied. ${\cal A}$ does not satisfy  (L),
  (Ex),  (**) (hence (BB)), (*) (hence (B)).\\
Hence, ${\cal A}$ is a {\bf proper oRM algebra}, with (D).
\end{ex}

\begin{ex}\em
Consider the set $A=\{a, b,c, 1\}$ with the following table of $\ra$: \\ 

\noindent  $\qquad \qquad \qquad \qquad$
\begin{tabular}{*{1}{c|}cccc}
$\ra$ & a & b & c & 1  \\
\hline
a & 1 & a & a & a  \\
b & a & 1 & a & a \\
c & b & 1 & 1 & 1 \\
1 & a & b & c & 1 
\end{tabular} 

\noindent\\
Properties (Re), (M), (An), (Tr) are satisfied. ${\cal A}$ does not satisfy  (L) for $x=a$;
  (Ex) for $a,b,a$;  (BB), (**) for $b,c,1$;  (B), (*) for $b,c,b$;  (D) for $a,b$.\\
Hence, ${\cal A}$ is a {\bf proper oRM algebra}, without (D).
\end{ex}

{\bf $\bullet$  Proper *aRM algebras} (P16)
\begin{ex}\em
Consider the set $A=\{a, b,c, 1\}$ with the following table of $\ra$: \\ 

\noindent  $\qquad \qquad \qquad \qquad$
\begin{tabular}{*{1}{c|}ccc}
$\ra$ & a & b& 1  \\
\hline
a& 1& a & a   \\
b& b & 1& 1 \\
1& a& b& 1 
\end{tabular} 

\noindent\\
Properties (Re), (M), (An), (*) (hence  (Tr)), (D) are satisfied. ${\cal A}$ does not satisfy (L), (Ex),
 (BB), (**),  (B).\\
Hence, ${\cal A}$ is a  {\bf proper *aRM algebra}, with (D).
\end{ex}

\begin{ex}\em
Consider the set $A=\{a, b,c, 1\}$ with the following table of $\ra$: \\ 

\noindent $\qquad \qquad \qquad \qquad$ 
\begin{tabular}{*{1}{c|}cccc}
$\ra$ & a & b & c & 1  \\
\hline
a & 1 & a & a & a  \\
b & a & 1 & a & a \\
c & a & a & 1 & 1 \\
1 & a & b & c & 1 
\end{tabular} 

\noindent\\
 Properties (Re), (M), (An), (*) (hence (Tr)) are satisfied.
 ${\cal A}$ does not satisfy  (L) for $x=a$,   (Ex) for $a,b,a$,
 (BB), (**) for $b,c,1$,  (B) for $a,b,c$,  (D) for $a,b$.\\
Hence, ${\cal A}$ is a  {\bf proper *aRM algebra}, without (D).
\end{ex}

{\bf $\bullet$  Proper aRM** algebras} (P17)

\begin{ex}\em
Consider the set $A=\{a, b,c, 1\}$ with the following table of $\ra$: \\ 

\noindent $\qquad \qquad \qquad \qquad$ 
\begin{tabular}{*{1}{c|}cccc}
$\ra$ & a & b & c & 1  \\
\hline
a & 1 & a & a & a  \\
b & a & 1 & a & b \\
c & a & b & 1 & 1 \\
1 & a & b & c & 1 
\end{tabular}  

\noindent\\
Properties (Re), (M), (An), (**) (hence  (Tr)) are satisfied. ${\cal A}$ does not satisfy (L) for $x=a$,
(Ex) for $a,b,a$,   (BB) for $a,b,c$, (B), (*) for $b,c,1$,
  (D) for $a,b$.\\
Hence, ${\cal A}$ is a {\bf proper aRM** algebra}, without (D).
\end{ex}

{\bf $\bullet$  Proper *aRM** algebras} (P18)

\begin{ex}\em

Consider the set $A=\{a, b,c, 1\}$ with the following table of $\ra$: \\ 

\noindent $\qquad \qquad \qquad \qquad$ 
\begin{tabular}{*{1}{c|}cccc}
$\ra$ & a & b & c & 1  \\
\hline
a & 1 & a & a & a  \\
b & a & 1 & a & a \\
c & a & a & 1 & a \\
1 & a & b & c & 1 
\end{tabular}  

\noindent\\
Properties (Re), (M), (An), (*), (**) (hence  (Tr)) are satisfied. ${\cal A}$ does not satisfy (L) for $x=a$,
  (Ex) for $a,b,a$,  (BB), (B) for $a,b,c$,   (D) for $a,b$.\\
Hence, ${\cal A}$ is a {\bf proper *aRM** algebra}, without (D).
\end{ex}

      \subsubsection{RML algebras without (Ex):\\ tRML, *RML, RML**, *RML**, pre-BBBCC,\\  oRML, *aRML, aRML**, *aRML**}
${}$


{\bf $\bullet$  Proper tRML algebras} (P19)

\begin{ex}\em
Consider the set $A=\{a, b,c, 1\}$ with the following table of $\ra$: \\ 

\noindent $\qquad \qquad \qquad \qquad$ 
\begin{tabular}{*{1}{c|}cccc}
$\ra$ & a & b & c & 1  \\
\hline
a & 1 & a & b & 1  \\
b & a & 1 & 1 & 1 \\
c & a & 1 & 1 & 1 \\
1 & a & b & c & 1 
\end{tabular}  

\noindent\\
Properties (Re), (M),  (L), (Tr), (D) are satisfied. ${\cal A}$ does not satisfy  (An) for $b,c$;   (Ex) for $a,b,b$;
   (pi) for $b,a$;  (BB), (**) for $b,a,1$; (B), (*) for $a,b,c$.\\
Hence, ${\cal A}$ is a {\bf proper tRML algebra}, with (D).
\end{ex}

\begin{ex}\em
Consider the set $A=\{0, a,b,c, 1\}$ with the following table of $\ra$: \\ 

\noindent $\qquad \qquad \qquad \qquad$ 
\begin{tabular}{*{1}{c|}ccccc}
$\ra$ & 0 & a & b & c & 1 \\
\hline
0 & 1 & 1 & 1 & 1 & 1  \\
a & c & 1 & 1 & 1 & 1 \\
b & b & {\bf a} & 1 & 1 & 1  \\
c &  a & {\bf a} & {\bf 1} & 1 &  1  \\
1 & 0 & a & b & c & 1
\end{tabular} 

\noindent\\
Then the bounded algebra $ {\cal A}=(A, \ra,0, 1)$ verifies properties (Re), (M), (L), (Tr) and (D),  (DN).
 It does not verify properties (An),  (Ex), (BB), (B), (*), (**),  (pi).\\
  Hence,  ${\cal A}$ is a {\bf proper tRML algebra}, with (DN),  and with (D).
\end{ex}

\begin{ex}\em
Consider the set $A=\{0, a,b,c, 1\}$ with the following table of $\ra$: \\ 

\noindent  $\qquad \qquad \qquad \qquad$
\begin{tabular}{*{1}{c|}ccccc}
$\ra$ & 0 & a & b & c & 1 \\
\hline
0 & 1 & 1 & 1 & 1 & 1  \\
a & c & 1 & 1 & 1 & 1 \\
b & b & {\bf a} & 1 & 1 & 1  \\
c &  a & {\bf b} & {\bf 1} & 1 &  1  \\
1 & 0 & a & b & c & 1
\end{tabular} 

\noindent\\
Then the bounded algebra $ {\cal A}=(A, \ra,0, 1)$ verifies properties (Re), (M), (L), (Tr), and (DN).
 It does not verify  (An) for $b,c$, (Ex) for $b,c,0$, (BB), (**) for $0,c,b$,  (B), (*) for $b,0,a$,
 (D) for $a,c$, (pi) for $0,a$.\\ 
  Hence,  ${\cal A}$ is a {\bf proper tRML algebra}, with (DN), without (D).
\end{ex}

{\bf $\bullet$  Proper *RML algebras} (P20)

\begin{ex}\em
Consider the set $A=\{a, b,c, 1\}$ with the following table of $\ra$: \\ 

\noindent $\qquad \qquad \qquad \qquad$ 
\begin{tabular}{*{1}{c|}cccc}
$\ra$ & a & b & c & 1  \\
\hline
a & 1 & a & a & 1  \\
b & a & 1 & 1 & 1 \\
c & a & 1 & 1 & 1 \\
1 & a & b & c & 1 
\end{tabular} 

\noindent\\
Properties (Re), (M), (L), (*) (hence  (Tr)), (D) are satisfied.
 ${\cal A}$ does not satisfy (An) for $b,c$;  (Ex) for $a,b,b$;
 (pi) for $b,a$;   (**), (BB) for $b,a,1$;  (B) for $a,1,b$.\\
Hence, ${\cal A}$ is a {\bf proper *RML algebra}, with (D).
\end{ex}

\begin{ex}\em
Consider the set $A=\{0, a,b,c, 1\}$ with the following table of $\ra$: \\ 

\noindent $\qquad \qquad \qquad \qquad$ 
\begin{tabular}{*{1}{c|}ccccc}
$\ra$ & 0 & a & b & c & 1 \\
\hline
0 & 1 & 1 & 1 & 1 & 1  \\
a & c & 1 & 1 & 1 & 1 \\
b & b & {\bf c} & 1 & 1 & 1  \\
c &  a & {\bf c} & {\bf 1} & 1 &  1  \\
1 & 0 & a & b & c & 1
\end{tabular} 

\noindent\\
Then the bounded algebra $ {\cal A}=(A, \ra,0, 1)$ verifies properties (Re), (M), (L), (*) (hence  (Tr))  and (D), and (DN).
 It does not verify properties  (An), (Ex), (BB), (B), (**), (pi).\\
  Hence, ${\cal A}$ is a {\bf proper *RML algebra}, with (DN), with (D).
\end{ex}

\begin{ex}\em
Consider the set $A=\{0, a,b,c, 1\}$ with the following table of $\ra$: \\ 

\noindent  $\qquad \qquad \qquad \qquad$
\begin{tabular}{*{1}{c|}ccccc}
$\ra$ & 0 & a & b & c & 1 \\
\hline
0 & 1 & 1 & 1 & 1 & 1  \\
a & c & 1 & 1 & 1 & 1 \\
b & b & {\bf c} & 1 & 1 & 1  \\
c &  a & {\bf a} & {\bf 1} & 1 &  1  \\
1 & 0 & a & b & c & 1
\end{tabular} 

\noindent\\
Then the bounded algebra $ {\cal A}=(A, \ra,0, 1)$ verifies properties (Re), (M), (L), (*) (hence (Tr)), and (DN).
 It does not verify properties (An), (Ex), (BB), (B),  (**), (D), (pi).\\
  Hence,  ${\cal A}$ is a {\bf proper *RML algebra}, with (DN),  without (D).
\end{ex}

{\bf $\bullet$ Proper RML** algebras} (P21)

\begin{ex}\em
Consider the set $A=\{a, b,c,d, 1\}$ with the following table of $\ra$: \\ 

\noindent $\qquad \qquad \qquad \qquad$  
\begin{tabular}{*{1}{c|}ccccc}
$\ra$ & a & b & c & d & 1 \\
\hline
a & 1 & a  & a  & b  &  1 \\
b & 1  & 1 & a  & b  &  1 \\
c &  1 &  1 & 1 &  1 & 1  \\
d & 1  &  1 &  1 & 1 &  1 \\
1 & a & b & c & d & 1
\end{tabular} 

\noindent\\
Properties (Re), (M), (L), (**) (hence  (Tr)) are satisfied. ${\cal A}$ does not satisfy 
(An) for $c,d$, (Ex) for $a,b,d$,  (pi) for $b,a$, (BB) for $d,a,c$, (B), (*) for $a,c,d$, (D) for $d,a$.\\
Hence, ${\cal A}$ is a {\bf proper RML** algebra}, without (D).
\end{ex}

{\bf $\bullet$  Proper *RML** algebras} (P22)

\begin{ex}\em
Consider the set $A=\{a, b,c,d, 1\}$ with the following table of $\ra$: \\ 

\noindent  $\qquad \qquad \qquad \qquad$ 
\begin{tabular}{*{1}{c|}ccccc}
$\ra$ & a & b & c & d & 1 \\
\hline
a & 1 & a  & a  & d  &  1 \\
b & 1  & 1 & 1  & a  &  1 \\
c &  1 &  1 & 1 &  a & 1  \\
d & 1  &  a &  a & 1 &  1 \\
1 & a & b & c & d & 1
\end{tabular} 

\noindent\\
Properties (Re), (M), (L), (*), (**) (hence  (Tr)) are satisfied.
 ${\cal A}$ does not satisfy (An), (Ex), (pi), (BB), (B),  (D).\\
Hence, ${\cal A}$ is a {\bf proper *RML** algebra}, without (D).
\end{ex}

{\bf $\bullet$  Proper pre-BBBCC algebras} (P23)

\begin{ex}\em
Consider the set $A=\{a, b,c,d, 1\}$ with the following table of $\ra$: \\ 

\noindent $\qquad \qquad \qquad \qquad$  
\begin{tabular}{*{1}{c|}ccccc}
$\ra$ & a & b & c & d & 1 \\
\hline
a & 1 & a  & c  & c  &  1 \\
b & 1  & 1 & d  & c  &  1 \\
c &  a &  b & 1 &  1 & 1  \\
d & a  &  b &  1 & 1 &  1 \\
1 & a & b & c & d & 1
\end{tabular} 

\noindent\\
Properties (Re), (M), (L), (BB), (B) (hence  (**), (*), (Tr), (D)) are satisfied. 
${\cal A}$ does not satisfy (An) for $x=c$, $y=d$,  (Ex) for $a,b,c$,  (pi) for $b,a$.\\
Hence, ${\cal A}$ is a {\bf proper pre-BBBCC algebra}.
\end{ex}


{\bf $\bullet$  Proper oRML algebras} (P24)

\begin{ex}\em
Consider the set $A=\{a, b,c, 1\}$ with the following table of $\ra$: \\ 

\noindent  $\qquad \qquad \qquad \qquad$
\begin{tabular}{*{1}{c|}cccc}
$\ra$ & a & b & c & 1  \\
\hline
a & 1 & a & a & 1  \\
b & a & 1 & b & 1 \\
c & 1 & a & 1 & 1 \\
1 & a & b & c & 1 
\end{tabular} 

\noindent\\
Properties (Re), (M), (L), (An), (Tr), (D) are satisfied. ${\cal A}$ does not satisfy:
 (Ex) for $a,b,b$;  (pi) for $b,a$;   (**), (BB)  for $b,c,1$;
   (*), (B) for $b,c,a$.\\
Hence, ${\cal A}$ is a {\bf proper oRML algebra}, with (D).
\end{ex}

\begin{ex}\em
Consider the set $A=\{0, a,b,c, 1\}$ with the following table of $\ra$: \\ 

\noindent  $\qquad \qquad \qquad \qquad$ 
\begin{tabular}{*{1}{c|}ccccc}
$\ra$ & 0 & a & b & c & 1 \\
\hline
0 & 1 & 1 & 1 & 1 & 1  \\
a & c & 1 & 1 & 1 & 1 \\
b & b & {\bf a} & 1 & 1 & 1  \\
c &  a & {\bf a} & {\bf a} & 1 &  1  \\
1 & 0 & a & b & c & 1
\end{tabular} 

\noindent\\
Then the bounded algebra $ {\cal A}=(A, \ra,0, 1)$ verifies properties (Re), (M), (L), (An), (Tr) and (D),  and (DN).
 It does not verify properties  (Ex), (BB), (B),  (**), (*),  (pi).\\
  Hence,  ${\cal A}$ is a {\bf proper oRML algebra}, with (DN) and with  (D).
\end{ex}

\begin{ex}\em
Consider the set $A=\{0, a,b,c, 1\}$ with the following table of $\ra$: \\ 

\noindent $\qquad \qquad \qquad \qquad$ 
\begin{tabular}{*{1}{c|}ccccc}
$\ra$ & 0 & a & b & c & 1 \\
\hline
0 & 1 & 1 & 1 & 1 & 1  \\
a & c & 1 & 1 & 1 & 1 \\
b & b & {\bf a} & 1 & 1 & 1  \\
c &  a & {\bf b} & {\bf a} & 1 &  1  \\
1 & 0 & a & b & c & 1
\end{tabular} 

\noindent\\
Then the bounded  algebra $ {\cal A}=(A, \ra,0, 1)$ verifies properties (Re), (M), (L), (An), (Tr),   and (DN).
 It does not verify properties (Ex), (BB), (B),  (**), (*), (D),  (pi).\\
  Hence,  ${\cal A}$ is a {\bf proper oRML algebra},  with (DN), without (D).
\end{ex}

{\bf $\bullet$  Proper *aRML algebras} (P25)

\begin{ex}\em
Consider the set $A=\{a, b,c, 1\}$ with the following table of $\ra$: \\ 

\noindent  $\qquad \qquad \qquad \qquad$
\begin{tabular}{*{1}{c|}cccc}
$\ra$ & a & b & c & 1  \\
\hline
a & 1 & a & a & 1  \\
b & a & 1 & a & 1 \\
c & a & a & 1 & 1 \\
1 & a & b & c & 1 
\end{tabular} 

\noindent\\
Properties (Re), (M), (L), (*) (hence  (Tr)) are satisfied. ${\cal A}$ does not satisfy:   (Ex) for $a,b,b$;
   (pi) for $b,a$;  (BB), (**) for $b,a,1$;  (B) for $a,1,b$;   (D) for $b,c$.\\
Hence, ${\cal A}$ is a {\bf proper *aRML algebra}, without (D).
\end{ex}

\begin{ex}\em
Consider the set $A=\{a,b, 1\}$ with the following table of $\ra$: \\ 

\noindent  $\qquad \qquad \qquad \qquad$
\begin{tabular}{*{1}{c|}ccc}
$\ra$ & a & b & 1   \\
\hline
a & 1 & a & 1  \\
b & a & 1 & 1 \\
1 & a & b & 1  
\end{tabular} 

\noindent\\ 
Then the algebra $ {\cal A}=(A, \ra, 1)$ verifies properties (Re), (M), (L), (An), (*) (hence (Tr)), and (D).  
It does not verify (Ex), (BB), (**), (B), (pi).\\
 Hence,  ${\cal A}$ is a {\bf proper *aRML algebra}, with (D). 
\end{ex}

\begin{ex}\em
Consider the set $A=\{0, a,b,c, 1\}$ with the following table of $\ra$: \\ 

\noindent  $\qquad \qquad \qquad \qquad$
\begin{tabular}{*{1}{c|}ccccc}
$\ra$ & 0 & a & b & c & 1 \\
\hline
0 & 1 & 1 & 1 & 1 & 1  \\
a & c & 1 & 1 & 1 & 1 \\
b & b & {\bf c} & 1 & 1 & 1  \\
c &  a & {\bf a} & {\bf a} & 1 &  1  \\
1 & 0 & a & b & c & 1
\end{tabular} 

\noindent\\
Then the  bounded algebra $ {\cal A}=(A, \ra,0, 1)$ verifies properties (Re), (M), (L), (An), (*) (hence (Tr)),   and (DN).
 It does not verify properties (Ex), (BB), (B),  (**),  (D), (pi).\\
  Hence, ${\cal A}$ is a {\bf proper *aRML algebra},  with (DN), without (D).
\end{ex}

{\bf $\bullet$  Proper aRML** algebras} (P26)

\begin{ex}\em
Consider the set $A=\{a, b,c, 1\}$ with the following table of $\ra$: \\ 

\noindent  $\qquad \qquad \qquad \qquad$
\begin{tabular}{*{1}{c|}cccc}
$\ra$ & a & b & c & 1  \\
\hline
a & 1 & a & c & 1  \\
b & 1 & 1 & 1 & 1 \\
c & a & a & 1 & 1 \\
1 & a & b & c & 1 
\end{tabular} 

\noindent\\
Properties (Re), (M), (L), (An),  (**) (hence  (Tr)) are satisfied. ${\cal A}$ does not satisfy:   (Ex) for $a,c,b$;
 (pi) for $b,a$;   (BB) for $b,1,c$;   (B), (*) for $a,b,c$;
  (D) for $b,c$.\\
Hence, ${\cal A}$ is a {\bf proper aRML** algebra}, without (D).
\end{ex}

\begin{ex}\em
Consider the set $A=\{0, a,b,c, 1\}$ with the following table of $\ra$: \\ 

\noindent  $\qquad \qquad \qquad \qquad$
\begin{tabular}{*{1}{c|}ccccc}
$\ra$ & 0 & a & b & c & 1 \\
\hline
0 & 1 & 1 & 1 & 1 & 1  \\
a & c & 1 & 1 & 1 & 1 \\
b & b & {\bf a} & 1 & 1 & 1  \\
c &  a & {\bf a} & {\bf b} & 1 &  1  \\
1 & 0 & a & b & c & 1
\end{tabular} 

\noindent\\
Then the bounded algebra $ {\cal A}=(A, \ra,0, 1)$ verifies properties (Re), (M), (L), (An), (**) (hence (Tr)),  (D)
  and (DN).
 It does not verify properties (Ex), (BB), (B),  (*),  (pi).\\
  Hence, ${\cal A}$ is a {\bf proper aRML** algebra}, with (DN), with (D).
\end{ex}

{\bf $\bullet$  Proper *aRML** algebras} (P27)

\begin{ex}\em
Consider the set $A=\{0, a,b,c, 1\}$ with the following table of $\ra$: \\ 

\noindent  $\qquad \qquad \qquad \qquad$
\begin{tabular}{*{1}{c|}ccccc}
$\ra$ & 0 & a & b & c & 1 \\
\hline
0 & 1 & 1 & 1 & 1 & 1  \\
a & c & 1 & 1 & 1 & 1 \\
b & b & {\bf b} & 1 & 1 & 1  \\
c &  a & {\bf a} & {\bf c} & 1 &  1  \\
1 & 0 & a & b & c & 1
\end{tabular} 

\noindent\\
Then the bounded  algebra $ {\cal A}=(A, \ra, 0,1)$ verifies properties (Re), (M), (L), (An), (*), (**) (hence (Tr)),  (D)
  and (DN).
 It does not verify properties (Ex), (BB), (B), (pi).\\
  Hence,  ${\cal A}$ is a {\bf proper  *aRML** algebra}, with (DN), 
 with (D).
\end{ex}

\begin{ex}\em
Consider the set $A=\{0, a,b,c, 1\}$ with the following table of $\ra$: \\ 

\noindent $\qquad \qquad \qquad \qquad$  
\begin{tabular}{*{1}{c|}ccccc}
$\ra$ & 0 & a & b & c & 1 \\
\hline
0 & 1 & 1 & 1 & 1 & 1  \\
a & c & 1 & 1 & 1 & 1 \\
b & b & {\bf c} & 1 & 1 & 1  \\
c &  a & {\bf a} & {\bf b} & 1 &  1  \\
1 & 0 & a & b & c & 1
\end{tabular} 

\noindent\\
Then the  bounded algebra $ {\cal A}=(A, \ra,0,1)$ verifies properties (Re), (M), (L), (An), (**), (*) (hence (Tr)),
  and (DN).
 It does not verify properties (Ex), (BB), (B), (D), (pi).\\
Hence, ${\cal A}$ is a  {\bf proper *aRML** algebra},
 with (DN), without (D).
\end{ex}

               \subsubsection{Examples of RM and RML algebras with  (Ex):\\
RME**, BCH**, BE**, aBE**}
${}$

{\bf $\bullet$  Proper RME** algebras} (P28)              

\begin{ex}\em
Consider the set $A=\{a, b,c,d, 1\}$ with the following table of $\ra$: \\ 

\noindent $\qquad \qquad \qquad \qquad$  
\begin{tabular}{*{1}{c|}ccccc}
$\ra$ & a & b & c & d & 1 \\
\hline
a & 1 & a  & a  & a  &  a \\
b & a  & 1 & b  & c  &  1 \\
c &  a &  1 & 1 &  1 & 1  \\
d & a  &  1 &  1 & 1 &  1 \\
1 & a & b & c & d & 1
\end{tabular} 

\noindent\\
Properties  (Re), (M), (Ex) (hence (D)), (**) (hence  (Tr)), (D) are  satisfied. ${\cal A}$ does not satisfy:
  (An) for $c,d$; (BB) for $d,b,c$;   (B), (*) for $b,c,d$.\\
Hence, ${\cal A}$ is a {\bf proper RME** algebra}.
\end{ex}

{\bf $\bullet$  Proper BCH** algebras} (P29)

\begin{ex}\em
Consider the set $A=\{a, b,c,d, 1\}$ with the following table of $\ra$: \\ 

\noindent  $\qquad \qquad \qquad \qquad$ 
\begin{tabular}{*{1}{c|}ccccc}
$\ra$ & a & b & c & d & 1 \\
\hline
a & 1 & a  & a  & a  &  a \\
b & a  & 1 & b  & d  &  1 \\
c &  a &  1 & 1 &  1 & 1  \\
d & a  &  b &  c & 1 &  1 \\
1 & a & b & c & d & 1
\end{tabular} 

\noindent\\
Properties  (Re), (M), (Ex) (hence (D)), (An), (**) (hence (Tr)) are satisfied. ${\cal A}$ does not satisfy:
   (BB) for $d,b,c$;
 (B), (*) for $b,c,d$.\\
Hence, ${\cal A}$ is a {\bf proper BCH** algebra}.
\end{ex}


{\bf $\bullet$  Proper BE** algebras} (P30)

\begin{ex}\em
Consider the set $A=\{a, b,c, 1\}$ with the following table of $\ra$: \\ 

\noindent $\qquad \qquad \qquad \qquad$ 
\begin{tabular}{*{1}{c|}cccc}
$\ra$ & a & b & c & 1  \\
\hline
a & 1 & a & b & 1  \\
b & 1 & 1 & 1 & 1 \\
c & 1 & 1 & 1 & 1 \\
1 & a & b & c & 1 
\end{tabular} 

\noindent\\
Properties (Re), (M), (L), (Ex) (hence (D)), (**) (hence  (Tr)) are satisfied. ${\cal A}$ does not satisfy:
(An) for $b,c$;
   (pi) for $b,a$;   (BB) for $c,a,b$; (B), (*) for $a,b,c$.\\
Hence, ${\cal A}$ is a {\bf proper BE** algebra}.
\end{ex}

{\bf $\bullet$  Proper aBE** algebras} (P31)

\begin{ex}\em
Consider the set $A=\{0, a,b, 1\}$ with the following table of $\ra$: \\ 

\noindent $\qquad \qquad \qquad \qquad$ 
\begin{tabular}{*{1}{c|}cccc}
$\ra$ & 0 & a & b & 1  \\
\hline
0 & 1 & 1 & 1 & 1  \\
a & 0 & 1 & b & 1 \\
b & b & a & 1 & 1  \\
1 & 0 & a & b & 1 \\

\end{tabular} 

\noindent\\ 
Then the algebra $ {\cal A}=(A, \ra, 1)$ verifies properties (Re), (M), (L), (Ex) (hence (D)), (An), (**) (hence (Tr))
 (see \cite{Ahn-So}, Example 3.9). 
It does not verify  condition (BB)  for $x=a$, $y=b$, $z=0$: 
$ b=b \ra 0=y \ra z \not \leq (z \ra x) \ra (y \ra x)=(0 \ra a) \ra (b \ra a)= 
1 \ra a=a$; 
it does not verify condition (pi) for $x=0$, $y=b$: 
$1= b \ra b=b \ra (b \ra 0)=y \ra (y \ra x) \not = y \ra x=b \ra 0= b$.
 The relation $\leq$
is  a lattice order.\\
 Hence,   ${\cal A}$ is a {\bf proper aBE** lattice}.

We shall represent the set $A$ and the  binary relation $\leq$ by the following Hasse diagram:\\

\begin{figure}[htbp]
\begin{center}
\begin{picture}(80,50)(0,0) 

\put(39,-2){\makebox(2,2){$\bullet$}}

\put(40,0){\line(-3,4){15}}
\put(40,0){\line(3,4){15}}
\put(23,18){\makebox(2,2){$\bullet$}}
\put(55,18){\makebox(2,2){$\bullet$}}
 
\put(25,20){\line(3,4){15}}
\put(55,20){\line(-3,4){15}}

\put(39,40){\makebox(2,2){$\bullet$}}

\put(35,-15){\makebox(10,10){0}}
\put(12,17){\makebox(10,10){a}}
\put(60,17){\makebox(10,10){b}}

\put(35,45){\makebox(10,10){1}}

\end{picture} 
 \end{center}
\caption{Proper aBE** lattice   }
 \label{fig: A1}
\end{figure}
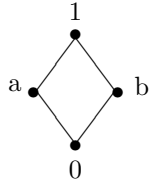
\end{ex}

{\bf Remark:} We have proved that condition (BB) implies (transitivity) (Tr); but  (Tr)  does not imply condition (BB). Indeed, the relation $\leq$ associated with
 above ${\cal A}$ is transitive, but ${\cal A}$ does not verify (BB).

                            \section {Examples of (proper) generalizations of Hilbert algebras}

                \subsection{Examples of pi-RML algebras from Hierachy 2H}

Recall that we have proved that only  RML algebras can have the properties  (pi)  and (pimpl).

              \subsubsection{ Examples of pi-RML algebras  without condition (Ex):\\ pi-RML, pi-pre-BCC, pi-aRML, pi-BCC }
${}$

{\bf $\bullet$  Proper pi-RML algebras} (P8-pi)

\begin{ex}\em
Consider the set $A=\{a,b,c,d,1\}$ with the following table of $\ra$: \\ 

\noindent $\qquad \qquad \qquad \qquad$  
\begin{tabular}{*{1}{c|}ccccc}
$\ra$ & a & b & c & d & 1 \\
\hline
a & 1 & b & b & b & 1 \\
b & a & 1 & a & a & 1 \\
c & a & a & 1 & 1 & 1 \\
d & a & 1 & 1 & 1 & 1 \\
1 & a & b & c & d & 1 
\end{tabular} 

\noindent\\
${\cal A}$ satisfies (Re), (M), (L), (pi). It does not satisfy:  (An) for $c,d$, 
 (Ex) for $a,c,b$, 
   (pimpl) for $b,c,a$,
 (BB), (**) for $b,c,d$,
 (B), (*), (Tr) for $c,d,b$, 
  (D) for $b,c$.\\
Hence, ${\cal A}$ is a {\bf proper pi-RML algebra}, without (D).
\end{ex}

{\bf $\bullet$  Proper pi-pre-BCC algebras} (P1-pi)

\begin{ex}\em

Consider the set $A=\{a,b,c,d,1\}$ with the following table of $\ra$: \\ 

\noindent $\qquad \qquad \qquad \qquad$  
\begin{tabular}{*{1}{c|}ccccc}
$\ra$ & a & b & c & d & 1 \\
\hline
a & 1 & b & b & b & 1 \\
b & a & 1 & c & c & 1 \\
c & a & 1 & 1 & 1 & 1 \\
d & a & 1 & 1 & 1 & 1 \\
1 & a & b & c & d & 1 
\end{tabular} 

\noindent\\
Properties  (Re), (M), (L),   (B) (hence  (*), (Tr), (**)) and (pi) are satisfied. 
${\cal A}$ does not satisfy (An) for $c,d$,  (Ex)  for $a,b,c$,   (pimpl) for $a,b,c$,
  (BB) for $c,b,a$,
D)  for $d,a$.\\
Hence, ${\cal A}$ is a {\bf proper pi-pre-BCC algebra},  without (D).
\end{ex}

{\bf $\bullet$  Proper pi-aRML algebras} (P9-pi)

\begin{ex}\em

Consider the set $A=\{a,b,c,d,1\}$ with the following table of $\ra$: \\ 

\noindent $\qquad \qquad \qquad \qquad$  
\begin{tabular}{*{1}{c|}ccccc}
$\ra$ & a & b & c & d & 1 \\
\hline
a & 1 & b & b & b & 1 \\
b & a & 1 & a & a & 1 \\
c & a & a & 1 & 1 & 1 \\
d & a & 1 & a & 1 & 1 \\
1 & a & b & c & d & 1 
\end{tabular} 

\noindent\\
Properties  (Re), (M), (L), (An), (pi) are satisfied. 
${\cal A}$ does not satisfy: (Ex)  for $a,c,b$,  (B) for $a,1,c$, (BB) for $b,c,d$,  (*) for $c,d,b$, (**) for $b,c,d$,
(Tr) for $c,d,b$,  (pimpl) for $b,c,a$,
  (D)  for $b,c$.\\
Hence, ${\cal A}$ is a {\bf proper pi-aRML algebra},   without (D).
\end{ex}

{\bf $\bullet$  Proper pi-BCC algebras}  (PO4-pi)

\begin{ex}\em
Consider the set $A=\{ a,b,c, 1\}$ with the following table of $\ra$: \\ 

\noindent  $\qquad \qquad \qquad \qquad$ 
\begin{tabular}{*{1}{c|}cccc}
$\ra$ & a & b & c & 1  \\
\hline
a & 1 & b & b & 1  \\
b & a & 1 & c & 1 \\
c & 1 & 1 & 1 & 1  \\
1 & a & b & c & 1 
\end{tabular} $\qquad$ 

\noindent\\ 
Then the algebra $ {\cal A}=(A, \ra, 1)$ verifies properties (Re), (M), (L), (An), (B) (hence (*), (Tr), (**)) and  (pi).
 It does not verify:  (Ex) for $a,b,c$,   (BB) for $c,b,a$,  (pimpl) for $a,b,c$,  (D) for $c,a$.\\
Hence, {\bf ${\cal A}$ is a  proper pi-BCC algebra},  without (D).
 \end{ex}

             \subsubsection{ Examples of pi-RML algebras  with condition (Ex): \\pi-BE, pi-pre-BCK and pimpl-pre-BCK,
 pi-aBE, pi-BCK = pimpl-BCK = Hilbert algebras }
${}$

{\bf $\bullet$  Proper pi-BE algebras} (PO6-pi)


\begin{ex}\em
Consider the set $A=\{ a,b,c, 1\}$ with the following table of $\ra$: \\ 

\noindent $\qquad \qquad \qquad \qquad$  
\begin{tabular}{*{1}{c|}cccc}
$\ra$ & a & b & c & 1  \\
\hline
a & 1 & 1 & c & 1  \\
b & 1 & 1 & 1 & 1 \\
c & a & b & 1 & 1  \\
1 & a & b & c & 1 
\end{tabular} $\qquad$ 

\noindent\\ 
Then the algebra $ {\cal A}=(A, \ra, 1)$ verifies properties (Re), (M), (L), (Ex) (hence (D)) and  (pi).
 It does not verify  (An) for $a,b$,  (BB), (**)  for 
$x=c, \; y=a, \; z=b$; (B), (*), (Tr) for $a,b,c$; (pimpl) 
for $x=a, \; y=b, \; z=c$.\\
Hence,  ${\cal A}$ is a  {\bf proper pi-BE algebra}.

We shall represent the set $A$ and the  binary relation $\leq$ by the following Hasse-type diagram:\\

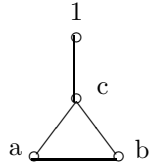
\begin{figure}[htbp]
\begin{center}
\begin{picture}(50,30)(0,20) 


\put(23,18){\makebox(2,2){$\circ$}}
\put(55,18){\makebox(2,2){$\circ$}}
 
           \put(23,18){\line(3,0){32} }
\put(25,20){\line(3,4){15}}
\put(55,20){\line(-3,4){15}}

\put(39,40){\makebox(2,2){$\circ$}}
                             \put(39,40){\line(0,4){25}}
\put(39,63){\makebox(2,2){$\circ$}}
\put(39,73){\makebox(2,2){1}}

\put(12,17){\makebox(10,10){a}}
\put(60,17){\makebox(10,10){b}}

\put(45,40){\makebox(10,10){c}}

\end{picture}
 \end{center}
\caption{Proper pi-BE algebra }
 \label{fig: A5}
\end{figure}
\end{ex}

{\bf $\bullet$  Proper pi-pre-BCK  (PO7-pi) and pimpl-pre-BCK (PO7-pimpl) algebras} 

\begin{ex}\em

Consider the set $A=\{a, b, c, d,1\}$ with the following table of $\ra$: \\ 

\noindent $\qquad \qquad \qquad \qquad$  
\begin{tabular}{*{1}{c|}ccccc}
$\ra$ & a & b & c & d & 1 \\
\hline
a & 1 & b & b & d & 1 \\
b & a & 1 & 1 & d & 1 \\
c & a & 1 & 1 & d & 1 \\
d & 1 & b & c & 1 & 1 \\
1 & a & b & c & d & 1 
\end{tabular} 

\noindent\\
Then the algebra $ {\cal A}=(A, \ra, 1)$ verifies properties (Re), (M), (L), (Ex) (hence (D)),
 (*)  (hence (B), (BB),  (**),  (Tr)) and (pi).  
It does not verify properties  (An) for $x=b$, $y=c$ and  (pimpl) for $x=d$, $y=a$, $z=c$.\\
 Hence,  ${\cal A}$ is a proper  pre-BCK algebra,  verifying (pi) and not verifying (pimpl),
Hence, it is a {\bf proper pi-pre-BCK algebra}.
\end{ex}

\begin{ex}\em (This example is taken from \cite{Busneag-Rudeanu}, from the proof of Proposition 1.1)

Consider the set $A_2=\{c,d, 1\}$ with the following table of $\ra$: \\ 

\noindent $\qquad \qquad \qquad \qquad$ 
\begin{tabular}{*{1}{c|}ccc}
$\ra$ & c & d & 1   \\
\hline
c & 1 & 1 & 1  \\
d & 1 & 1 & 1 \\
1 & c & d & 1  
\end{tabular} 

\noindent\\ 
Then the algebra $ {\cal A}=(A, \ra, 1)$ verifies properties 
(Re), (M), (L), (Ex) (hence (D)), (*)  (hence (B), (BB), (**),  (Tr)) and  (pimpl) (hence  (pi)).  It does not verify (An).
 $\leq$ is only a pre-order (\cite{Busneag-Rudeanu}, the proof of  Proposition 1.1.). \\
 Hence, ${\cal A}$ is a {\bf proper pimpl-pre-BCK algebra},
  with the  minimum number of elemements, three.\\
We shall represent the set $A$ and the  binary relation $\leq$ by the following Hasse-type diagram:\\

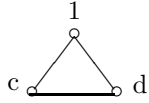
\begin{figure}[htbp]
\begin{center}
\begin{picture}(50,30)(0,20) 


\put(23,18){\makebox(2,2){$\circ$}}
\put(55,18){\makebox(2,2){$\circ$}}
 
           \put(23,18){\line(3,0){32} }
\put(25,20){\line(3,4){15}}
\put(55,20){\line(-3,4){15}}

\put(39,40){\makebox(2,2){$\circ$}}

\put(12,17){\makebox(10,10){c}}
\put(60,17){\makebox(10,10){d}}

\put(35,45){\makebox(10,10){1}}

\end{picture}
 \end{center}
\caption{The  proper pimpl-pre-BCK algebra   with 3 elements  }
 \label{fig: A2}
\end{figure}
\end{ex}


\begin{ex}\em
Consider the set $A=\{ a,b,c, 1\}$ with the following table of $\ra$: \\ 

\noindent $\qquad \qquad \qquad \qquad$
\begin{tabular}{*{1}{c|}cccc}
$\ra$ & a & b & c & 1  \\
\hline
a & 1 & 1 & 1 & 1  \\
b & 1 & 1 & 1 & 1 \\
c & 1 & 1 & 1 & 1  \\
1 & a & b & c & 1 
\end{tabular} 

\noindent\\ 
Then the algebra $ {\cal A}=(A, \ra, 1)$ verifies properties  (Re), (M), (L), (Ex) (hence (D)), (*) 
(hence (B), (BB), (**),  (Tr)) and  (pimpl) (hence  (pi)).
It does not verify (An).
 The relation $\leq$
is a pre-order relation.\\
Hence,  ${\cal A}$ is a {\bf proper pimpl-pre-BCK algebra}, with four elements.

We shall represent the set $A$ and the  binary relation $\leq$ by the following Hasse-type diagram:\\

\begin{figure}[htbp]
\begin{center}
\begin{picture}(30,10)(0,20) 
        
\put(10,20){\makebox(2,2){$\circ$}}
           
       \put(10, 7){\makebox(2,2){b}}
        \put(10,20){\line(0,4){20}}

\put(-5,20){\makebox(2,2){$\circ$}}
\put(-5,20){\line (4,0) {30}}
\put(25,20){\makebox(2,2){$\circ$}}

\put(-5,20){\line(3,4){15}}
\put(25,20){\line(-3,4){15}}

\put(10,37){\makebox(2,2){$\circ$}}

\put(-18,7){\makebox(10,10){a}}
\put(30,7){\makebox(10,10){c}}

\put(6,42){\makebox(10,10){1}}
\end{picture}
 \end{center}
\caption{The  proper pimpl-pre-BCK algebra  with 4 elements } \label{fig: A4}
\end{figure}
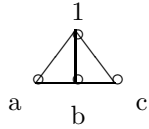
\end{ex}

{\bf $\bullet$ Proper pi-aBE algebras} (P2-pi)

\begin{ex}\em
Consider the set $A=\{ a,b,c, 1\}$ with the following table of $\ra$: \\ 

\noindent  $\qquad \qquad \qquad \qquad$
\begin{tabular}{*{1}{c|}cccc}
$\ra$ & a & b & c & 1  \\
\hline
a & 1 & 1 & c & 1  \\
b & a & 1 & 1 & 1 \\
c & a & b & 1 & 1  \\
1 & a & b & c & 1 
\end{tabular} 

\noindent\\ 
Then the algebra $ {\cal A}=(A, \ra, 1)$ verifies properties (Re), (M), (L), (Ex) (hence (D)), (An) and  (pi).
 It does not verify  (BB), (**)  for 
$x=c, \; y=a, \; z=b$, (B), (*), (Tr) for $a,b,c$,
 (pimpl) for $x=a, \; y=b, \; z=c$.\\
Hence,  ${\cal A}$ is a proper aBE algebra, verifying (pi) and not verifying (pimpl), which is linearly-ordered; hence it
is a {\bf proper pi-aBE algebra}.
We shall represent the set $A$ and the  binary relation $\leq$ by the following Hasse-type diagram:\\

\begin{figure}[htbp]
\begin{center} 
\begin{picture}(50,40)(0,10) 

\put(39,20){\makebox(2,2){$\circ$}}
\put(39,0){\makebox(2,2){$\circ$}}
 
\put(39,20){\line(0,4){23}}
\put(39,0){\line(0,4){23}}

                            \put(39,40){\makebox(2,2){$\circ$}}
                             \put(39,40){\line(0,4){25}}
                            \put(39,63){\makebox(2,2){$\circ$}}
                                 \put(39,73){\makebox(2,2){1}}

\put(45,20){\makebox(2,2){b}}
\put(45,0){\makebox(2,2){a}}

                               \put(45,40){\makebox(10,10){c}}

\end{picture}
 \end{center}
\caption{Proper linearly-ordered pi-aBE algebra   }
 \label{fig: A6}
\end{figure}
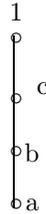
\end{ex}

\begin{ex}\em
Consider the set $A=\{0,a,b,c, 1\}$ with the following table of $\ra$: \\ 

\noindent  $\qquad \qquad \qquad \qquad$
\begin{tabular}{*{1}{c|}ccccc}
$\ra$ & 0 & a & b & c & 1 \\
\hline
0 & 1 & 1 & b & 1 & 1  \\
a & 0 & 1 & 1 & c & 1 \\
b & 0 & a & 1 & c & 1  \\
c & 0 & a & b & 1 & 1  \\
1 & 0 & a & b & c & 1
\end{tabular} 

\noindent\\
Then the algebra $ {\cal A}=(A, \ra, 1)$ verifies the properties (Re), (M), (L), (Ex) (hence (D)), (An) and  (pi).
It does not verify  (BB), (**) for  $x=b, \; y=0, \; z=a$, (B), (*), (Tr) for $0,a,b$,
(pimpl) for $x=0, \; y=a, \; z=b$.\\
  Hence,  ${\cal A}$ is a {\bf proper pi-aBE algebra}.

We shall represent the set $A$ and the  binary relation $\leq$ by the following Hasse-type diagram:\\

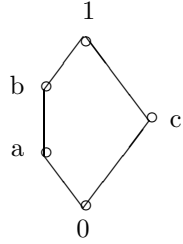
\begin{figure}[htbp]
\begin{center}
\begin{picture}(300,100)(-80,-20) 

\put(40,0){\makebox(2,2){$\circ$}}

\put(40,0){\line(3,4){25}}
\put(40,0){\line(-3,4){15}}

\put(65,33){\line(-3,4){24}}
\put(25,45){\line(3,4){15}}
\put(65,33){\makebox(2,2){$\circ$}}
                    \put(25,20){\line(0,1) {25}}
                      \put(25,20){\makebox(2,2){$\circ$}}
                     \put(25,45){\makebox(2,2){$\circ$}}

\put(40,62){\makebox(2,2){$\circ$}}

\put(35, -10){\makebox(10,5){0}}
\put(70,28){\makebox(10,10){c}}
\put(10,42){\makebox(10,10){b}}
\put(10,17){\makebox(10,10){a}}
\put(37,70){\makebox(10,10){1}}

\end{picture}
 \end{center}
\caption{Proper non-linearly-ordered  pi-aBE algebra } \label{fig: A3}
\end{figure}
\end{ex}

{\bf Remark:} We have proved that condition (pimpl) implies condition (pi); the  inverse implication is not true.
 Indeed, the two above aBE algebras   verify (pi), but 
do not verify (pimpl).\\

{\bf $\bullet$ Hilbert algebras}

\begin{ex}\em
Consider the set $A=\{a,b,c, 1\}$ with the following table of $\ra$: \\ 

\noindent  $\qquad \qquad \qquad \qquad$
\begin{tabular}{*{1}{c|}ccccc}
$\ra$ &  a & b & c & 1 \\
\hline
a & 1 & b & b & 1 \\
b & a & 1 & a & 1 \\
c & 1 & 1 & 1 & 1 \\
1 & a & b & c & 1 
\end{tabular}  

\noindent\\
Properties  (Re), (M), (L), (Ex) (hence (D)), (An),  (BB) (hence  (B), (**), (*), (Tr)) and (pimpl) (hence  (pi))
 are satisfied.\\
Hence, ${\cal A}$ is a {\bf pi-BCK = pimpl-BCK = Hilbert algebra}.
\end{ex}

                          \subsection{Examples of pi-RML algebras from Hierachy 2.1 H, 2.2 H and 2.3 H}

              \subsubsection{Examples of  pi-RML algebras without (Ex): \\
pi-tRML, pi-*RML, pi-RML**, pi-*RML**, pi-pre-BBBCC, pimpl-pre-BCC, \\pi-oRML, pi-*aRML, pi-aRML**, pi-*aRML**}
${}$


{\bf $\bullet$ Proper pi-tRML algebras} (P19-pi)

\begin{ex}\em
Consider the set $A=\{a, b,c,d, 1\}$ with the following table of $\ra$: \\ 

\noindent $\qquad \qquad \qquad \qquad$  
\begin{tabular}{*{1}{c|}ccccc}
$\ra$ & a & b & c & d & 1 \\
\hline
a & 1 & b  & b  & b  &  1 \\
b & a  & 1 & a  & d  &  1 \\
c &  a &  a & 1 &  1 & 1  \\
d & a  &  a &  1 & 1 &  1 \\
1 & a & b & c & d & 1
\end{tabular} 

\noindent\\
Properties (Re), (M),  (L), (Tr), (D) and (pi)  are satisfied. ${\cal A}$ does not satisfy  (An) for $c,d$;  
 (Ex) for $a,b,d$;
  (pimpl) for $a,b,d$;  (BB), (**) for $b,c,1$; (B) for $a,1,c$; (*) for $b,c,d$.\\
Hence, ${\cal A}$ is a {\bf proper pi-tRML algebra}, with (D).
\end{ex}

{\bf $\bullet$ Proper pi-*RML algebras} (P20-pi)

\begin{ex}\em
Consider the set $A=\{a, b,c,d, 1\}$ with the following table of $\ra$: \\ 

\noindent $\qquad \qquad \qquad \qquad$  
\begin{tabular}{*{1}{c|}ccccc}
$\ra$ & a & b & c & d & 1 \\
\hline
a & 1 & b  & b  & b  &  1 \\
b & a  & 1 & a  & a  &  1 \\
c &  a &  a & 1 &  1 & 1  \\
d & a  &  a &  1 & 1 &  1 \\
1 & a & b & c & d & 1
\end{tabular} 

\noindent\\
Properties (Re), (M), (L), (*) (hence (Tr)), (pi) are satisfied. ${\cal A}$ does not satisfy:  (Ex) for $a,c,b$;  (pimpl) for $b,c,a$;
(BB), (**) for $b,c,1$; (B) for $a,1,c$;  (D) for $b,c$.\\
Hence, ${\cal A}$ is a {\bf proper pi-*RML algebra}, without (D).
\end{ex}

{\bf $\bullet$ Proper pi-RML** algebras} (P21-pi)

\begin{ex}\em
Consider the set $A=\{a, b,c,d, 1\}$ with the following table of $\ra$: \\ 

\noindent  $\qquad \qquad \qquad \qquad$ 
\begin{tabular}{*{1}{c|}ccccc}
$\ra$ & a & b & c & d & 1 \\
\hline
a & 1 & b  & b  & b  &  1 \\
b & a  & 1 & a  & d  &  1 \\
c &  1 &  1 & 1 &  1 & 1  \\
d & 1  &  1 &  1 & 1 &  1 \\
1 & a & b & c & d & 1
\end{tabular} 

\noindent\\
Properties (Re), (M), (L),  (**) (hence (Tr)), (pi) are satisfied. 
${\cal A}$ does not satisfy: (An) for $c,d$; (Ex), (pimpl)  for $a,b,d$; 
(BB) for $d,b,a$; (B), (*) for $b,c,d$; (D) for $d,a$.\\
Hence, ${\cal A}$ is a {\bf proper pi-RML** algebra}, without (D).
\end{ex}   

{\bf $\bullet$  Proper pi-*RML** algebras} (P22-pi)

There are no proper pi-*RML** algebras with three,
 four or  five elements  (we have run
PASCAL programs searching for such an algebra). Finding  an example of pi-*RML** algebra with  minimum six 
elements, if it exists,  was announced as an open problem in our preprint on arXiv.
 Professor Michael Kinyon, Department of Mathematics, University of Denver,  communicated us the 
following example of proper pi-*RML** algebra  with six elements, found by using the finite model finder Mace4:
 
\begin{ex}\em (Michael Kinyon)

Consider the set $A=\{a, b,c,d,e, 1\}$ with the following table of $\ra$: \\ 

\noindent  $\qquad \qquad \qquad \qquad$ 
\begin{tabular}{*{1}{c|}cccccc}
$\ra$ & a & b & c & d & e & 1 \\
\hline
a & 1 & b  & c  & e  & e & 1 \\
b & a  & 1 & c  & a  &  a &1 \\
c &  1 &  1 & 1 &  1 & 1 & 1  \\
d & 1  &  b &  b & 1 &  1 & 1\\
e & 1  &  b &  b & 1 &  1 & 1\\
1 & a & b & c & d & e &1
\end{tabular} 

\noindent\\
Properties (Re), (M), (L), (*), (**) and (pi) are satisfied. 
${\cal A}$ does not satisfy: (An) for $d,e$; (Ex), (pimpl)  for $a,b,d$ (note that (Ex) and (pimpl) are not satisfied 
for the same elements); (BB) for $c,a,d$;  (B) for $b,d,c$; (D) for $c,d$.\\
Hence, ${\cal A}$ is a {\bf proper pi-*RML** algebra}.
\end{ex}

{\bf $\bullet$ Proper pi-pre-BBBCC  (P23-pi) and pimpl-pre-BBBCC (P23-pimpl) algebras}

\begin{ex}\em
Consider the set $A=\{a, b,c,d, 1\}$ with the following table of $\ra$: \\ 

\noindent  $\qquad \qquad \qquad \qquad$ 
\begin{tabular}{*{1}{c|}ccccc}
$\ra$ & a & b & c & d & 1 \\
\hline
a & 1 & b  & b  & d  &  1 \\
b & a  & 1 & 1  & d  &  1 \\
c &  a &  1 & 1 &  d & 1  \\
d & a  &  c &  c & 1 &  1 \\
1 & a & b & c & d & 1
\end{tabular} 

\noindent\\
Properties (Re), (M), (L), (BB) (hence (D)), (B) (hence  (**), (*), (Tr)), (pi) are satisfied. 
${\cal A}$ does not satisfy: (An) for $b,c$; (Ex), (pimpl)  for $a,d,b$.\\
Hence, ${\cal A}$ is a {\bf proper pi-pre-BBBCC algebra}.
\end{ex}

\begin{ex}\em
Consider the set $A=\{a, b,c,d, 1\}$ with the following table of $\ra$: \\ 

\noindent $\qquad \qquad \qquad \qquad$   
\begin{tabular}{*{1}{c|}ccccc}
$\ra$ & a & b & c & d & 1 \\
\hline
a & 1 & b  & b  & 1  &  1 \\
b & a  & 1 & 1  & a  &  1 \\
c &  a &  1 & 1 &  a & 1  \\
d & 1  &  c &  c & 1 &  1 \\
1 & a & b & c & d & 1
\end{tabular} 

\noindent\\
Properties (Re), (M), (L), (BB) (hence (D)), (B) (hence  (**), (*), (Tr)) and  (pimpl)  (hence (pi)) are satisfied. 
${\cal A}$ does not satisfy: (An) for $b,c$; (Ex) for $a,d,b$.\\
Hence, ${\cal A}$ is a {\bf proper pimpl-pre-BBBCC algebra}.
\end{ex}


{\bf $\bullet$ Proper pi-oRML algebras} (P24-pi)

\begin{ex}\em
Consider the set $A=\{a, b,c,d, 1\}$ with the following table of $\ra$: \\ 

\noindent $\qquad \qquad \qquad \qquad$  
\begin{tabular}{*{1}{c|}ccccc}
$\ra$ & a & b & c & d & 1 \\
\hline
a & 1 & b  & b  & b  &  1 \\
b & a  & 1 & a  & a  &  1 \\
c &  a &  a & 1 &  d & 1  \\
d & a  &  1 &  a & 1 &  1 \\
1 & a & b & c & d & 1
\end{tabular} 

\noindent\\
Properties (Re), (M), (L), (An),  (Tr) and  (pi)  are satisfied. 
${\cal A}$ does not satisfy: (Ex) for $a,c,b$; (pimpl) for $a,c,d$; (BB) for $b,c,d$; (**) for $b,c,1$; 
(B) for $a,1,c$; (*) for $c,d,b$; (D) for $b,c$.\\
Hence, ${\cal A}$ is a {\bf proper pi-oRML algebra}, without (D).
\end{ex}

{\bf $\bullet$ Proper pi-*aRML algebras} (P25-pi)

\begin{ex}\em
Consider the set $A=\{ a,b,c,d, 1\}$ with the following table of $\ra$: \\ 

\noindent  $\qquad \qquad \qquad \qquad$ 
\begin{tabular}{*{1}{c|}ccccc}
$\ra$ & a & b & c & d & 1 \\
\hline
a & 1 & b  & b  & b  &  1 \\
b & a  & 1 & a  & a  &  1 \\
c &  a &  a & 1 &  a & 1  \\
d & a  &  a &  a & 1 &  1 \\
1 & a & b & c & d & 1
\end{tabular} 

\noindent\\
Then the algebra $ {\cal A}=(A, \ra, 1)$ verifies properties (Re), (M), (L), (An), (*) (hence (Tr)) and (pi).
 It does not verify: (Ex) for $a,c,b$, (BB) for $b,c,1$, (B) for $a,1,c$,  (**) for $b,c,1$,  (D) for $b,c$, 
(pimpl) for $b,c,a$.\\
  Hence, ${\cal A}$ is a {\bf proper pi-*aRML algebra},  without (D).
\end{ex}

{\bf $\bullet$ Proper pi-aRML** algebras} (P26-pi)

\begin{ex}\em
Consider the set $A=\{ a,b,c,d, 1\}$ with the following table of $\ra$: \\ 

\noindent  $\qquad \qquad \qquad \qquad$ 
\begin{tabular}{*{1}{c|}ccccc}
$\ra$ & a & b & c & d & 1 \\
\hline
a & 1 & b  & b  & b  &  1 \\
b & a  & 1 & a  & d  &  1 \\
c &  1 &  1 & 1 &  1 & 1  \\
d & a  &  1 &  a & 1 &  1 \\
1 & a & b & c & d & 1
\end{tabular} 

\noindent\\
Then the algebra $ {\cal A}=(A, \ra, 1)$ verifies properties (Re), (M), (L), (An), (**) (hence (Tr)) and (pi).
 It does not verify: (Ex) for $a,b,d$, (BB) for $d,b,a$, (B) for $b,c,d$,  (*) for $b,c,d$, (pimpl) for $a,b,d$,
(D) for $d,a$.\\
  Hence, ${\cal A}$ is a {\bf proper pi-aRML** algebra}, without (D).
\end{ex}

{\bf $\bullet$ Proper pi-*aRML** algebras} (P27-pi)

\begin{ex}\em
Consider the set $A=\{ a,b,c,d, 1\}$ with the following table of $\ra$: \\ 

\noindent  $\qquad \qquad \qquad \qquad$ 
\begin{tabular}{*{1}{c|}ccccc}
$\ra$ & a & b & c & d & 1 \\
\hline
a & 1 & b  & b  & d  &  1 \\
b & a  & 1 & c  & d  &  1 \\
c &  a &  1 & 1 &  a & 1  \\
d & 1  &  b &  b & 1 &  1 \\
1 & a & b & c & d & 1
\end{tabular} 

\noindent\\
Then the algebra $ {\cal A}=(A, \ra, 1)$ verifies properties (Re), (M), (L), (An), (*), (**) (hence (Tr)),
  and (pi).
 It does not verify: (Ex) for $a,b,c$, (BB) for $c,b,a$, (B) for $a,c,d$, (D) for $c,a$, 
and it does not verify condition (pimpl) for $a,b,c$.\\
Hence, ${\cal A}$ is a  {\bf proper pi-*aRML** algebra},
  without (D).
\end{ex}

               \subsubsection{Examples of pi-RML algebras with  (Ex):\\ pi-BE**, pi-aBE** algebras}
${}$

{\bf $\bullet$ Proper pi-BE** algebras} (P30-pi)

\begin{ex}\em
Consider the set $A=\{a, b,c,d, 1\}$ with the following table of $\ra$: \\ 

\noindent  $\qquad \qquad \qquad \qquad$ 
\begin{tabular}{*{1}{c|}ccccc}
$\ra$ & a & b & c & d & 1 \\
\hline
a & 1 & b  & b  & d  &  1 \\
b & a  & 1 & a  & d  &  1 \\
c &  1 &  1 & 1 &  1 & 1  \\
d & 1  &  1 &  1 & 1 &  1 \\
1 & a & b & c & d & 1
\end{tabular} 

\noindent\\
Properties (Re), (M), (L), (Ex) (hence (D)), (**) (hence  (Tr)) and (pi) are satisfied. 
${\cal A}$ does not satisfy:  (An) for $c,d$,  
(BB)  for $d, a, c$,   (B) for $a, c , d$, (*) for $a, c, d$, (pimpl) for $a, c, d$.\\
Hence, ${\cal A}= (A, \ra,1)$  is a {\bf proper pi-BE** algebra}.
\end{ex}

{\bf  $\bullet$ Proper pi-aBE** algebras} (P31-pi)

\begin{ex}\em 
Consider the set $A=\{ a,b,c,d, 1\}$ with the following table of $\ra$: \\ 

\noindent  $\qquad \qquad \qquad \qquad$ 
\begin{tabular}{*{1}{c|}ccccc}
$\ra$ & a & b & c & d & 1 \\
\hline
a & 1 & b  & b  & d  &  1 \\
b & a  & 1 & a  & d  &  1 \\
c &  1 &  1 & 1 &  1 & 1  \\
d & a  &  b &  c & 1 &  1 \\
1 & a & b & c & d & 1
\end{tabular} 
 
\noindent\\
Then the algebra $ {\cal A}=(A, \ra, 1)$ verifies properties (Re), (M), (L), (Ex), (An), (**) (hence (Tr)) and (pi).
 It does not verify  condition (BB)  for $x=d$, $y=a$, $z=c$, condition (B) for $a,c,d$, condition (*) for $a,c,d$;
 it does not verify (pimpl) for $x=a$, $y=c$, $z=d$.
 The relation $\leq$
is reflexive, antisymmetrique and tranzitive, hence is an order relation, namely is a lattice order.\\
 Hence,   ${\cal A}$ is a {\bf proper pi-aBE** lattice}.
\end{ex}

           \section{Final remarks}

\begin{rems}\em ${}$

(i) There are only $2^1 = 2$  RM algebras $(A=\{a,1\}, \ra, 1)$ with two elements: \\

\noindent $\qquad \qquad \qquad \qquad$ 
\begin{tabular}{*{1}{c|}cc}
$\ra$ & a &  1   \\
\hline
a & 1 & $\cdot$  \\
1 & a &  1  
\end{tabular} 

\noindent\\ 
One is a special case of  BCI algebra, namely the {\bf p-semisimple BCI algebra with two elements}:\\

\noindent  $\qquad \qquad \qquad \qquad$
\begin{tabular}{*{1}{c|}cc}
$\ra$ & a &  1   \\
\hline
a & 1 & a  \\
1 & a &  1  
\end{tabular} 

\noindent\\ 
  the other is a special case of BCK algebra (hence of RML algebra), 
namely the particular case of Hilbert algebra which is the {\bf Boolean algebra with two elements}:\\

\noindent  $\qquad \qquad \qquad \qquad$
\begin{tabular}{*{1}{c|}cc}
$\ra$ & a &  1   \\
\hline
a & 1 & 1  \\
1 & a &  1  
\end{tabular} 

\noindent\\

(ii) There are $3^4 = 81$ RM algebras $(A=\{a,b,1\}, \ra,1)$ with three elements:\\

\noindent  $\qquad \qquad \qquad \qquad$
\begin{tabular}{*{1}{c|}ccc}
$\ra$ & a & b& 1  \\
\hline
a& 1& $\cdot$ & $\cdot$   \\
b& $\cdot$ & 1& $\cdot$ \\
1& a& b& 1 
\end{tabular} 

\noindent\\
and among them there are $3^2=9$ RML algebras with three elements:\\

\noindent  $\qquad \qquad \qquad \qquad$
\begin{tabular}{*{1}{c|}ccc}
$\ra$ & a & b& 1  \\
\hline
a& 1& $\cdot$ & 1   \\
b& $\cdot$ & 1& 1 \\
1& a& b& 1 
\end{tabular} 

A Pascal programme has determined all the 81 RM algebras. They are: \\
- 4 proper RM algebras (2 with (D), 2 without (D)),  2 proper pre-BBBZ algebras,
 2 proper pre-BCI algebras (with (D));\\
- 3 proper BCI algebras (1 with (D), 2 without (D)), 8 proper aRM algebras (4 with (D), 4 without (D)), 
8 proper *aRM algebras (2 with (D), 6 without (D)), 24 proper oRM algebras (8 with (D), 16 without (D)),
4 proper aRM** algebras (all without (D)), 17 proper *aRM** algebras (1 with (D), 16 without (D));\\
- (the 9 RML algebras:) 1 proper pimpl-pre-BCK algebra, 3 proper *aRML algebras (all with (D)), 
2 proper BCK algebras, 3 Hilbert algebras.

(iii) There are $4^9 = 262.144$ RM algebras $(A=\{a,b,c,1\}, \ra, 1)$  with four elements:\\

\noindent  $\qquad \qquad \qquad \qquad$
\begin{tabular}{*{1}{c|}cccc}
$\ra$ &  a & b & c & 1 \\
\hline
a & 1 & $\cdot$ & $\cdot$ & $\cdot$ \\
b & $\cdot$ & 1 & $\cdot$ & $\cdot$ \\
c & $\cdot$ & $\cdot$ & 1 & $\cdot$ \\
1 & a & b & c & 1 
\end{tabular}  

\noindent\\
and among them there are $4^6=4.096$ RML algebras with four elements:\\

\noindent  $\qquad \qquad \qquad \qquad$
\begin{tabular}{*{1}{c|}cccc}
$\ra$ &  a & b & c & 1 \\
\hline
a & 1 & $\cdot$ & $\cdot$ & 1 \\
b & $\cdot$ & 1 & $\cdot$ & 1 \\
c & $\cdot$ & $\cdot$ & 1 & 1 \\
1 & a & b & c & 1 
\end{tabular}  

\noindent\\

Different Pascal programmes can provide the algebras we look for. We can say that there are many proper pre-BBBZ algebras,
but there are no proper pre-BZ algebras, proper pre-BCC algebras, proper pre-BBBCC algebras with four elements.

(iv) There are $5^{16} = 152.587.890.625$ RM algebras $(A=\{a,b,c,d,1\}, \ra, 1)$  with five elements:\\

\noindent  $\qquad \qquad \qquad \qquad$ 
\begin{tabular}{*{1}{c|}ccccc}
$\ra$ & a & b & c & d & 1 \\
\hline
a & 1 & $\cdot$  & $\cdot$  & $\cdot$  &  $\cdot$ \\
b & $\cdot$ & 1 & $\cdot$  & $\cdot$  &  $\cdot$ \\
c &  $\cdot$ &  $\cdot$ & 1 &  $\cdot$ & $\cdot$  \\
d & $\cdot$  &  $\cdot$ &  $\cdot$ & 1 &  $\cdot$ \\
1 & a & b & c & d & 1
\end{tabular} 
 
\noindent\\
and among them there are $5^{12}=244.140.625$ RML algebras with five elements.

A Pascal programme has found, for example, all the  proper pimpl-pre-BBBCC algebras with five elements, 
which are in number of 60. 

(v) There are $6^{20}=3.656.158.440.062.976$ RML algebras  $(A=\{a,b,c,d,e,1\}, \ra, 1)$ with six elements :\\

\noindent $\qquad \qquad \qquad \qquad$  
\begin{tabular}{*{1}{c|}cccccc}
$\ra$ & a & b & c & d &e & 1 \\
\hline
a & 1 & $\cdot$  & $\cdot$  & $\cdot$ & $\cdot$ &  1 \\
b & $\cdot$  & 1 & $\cdot$  & $\cdot$  &  $\cdot$ & 1 \\
c & $\cdot$ &  $\cdot$ & 1 & $\cdot$ & $\cdot$ & 1  \\
d & $\cdot$ & $\cdot$ & $\cdot$ & 1  & $\cdot$ &  1 \\
e & $\cdot$ & $\cdot$ & $\cdot$ & $\cdot$ & 1  &  1  \\
1 & a & b & c & d & e & 1
\end{tabular} 

\noindent\\
Hence, it is  difficult to look for an example of RML algebra with six elements. 
\end{rems}

\begin{rems}\em ${}$

1). The properties (*) and (**) are independent, i.e. there are algebras verifying (*) and not verifying (**) (for example,
the proper *RM algebras)
and there algebras verifying (**) and not verifying (*) (for example, the proper RM** algebras). 

2). Properties  (BB) and (B) are dependent, namely  (BB) $\Rightarrow$ (B). We have examples of  RM algebras verifying (B)
 and not verifying (BB) (namely the proper pre-BZ, BZ, pre-BCC, BCC algebras), we have examples of RM algebras verifying
both (BB) and (B) (the proper pre-BBBZ algebras, the proper pre-BBBCC algebras).

3). Properties  (B) and (**) are dependent, namely (B) $\Rightarrow$ (**). 
We have examples of RM algebras verifying (**) and not verifying (B)
(for example, the proper RM** algebras),
we have examples of  RM algebras verifying both (B) and (**) 
(the proper  pre-BZ, BZ, pre-BCC, BCC algebras).
 
4).  We found very  many examples of pre-BBBZ algebras  and  few examples of  pre-BBBCC algebras. 
We have (see Figure \ref{fig:fig23sT}):  {\bf pre-BBBCC} = {\bf pre-BBBZ} + (L);
 {\bf pre-BBBZ} + (An) = {\bf BCI} and  {\bf pre-BBBCC} + (An) = {\bf BCK}.
\end{rems}

\end{document}